\documentclass[11pt]{article}

\usepackage[a4paper]{geometry}
\usepackage{mathrsfs}
\usepackage[all]{xy}
\usepackage[T1]{fontenc}
\usepackage{amsmath}
\usepackage{amsthm}
\usepackage{amssymb}
\usepackage{array}
\usepackage{ifthen}
\usepackage{pb-diagram}
\usepackage{tikz}
\usepackage{verbatim}
\usepackage{amscd}

\pdfoutput=1

\swapnumbers
\theoremstyle{definition}
\newtheorem{absatz}{}[section]
\newtheorem{absatz4}{}[subsection]
\theoremstyle{theorem}
\newtheorem{cor}[absatz]{Corollary}
\newtheorem{cor4}[absatz4]{Corollary}
\newtheorem{thm}[absatz]{Theorem}

\newtheorem{prop}[absatz]{Proposition}
\newtheorem{prop4}[absatz4]{Proposition}

\newtheorem{lem}[absatz]{Lemma}
\newtheorem{lem4}[absatz4]{Lemma}
\newtheorem{defn4}[absatz4]{Definition}
\theoremstyle{definition}
\theoremstyle{remark}
\newtheorem{rem}[absatz]{Remark}

\theoremstyle{definition}

\theoremstyle{theorem}
\newtheorem*{theorem*}{Theorem}

\def\CCC{{\mathbb C}}

\def\HHH{{\mathbb H}}

\def\NNN{{\mathbb N}}

\def\PPP{{\mathbb P}}
\def\QQQ{{\mathbb Q}}
\def\RRR{{\mathbb R}}

\def\ZZZ{{\mathbb Z}}

\def\({\left(}
\def\){\right)}

\def\<{\left\langle}
\def\>{\right\rangle}

\DeclareMathOperator{\tr}{tr}
\DeclareMathOperator{\id}{id}

\DeclareMathOperator{\im}{Im}
\DeclareMathOperator{\re}{Re}

\DeclareMathOperator{\Spec}{Spec}

\DeclareMathOperator{\res}{res}

\DeclareMathOperator{\vol}{vol}
\DeclareMathOperator{\discr}{disc}
 
\DeclareMathOperator{\fin}{fin}

\DeclareMathOperator{\jac}{{\rm{Jac}}}

\newcommand{\hyp}{\mu_{\rm{hyp}}}
\newcommand{\can}{\mu_{\rm{can}}}

\newcommand{\dif}{{\rm{d}}}

\begin{document}

\title{Self-intersection of the relative dualizing sheaf on modular curves $X_1(N)$}
\author{Hartwig Mayer}
\date{\today}

\maketitle

\abstract{
Let $N$ be an odd and squarefree positive integer divisible by at least two relative prime integers bigger or equal than $4$.
Our main theorem is  an asymptotic formula solely in terms of $N$ for the stable arithmetic self-intersection number of the relative dualizing sheaf  for  modular curves $X_1(N)/\QQQ$.
 From our main theorem we obtain an asymptotic formula for the stable Faltings height of the Jacobian $J_1(N)/\QQQ$  of $X_1(N)/\QQQ$, 
and, for sufficiently large $N$,  an effective version of  Bogomolov's conjecture for  $X_1(N)/\QQQ$.}

\bibliographystyle{amsalpha}

\section{Introduction}

Let $K$ be a number field and   $\mathcal O_K$ its ring of integers.  Let $\mathcal X /\mathcal O_K$ be the minimal regular model of a smooth projective curve $X/K$ of genus $g_X>0$. We call 
 $\mathcal X /\mathcal O_K$  an arithmetic surface.
In \cite{arakelov}, S. J. Arakelov introduced an intersection theory for metrized invertible sheaves on $\mathcal X /\mathcal O_K$. G. Faltings established in his work \cite{faltings} many fundamental results 
 to the theory of Arakelov.
Within this framework we may attach two important invariants  to   the curve $X/K$: 
Let $\overline{\omega}_{\mathcal X / \mathcal O_K}$ be
 the relative dualizing sheaf on $\mathcal X /\mathcal O_K$ equipped with the Arakelov metric. The first invariant is the stable arithmetic self-intersection number
$\displaystyle  \frac{1}{\lbrack K: \QQQ \rbrack} \overline{\omega}_{\mathcal X/ \mathcal O_K}^2$  which is independent of the field $K$ as long as  $X/K$ has semistable reduction over $\mathcal O_K$.
The second invariant is the arithmetic degree of the direct image of $\overline{\omega}_{\mathcal X / \mathcal O_K}$ which is, in other words, the stable Faltings height of 
the Jacobian ${\mathrm{Jac}}(X)/K$   of $X/K$.

The arithmetic significance of the stable arithmetic self-intersection number was given in \cite{szpiro} where it was shown that its strict positivity is equivalent to Bogomolov's conjecture (finally proven by E. Ullmo after
partial results by Burnol, Szpiro, and Zhang). Recall that this conjecture claims
that the set of algebraic points of the curve $X/K$ embedded into its Jacobian  are discretely distributed with respect to the ``N$\acute{\mathrm{e}}$ron-Tate topology'' supposing that the genus of the curve is bigger one. The second
invariant is particularly interesting in the situation of modular curves. E.g.,  the stable Faltings height of the Jacobian of the modular curve $X_1(N)/\QQQ$ plays an important role in the recent 
work \cite{edixhoven} on computing Fourier coefficients of modular forms.

The only cases so far in which the stable arithmetic self-intersection number of the relative dualizing sheaf is  known are arithmetic surfaces whose generic fiber is a curve of genus one (\cite{faltings}),
a curve of genus two (\cite{bostmestremoret-bailly}), or a modular curve $X_0(N)$, where $N$ is squarefree and $2,3 \nmid N$, (\cite{abbesullmo}, \cite{michelullmo}). More recently, upper bounds in the cases
of modular curves $X(N)$ and Fermat curves were found (see \cite{curillakuehn}). 
 In the case of the modular curves
$X_0(N)$, the stable Faltings height is asymptotically determined in \cite{jorgensonkramer1}. Asymptotics in the case of the modular curve $X_1(N)$ are already given in \cite{edixhoven} (cf. remark \ref{rem:edixhoven}).

\begin{absatz} {\it{Arakelov theory on arithmetic surfaces}}.
 Let $\mathcal X /\mathcal O_K$ be the minimal regular model of a smooth projective curve $X/K$ of genus $g_X>0$.
Let $D$ be a divisor on $\mathcal X$ and $\mathcal L = \mathcal O_{\mathcal X}(D)$ the corresponding line bundle on $\mathcal X /\mathcal O_K$. For every embedding
$\sigma: K \longrightarrow \CCC$ we equip the induced line bundle $L_{\sigma}$ on the compact Riemann surface $X_{\sigma}(\CCC)$, $X_{\sigma}:= X \times_{\sigma} {\mathrm{Spec}}(\CCC)$, with the unique admissible metric
(Arakelov metric) with respect to the canonical volume form $\can$ (cf. \cite{soule}, p. 332).
A line bundle $\mathcal L$ equipped with these  metrics for all embeddings $\sigma$ will be denoted by $\overline{\mathcal L}$.
For  the relative dualizing sheaf $\overline{\omega}_{\mathcal X / \mathcal O_K}$  the Arakelov metric has the following interpretation:
the residual maps
 \begin{align*}
 \Omega_{X_{\sigma}(\CCC)}^1 \otimes \mathcal O_{X_{\sigma}(\CCC)}(P)|_P \longrightarrow \CCC 
 \end{align*}
are isometries for all points $P \in X_{\sigma}(\CCC)$ and all embeddings $\sigma: K \longrightarrow \CCC$, where $\mathcal O_{X_{\sigma}(\CCC)}(P)$ is equipped with the Arakelov metric and $\CCC$ with the standard hermitian metric (see \cite{soule}, p. 333).
\newline
Let $\overline{\mathcal L} = \overline{ \mathcal O}_{\mathcal X}(P) $ and $\overline{\mathcal M}=  \overline{\mathcal O}_{\mathcal X}(Q)$ be  two  metrized line bundles, where $P,Q$ are two horizontal prime divisors on $\mathcal X$
with no common component. The intersection product of $\overline{\mathcal L}$ and $\overline{\mathcal M}$  is given by
\begin{align}
\overline{\mathcal L} \cdot \overline{\mathcal M} = (P,Q)_{\fin} - \sum_{\sigma : K \rightarrow \CCC} g_{\rm{can}}^{\sigma}\big(P_{\sigma}, Q_{\sigma}),
\end{align}
where $ (P,Q)_{\fin}$ is their local intersection product on $\mathcal X / \mathcal O_K$ (see \cite{soule}, p. 332),  $P_{\sigma}, Q_{\sigma}$ the induced points on $X_{\sigma}(\CCC)$,
 and  $g^{\sigma}_{\rm{can}}$ is the canonical Green's function  on 
$X_{\sigma}(\CCC) \times X_{\sigma}(\CCC) \setminus \Delta_{ X_{\sigma}(\CCC)}$ (for the definition see section \ref{absatz:Green}).

\end{absatz}

\begin{absatz} {\it Main results}.
Let $X_1(N)/ \QQQ$  be the smooth projective algebraic curve over $\mathbb{Q}$ that classifies elliptic curves equipped with a point of exact order $N$.
Let $N$ be  an odd and squarefree integer  of the form $N=N^{'} qr>0$ with $q$ and $r$ relative prime integers satisfying  $ q,r \geq 4$. Then, the minimal regular model
$\mathcal X_1(N) / \ZZZ \lbrack \zeta_N \rbrack$ is semistable (cf. proposition \ref{prop:katzmazur}) and the genus of $X_1(N)/\QQQ$ denoted by $g_N$ satisfies $g_N >0$.
With the notation $\overline{\omega}_N^2 = \displaystyle  \frac{1}{\lbrack K: \QQQ \rbrack} \overline{\omega}_{\mathcal X_1(N) / \mathcal O_K}^2$ our main theorem (cf. theorem \ref{thm:main}) is the following:
\begin{theorem*}
 Let $N$ be  an odd and squarefree integer  of the form $N=N^{'} qr>0$ with $q$ and $r$ relative prime integers satisfying  $ q,r \geq 4$. Then, we have
\begin{align*} 
\overline{\omega}_N^2 = 3 g_N \log(N) + o\big(g_N \log(N) \big).
\end{align*}
\end{theorem*}
Our first arithmetic application is the following asymptotics
 for the stable Faltings height $h_{\rm{Fal}}\big(J_1(N) \big)$
 of the Jacobian $J_1(N)/\QQQ$ of $X_1(N)/ \QQQ$ (cf. theorem \ref{thm:falheight}).
\begin{theorem*}
 Let $N$ be  an odd and squarefree integer  of the form $N=N^{'} qr>0$ with $q$ and $r$ relative prime integers satisfying  $ q,r \geq 4$. Then, we have
\begin{align*}
h_{\rm{Fal}}\big(J_1(N) \big) = \frac{g_{N}}{4} \log (N) + o\big(g_{N} \log (N) \big).
\end{align*}
\end{theorem*}
We also obtain an asymptotic formula   for the admissible self-intersection number of the relative dualizing sheaf $\overline{\omega}_{a,N}^2$ in the sense of the theory of Zhang in \cite{zhang}. From this
we can deduce, for large $N$,  an effective version of  Bogomolov's conjecture of the modular curve $X_1(N)/\QQQ$ (cf. theorem \ref{thm:bogomolov}).
\end{absatz}
\begin{absatz} {\it Outline}. The main structure of this article comes from the later section 7, namely, from proposition \ref{prop:geometricpart},  providing us the formula
\begin{align}
\overline{\omega}^{2}_N =   4 g_{N} (g_N -1)  g_{{\rm{can}}}(0, \infty) + 
\frac{1}{\varphi(N)} \frac{g_N+1}{g_N - 1}  \left( V_{0}, V_{\infty} \right)_{\rm{fin}},
\end{align}
for some explicit vertical divisors $V_0, V_{\infty}$ on $\mathcal X_1(N) / \ZZZ \lbrack \zeta_N \rbrack$. Hence, to achieve our main theorem, we first determine 
 the analytic part $4 g_{N} (g_N -1)  g_{{\rm{can}}}(0, \infty)$,  and then
 compute the algebraic part $\frac{1}{\varphi(N)} \frac{g_N+1}{g_N - 1}  \left( V_{0}, V_{\infty} \right)_{\rm{fin}}$ of the stable arithmetic self-intersection number $\overline{\omega}^{2}_N$. 
 \newline
For the analytic part we follow the strategy of \cite{abbesullmo}.
 In the second section, after recalling some basic facts of the (compactified) modular curves, we present the spectral
expansion of the automorphic kernels of weight $0$ and $2$. We conclude this section with a basic formula expressing the arithmetic average 
\begin{align*}
F(z):=\frac{1}{g_{N}} \,  \,\sum_{j=1}^{g_{N}} y^2  |f_j(z)|^2 \hspace{1,5cm} (z=x + iy \in \HHH)
\end{align*}
of an orthonormal basis $\{ f_j \}_{j=1}^{g_N}$ of holomorphic cusp forms of weight $2$ with respect to the congruence subgroup $\Gamma_1(N)$ in terms of the hyperbolic,
parabolic, and spectral parts of the expansions of the above automorphic kernels.  In section three we obtain a formula for 
$ g_{{\rm{can}}}(0, \infty)$, based on previous work of A. Abbes and E. Ullmo as well as J. Jorgenson and J. Kramer, essentially in terms of the constant term in the Laurent
expansion at $s=1$ of the Rankin-Selberg transform $R_F(s)$, $ s \in \CCC$, of the function $F(z)$ which we denote by  $C_F$. 
In order to compute the constant $C_F$ we exploit our basic formula and determine in section four the
 contribution of the  Rankin-Selberg transform of the  hyperbolic and in section five the
 contribution of the  Rankin-Selberg transform of the  parabolic and spectral
 parts to $C_F$  (cf. remark \ref{rem:strategy}). In section 6 we are then ready to determine the analytic part, and after computing the algebraic part of the stable 
arithmetic self-intersection number $\overline{\omega}^{2}_N$, we obtain in section 7 our main theorem. In section 8 we deduce the arithmetic applications mentioned above. Finally,
 we study in the appendix an Epstein zeta function naturally arising in section 4.
\end{absatz}


\begin{absatz}{\it{Acknowledgement.}}
The results of this article are essentially those of my thesis \cite{mayer}. I am much indebted to my advisor J\"urg Kramer for his support and very valuable concrete comments on this work. Without his help
this article would not have appeared.  I want to thank Bas Edixhoven who
pointed out to me two mistakes in a preliminary version of this article. To 
Anna von Pippich  I am very grateful  for discussions concerning the analytic part of this article.
\end{absatz}

\section{Background Material}
Let us collect some basic material of modular curves $X\big(\Gamma_1(N) \big)$ and their spectral theory. Our main references are \cite{diamondshurman}, \cite{iwaniec}, and \cite{roelcke,roelcke2}. The spectral theory in
\cite{iwaniec} and \cite{roelcke, roelcke2} is formulated in a more general framework, namely for Fuchsian subgroups of the first kind. In order to keep this exposition short, we restrict the discussion to the congruence
subgroup $\Gamma_1(N)$.  
 \begin{absatz} {\it{The upper half-plane}}. 
Let $\HHH := \{z= x+iy \in \CCC \, | \, \im(z)=y >0 \}$ be the upper half-plane equipped with the hyperbolic metric
\begin{align} \label{align:hypmetric} 
\dif s_{\mathrm{hyp}}^2(z):= \frac{\dif x^2 + \dif y^2 }{y^2}
\end{align}
giving $\HHH$ the structure of a $2$-dimensional Riemannian manifold of constant negative curvature equal to $-1$. The hyperbolic metric \eqref{align:hypmetric} induces the distance function
$\rho$ on $\HHH$ defined  by
\begin{align*}
\cosh \big( \rho(z,w) \big) = 1 + 2 u(z,w),
\end{align*}
 where 
 \begin{align} \label{align:hypdistance}
 u(z,w)= \frac{\mid z-w \mid^2}{4 \im (z) \im (w)} \hspace{1cm} (z,w \in \HHH),
 \end{align}
and the hyperbolic volume form on $\HHH$ given by
\begin{align} \label{align:hypvolume}
\hyp(z) := \frac{\dif x  \wedge \dif y}{y^2}.
\end{align}
\end{absatz} 

 \begin{absatz}{\it{Modular curves $X\big(\Gamma_1(N) \big)$.}}
     Let $N \geq 1$ be a positive integer and $\Gamma_1(N) \subseteq {\rm{SL}}_2(\ZZZ)$  the congruence subgroup defined by
\begin{align*}
\Gamma_1(N) := \left \{ \gamma=  \begin{pmatrix} a & b \\ c & d \end{pmatrix} \in  {\rm{SL}}_2(\ZZZ) \mid a \equiv d \equiv 1 \, \mbox{ mod } N, \,  c \equiv 0 \, \mbox{ mod } N \right \}
\end{align*}
acting by fractional linear transformation $z \mapsto \gamma z := \tfrac{az +b}{cz +d}$, $\gamma= \left( \begin{smallmatrix} a & b \\ c & d \end{smallmatrix} \right) \in  \Gamma_1(N)$,
 on the upper half-plane $\HHH$. The quotient space $\Gamma_1(N) \backslash \HHH$ is denoted by $Y\big( \Gamma_1(N) \big)$.
  Letting $\Gamma_1(N)$ act on the projective line $\PPP^1_{\QQQ}$ by $\gamma (s:t) : = (as + bt : cs + dt)$, $(s:t) \in  \PPP^1_{\QQQ}$,
 the (compactified) modular curve $X\big(\Gamma_1(N) \big)$  associated to  $\Gamma_1(N)$ is defined as   the quotient space 
 \begin{align*}
X\big(\Gamma_1(N) \big):= \Gamma_1(N) \backslash \big( \HHH \cup \PPP^1_{\QQQ} \big),
 \end{align*}
which can be endowed with a natural
 topology making the quotient space into a compact Riemann surface (see \cite{diamondshurman}, chap. 2.4). 
\newline
The finite set  $P_{\Gamma_1(N)}:=X\big( \Gamma_1(N) \big) \backslash Y \big(\Gamma_1(N) \big)$ is called the set of (inequivalent) cusps  of $\Gamma_1(N)$ and is represented by
 the set of elements   $(a,c)$ of  $(\ZZZ/ N \ZZZ)^2 $ of order $N$ modulo the equivalence relation $(a:c) \equiv (a',c')$ if and only if $(a',c') = (a + n c, c)$ for some $n \in \ZZZ/N \ZZZ$. 
 Moreover, the hyperbolic volume form $\hyp$ of \eqref{align:hypvolume} descends to a volume form on $X\big(\Gamma_1(N) \big)$  (see 
 \cite{diamondshurman}, p. 181) which we will denote again by $\hyp$. For the moduli interpretation of the modular curve $X\big(\Gamma_1(N) \big)$ we refer the reader to section \ref{section:geometry}.
\end{absatz}
 \begin{absatz}{{\it{Genus and volume formula}}.}
 We assume $N \geq 5$ to avoid elliptic fixed points and to have uniform formulas for the following quantities.
 Let $g_N$ be the genus and  $v_N$ the hyperbolic volume of $X\big(\Gamma_1(N) \big)$,   then we have
\begin{align} \label{align:genus}
  g_N = 1 + \frac{1}{24}  \varphi(N) N \prod_{p \mid N} \left(1+ \frac{1}{p} \right)  - \frac{1}{4} \sum_{d \mid N} \varphi(d) \varphi( N / d)
  \end{align}
  and
  \begin{align} \label{align:volume}
v_N =\frac{\pi}{6}  \varphi(N) N \prod_{p |N} \left(1+ \frac{1}{p} \right),
\end{align}
where $\varphi( \, \cdot \,)$ is Euler's phi function (see \cite{diamondshurman}, theorem 3.1.1, p. 68 and formula (5.15), p. 182). In particular, we have
$g_N \geq 1$ for $N=11$ or $N \geq 13$.
\end{absatz}

%
 
  \begin{absatz}{\it{The Hilbert space  $L^2\big( \Gamma_1(N) \backslash \HHH,k \big)$.}} 
Recall that a function $f: \HHH \longrightarrow \CCC$ is called automorphic  of weight $k$, $k \in \NNN$,
 with respect to $\Gamma_1(N)$, if it satisfies   $f(\gamma z) =  j_{k,\gamma}(z) f(z)$ for all 
$\gamma= \left(\begin{smallmatrix} a & b \\ c & d \end{smallmatrix} \right) \in \Gamma_1(N)$,
where $j_{k,\gamma}(z) :=  \frac{(cz+d)^k}{|cz+d|^k} $. 
 We denote by  $L^2\big(\Gamma_1(N) \backslash \HHH, k\big)$ the Hilbert space of square-integrable automorphic functions of weight $k$ 
with respect to $\Gamma_1(N)$ with scalar product given by
\begin{align*}
\langle f, g \rangle := \int_{\Gamma_1(N)\backslash \HHH}  f(z) \overline{g}(z) \hyp(z) {\hspace{1cm}}  \left(  f,g \in  L^2 \big(\Gamma_1(N) \backslash \HHH, k \big) \right).
\end{align*}
The hyperbolic Laplacian of weight $k$
\begin{align} \label{align:scalar}
\Delta_k := - y^2 \left( \frac{\partial^2}{\partial x^2} + \frac{\partial^2}{\partial y^2}  \right) + iky  \frac{\partial}{\partial x},
\end{align}
 acting as a non-negative self-adjoint operator on $ L^2\big(\Gamma_1(N) \backslash \HHH, k\big)$ (in fact, it is the unique self-adjoint extension of $\Delta_k$ acting
on the subspace of smooth and compactly supported automorphic functions of weight $k$; 
 see \cite{roelcke}, pp. 309--310), gives rise to the spectral decomposition 
\begin{align*}
  L^2 \big(\Gamma_1(N) \backslash \HHH,k \big) = L^2 \big(\Gamma_1(N) \backslash \HHH,k \big)_{0} \oplus L^2 \big(\Gamma_1(N) \backslash \HHH,k \big)_{\mathrm{r}}
  \oplus  L^2 \big(\Gamma_1(N) \backslash \HHH,k \big)_{\mathrm{c}};
\end{align*}
here $ L^2 \big(\Gamma_1(N) \backslash \HHH,k \big)_{0}$ is the space of cusp forms of weight $k$, i.e., of automorphic functions of weight $k$ with vanishing $0$-th Fourier coefficients in its
Fourier expansions with respect to the diverse cusps of $\Gamma_1(N)$,  and is discrete,
 $ L^2 \big(\Gamma_1(N) \backslash \HHH,k \big)_{\mathrm{r}}$ is the discrete part of $ L^2 \big(\Gamma_1(N) \backslash \HHH,k \big)_{0}^{\perp}$,  given by residues of Eisenstein series,
and $ L^2 \big(\Gamma_1(N) \backslash \HHH,k \big)_{\mathrm{c}} $ forms the continuous part, given by integrals of Eisenstein series.
\end{absatz}

 \begin{absatz}{\it{Eisenstein series and spectral expansion.}} 
  Let $\mathfrak a \in P_{\Gamma_1(N)}$ be a cusp of $\Gamma_1(N)$ and $\sigma_{\mathfrak a} \in {\mathrm{SL}}_2(\RRR)$ a scaling matrix of $\mathfrak a$,  i.e., $\sigma_{\mathfrak a} \infty= \frak a $ and 
\begin{align*}
\sigma_{\frak a}^{-1} {\Gamma_1(N)}_{\frak a } \,  \sigma_{\frak a} \cong \left\{  \begin{pmatrix} 1 & 1 \\ 0 & 1 \end{pmatrix}^m \big| \, m \in \ZZZ \right\}
\end{align*}
for ${\Gamma_1(N)}_{\frak a }$ the stabilizer group of the cusp $\mathfrak a$.
 The Eisenstein series for the cusp $ \mathfrak a$  of weight $k$ with respect to $\Gamma_1(N)$ is defined by (cf. \cite{roelcke2}, p. 291)
\begin{align*}
E_{\mathfrak a,k}(z,s) :=  \sum_{\gamma \in {\Gamma_1(N)}_{\mathfrak a} \backslash \Gamma_1(N)}
 {\im}(\sigma_{\mathfrak a}^{-1} \gamma z)^s j_{k,  \sigma_{\mathfrak a}^{-1} \gamma}(z)^{-1} \hspace{1cm} (s \in \CCC, \re(s) > 1).
\end{align*}
The  Eisenstein series $E_{\mathfrak a,k}(z,s)$ defines an holomorphic function for $\re(s) >1$ and is in this range an automorphic function of weight $k$ in the $z$-variable.  Moreover, it possesses
 a meromorphic continuation to the whole $s$-plane. 
The meromorphically continued Eisenstein serie $E_{\mathfrak a,k}(z,s)$ is holomorphic for $\re(s)= \tfrac{1}{2}$, and the poles at $s$ with $\re(s) > \tfrac{1}{2}$ lie all in the interval $(1/2, 1]$
(cf.  \cite{roelcke2},  Satz 10.3,  Satz 10.4, and  Satz 11.2). 
\newline
If $\{ u_j \}_{j=0}^{\infty}$ is an orthonormal basis of the discrete part of $ L^2 \big(\Gamma_1(N) \backslash \HHH,k \big)$, i.e.,
$ \Delta_k u_j = \lambda_j u_j $, $0= \lambda_0 < \lambda_1 \leq \lambda_2, \ldots $,
 then every $f \in L^2 \big(\Gamma_1(N) \backslash \HHH,k \big) $ has the spectral expansion
\begin{align} \label{align:spectral0}
f(z) =\sum_{j=0}^{\infty} \langle f, u_j \rangle u_j(z) +  
 \sum_{\mathfrak a \in P_{\Gamma_1(N)}} \frac{1}{4 \pi} \int_{- \infty}^{+\infty}   \big\langle f, E_{\mathfrak a,k}\big( \cdot, \tfrac{1}{2} +ir \big) \big\rangle E_{\mathfrak a,k}\big(z,  \tfrac{1}{2}  +ir \big) \dif r,
\end{align}
which converges in the norm topology. If furthermore, $f$ is smooth and bounded, then the sums in \eqref{align:spectral0}  are uniformly convergent on compacta of $\HHH$ (see
 \cite{roelcke2},  Satz 7.2, Satz 12.2, and  Satz 12.3).
\end{absatz}
\begin{absatz}{\it{Shifting operator and eigenspaces of $\lambda_0=0$}.}
Noting that  $\Delta_0$ and $\Delta_2$ have the same eigenvalues 
(see \cite{roelcke}, lemma 3.2), we  define $L^2_{\lambda_j}\big(\Gamma_1(N) \backslash \HHH, k \big)$, $k=0,2$, to be the eigenspace corresponding to the eigenvalue $\lambda_j$.
The differential operator (loc. cit. denoted by $K_0$)
\begin{align}
\Lambda_0:=  i y \frac{\partial}{\partial x}+ y \frac{\partial}{\partial y} :  L^2_{\lambda_j}\big(\Gamma_1(N) \backslash \HHH, 0 \big) \longrightarrow  L^2_{\lambda_j}\big(\Gamma_1(N) \backslash \HHH, 2 \big)
\end{align}
induces for $\lambda_j \not=0$ a bijection
satisfying 
 \begin{align} \label{align:discrete}
 \langle \Lambda_0(f), \Lambda_0(g) \rangle = \lambda_j \langle f, g \rangle \hspace{1cm} \big(f, g \in   L^2_{\lambda_j} \big(\Gamma_1(N) \backslash \HHH, 0 \big) \big).
 \end{align}
 (see \cite{roelcke}, lemma 6.1).  For $\lambda_j=0$ and $k=0$, we have that $ \mathcal L^2_{0}\big(\Gamma_1(N) \backslash \HHH, 0 \big)$ is one dimensional generated by the only residue  of the Eisenstein series at $s=1$ given
 by $v_N^{-1}$ (see \cite{iwaniec}, theorem 11.3).  
For $\lambda_j=0$ and $k=2$,   we have an isometry
 \begin{align} \label{align:cusp2hol}
  L^2_{0}\big(\Gamma_1(N) \backslash \HHH, 2 \big) \longrightarrow S_2\big( \Gamma_1(N) \big)
 \end{align}
by sending $f \mapsto y^{-1} f$(see \cite{roelcke}, Satz 6.3), where $S_2\big( \Gamma_1(N) \big)$ denotes the space of holomorphic cusp forms of weight $2$ for $\Gamma_1(N)$. 
For later purposes we mention that
\begin{align} \label{align:spectral}
\Lambda_0\big(E_{\mathfrak a,0}(z,s)\big) =  s  E_{\mathfrak a, 2}(z,s)
\end{align}
(see \cite{roelcke2}, p. 292, equation $(10.8)$).
\end{absatz}

 \begin{absatz}{\it{Automorphic kernels of weight $0$ and $2$}.} 
\label{absatz:growthcond}
Let $h: \RRR \longrightarrow \CCC$ be a test function, i.e., an even function satisfying 
\begin{enumerate}
\item[(i)]$h(r)$ can be extended holomorphically to the strip $|\im(r)|\! < \! \frac{A}{2}$ $(A \in \! \RRR, A>1),$
\item[(ii)]$h(r) < \! \! < (|r|+1)^{-2-\frac{A-1}{2}}$ for $|\im(r)| < \frac{A}{2}$.
\end{enumerate}
The  inverse  Selberg transform $k_0$ of $h$  of weight $0$ is given by the following three equations  (see \cite{iwaniec}, p. 32):
\begin{align*}
&g(w)= \frac{1}{2 \pi} \int_{- \infty}^{+ \infty} h(r) \exp(-iwr) \dif r &(w \in \RRR),\\
&q\big( e^{v} + e^{-v} -2 \big)=  g(v) & (v \in \RRR), \\
&k_{0}(u)= - \frac{1}{\pi} \int_{- \infty}^{+ \infty} q^{\prime}(u+v^2) \dif v &(u \ge 0).
\end{align*}
The  inverse Selberg transform $k_2$ of $h$ of weight $2$ is given in the similar way with the only change in the last step (see \cite{hejhal0}, p. 402 and p. 455):
\begin{align*}
k_{2}(u):= &- \frac{1}{\pi} \int_{- \infty}^{+ \infty} q^{\prime}(u+v^2) \frac{\sqrt{u+4+v^2}-v}{\sqrt{u+4+v^2}+v} \dif v  & (u \ge 0).
\end{align*}
The automorphic kernel of weight $0$ with respect to $\Gamma_1(N)$ on $\HHH \times \HHH$ is defined by
\begin{align}
K_0(z,w) :=  \sum_{\gamma \in \Gamma_1(N)} k_0\big(u(z, \gamma w)\big),
\end{align}
which is an automorphic function of weight $0$ in the $z$-variable; here the function $u$ is defined as in \eqref{align:hypdistance}. 
Similarly,
the  automorphic kernel  of weight $2$ with respect to $\Gamma_1(N)$ on $\HHH \times \HHH$ is defined by
\begin{align}
K_2(z,w) :=  \sum_{\gamma \in \Gamma_1(N)} k_2\big(u(z, \gamma w)\big) \, \frac{\gamma w- \overline{z}}{z -  \gamma \overline{w}} \, \, j_{2,\gamma}(w),
\end{align}
which is an automorphic function of weight $2$ in the $z$-variable.
\end{absatz}

  \begin{absatz}{\it{Spectral expansion of automorphic kernels.}}
From the orthogonal projections 
 \begin{align*}
 \langle K_0(z,w), u_j(z) \rangle = h_0(r_j)  \overline{u}_j(w)
  \end{align*}
  and 
  \begin{align*}
 \langle K_0(z,w), E_{\frak a,0}(z, \tfrac{1}{2} +ir) \rangle = h_0(r)   \overline{E}_{\frak a,0}(w, \tfrac{1}{2} +ir)
  \end{align*} 
   (see  \cite{iwaniec}, theorem~7.4), 
using the convention to write an eigenvalue $\lambda_j$ as $\lambda_j = \frac{1}{4} + r_j^2$ with $r_j \in \CCC$, we find from  \eqref{align:spectral0}  the following spectral expansion
  \begin{align}  \label{align:spectraldecomp1}
  K_0(z,w)=&\frac{h(\tfrac{i}{2})}{ v_N} + \sum_{j =1}^{\infty} h(r_j) u_j(z) \overline{u}_j(w) + \\ \notag 
  & \frac{1}{4 \pi} \sum_{\frak a \in P_{\Gamma_1(N)}} \int _{- \infty}^{+ \infty} h_0(r)  E_{\frak a,0}(z, \tfrac{1}{2} +ir) \overline{E}_{\frak a,0}(w, \tfrac{1}{2} +ir) \dif r.
  \end{align}
Noting that
  \begin{align*}
 \langle K_2(z,w), u_j(z) \rangle = h_2(r_j)  \overline{u}_j(w)
  \end{align*}
  and 
  \begin{align*}
 \langle K_2(z,w), E_{\frak a,2}(z, \tfrac{1}{2} +ir) \rangle = h_2(r)   \overline{E}_{\frak a,2}(w, \tfrac{1}{2} +ir)
  \end{align*} 
 (cf. \cite{abbesullmo}, lemma 3.1.1), we find from \eqref{align:spectral0} and observations \eqref{align:discrete} and \eqref{align:cusp2hol} the following spectral expansion
 \begin{align} \label{align:spectraldecomp2}
 K_2(z,w)= &h(\tfrac{i}{2}) \sum_{j=1}^{g_N} \im(z) f_j(z)  \im(w) \overline{f}_j(w) + \sum_{j =1}^{\infty} \frac{h(r_j)}{\lambda_j}\Lambda_0(u_j)(z) \overline{\Lambda_0(u_j)}(w) +  \notag \\
& \frac{1}{4 \pi} \sum_{\frak a \in P_{\Gamma_1(N)}} \int_{- \infty}^{+ \infty} h(r) E_{\frak a, 2}(z, \tfrac{1}{2} + i r)  \overline{E}_{\frak a, 2}(w, \tfrac{1}{2}  + i r) \dif r,
\end{align}
where $\{f_1,  \ldots , f_{g_{N}} \}$ is an orthonormal basis of $S_2\big(\Gamma_1(N)\big)$.
\end{absatz}

\begin{absatz}
Let $h(t,r):= \exp \left( -t \left( \tfrac{1}{4} + r^2 \right) \right)$ be the test function  with parameter $t\in \RRR_{>0}$.
The inverse Selberg transform of $ h(t,r)$ of weight $0$  defines  the function $k_0(t,u)$ for $u \geq 0$;  for a fixed $\gamma \in \Gamma_1(N)$, we set  
$$
K_{0, \gamma}(t,z):= k_0\big(t,u(z, \gamma z)\big) \hspace{1cm} \left(z \in \HHH \right).
$$
 Similarly, the inverse Selberg transform of $h(t,r)$ of weight $2$ defines the function  $k_2(t,u)$ for $u \geq 0$;  for  a fixed $\gamma \in \Gamma_1(N)$, we set  
 $$
 K_{2, \gamma}(t,z):= k_2\big(t,u(z, \gamma z)\big) \frac{\gamma z - \overline{z}}{z - \gamma \overline{z}} \, j_{2, \gamma}(z) \hspace{1cm} \left( z \in \HHH \right).
 $$
\end{absatz}
\begin{lem} \label{lem:kerneltrans}
Let be $ l \in \ZZZ$. With the above notation the following  series
\begin{align*}
\sum_{\substack{ \gamma \in \Gamma_1(N) \\ \tr( \gamma) =l }} K_{0, \gamma}(t,z) \hspace{1cm} \mbox{and} \hspace{1cm} 
\sum_{\substack{ \gamma \in \Gamma_1(N) \\  \tr( \gamma) =l }} K_{2, \gamma}(t,z) \hspace{1cm} (z \in \HHH)
\end{align*}
are automorphic functions of weight $0$ with respect to $\Gamma_1(N)$.
\begin{proof}
This follows from  the fact that
$K_{k,\gamma}(t,\delta z) = K_{k,\delta^{-1} \gamma \delta}(t, z)$ for $k=0,2$ and any $\delta= \left(\begin{smallmatrix} a & b \\ c & d \end{smallmatrix} \right) \in \Gamma_1(N)$.  
\end{proof}
\end{lem}
\begin{absatz}{{\it{Notation}}.}
We define for $t>0$, $k=0,2$, and $l \in \ZZZ$ with $|l|>2$
 \begin{align} \label{align:hyperpart}
  H_{k,l}(t,z) := \sum_{\substack{\gamma \in \Gamma_1(N) \\  \tr(\gamma)=l}} K_{ k, \gamma}(t,z) 
\hspace{1cm} H_k(t,z):= \sum_{\substack{l \in \ZZZ \\ |l|> 2}} H_{k,l}(t,z),
\end{align}
  \begin{align} \label{align:parapart}
 P_{k}(t,z) := \sum_{\substack{\gamma \in \Gamma_1(N) \\  |\tr(\gamma)| =2}} K_{ k, \gamma}(t,z), 
 \end{align}
 \begin{align} \label{align:specdecomp2}
C_k(t,z):=\, - \frac{1}{4 \pi}  \sum_{\frak a \in P_{\Gamma_1(N)}} \int_{- \infty}^{+ \infty} h(t,r) | E_{\frak a,k}(z, \tfrac{1}{2} +ir) |^2 \dif r - \frac{2-k}{2} \frac{1}{v_N},
\end{align}
and
 \begin{align} \label{align:parapart}
 D_0(t,z) :=    \sum_{j=1}^{\infty} h(t,r_j) |u_j(z)|^2 \hspace{1cm} D_2(t,z):=\sum_{j=1}^{\infty} \frac{h(t,r_j)}{\lambda_j} |\Lambda_0(u_j)(z)|^2. 
 \end{align}
 \end{absatz}
\begin{absatz} {{\it{Basic formula}.}}
 We assume that  $N=11$ or $N \geq 13$. We define the  $\Gamma_1(N)$-invariant functions (cf. lemma  \ref{lem:kerneltrans}) 
\begin{align} \label{align:wedors}
&H(t,z):= H_2(t,z) - H_0(t,z) \hspace{1cm} P(t,z):= P_2(t,z) - P_0(t,z) \notag \\ 
&C(t,z):=C_2(t,z) - C_0(t,z) \hspace{1,1cm } D(t,z):= D_2(t,z) - D_0(t,z),
\end{align}
such that  observations  \eqref{align:spectraldecomp1} and \eqref{align:spectraldecomp2}, taking the difference of $K_2(t,z)$ and $K_0(t,z)$, imply
\begin{align} \label{align:basicformula}
 &g_N F(z) +  D(t,z) = H(t,z) + P(t,z) + C(t,z), 
\end{align}
where 
\begin{align} \label{align:comparisonhypcan}
F(z):=\frac{1}{g_{N}} \,  \,\sum_{j=1}^{g_{N}} y^2  |f_j(z)|^2
\end{align}
 with $\{ f_j \}_{j=1}^{g_N}$ an orthonormal basis of $S_2\big( \Gamma_1(N) \big)$.
Note that we have $h(t,\tfrac{i}{2})=1$
and that  there is no elliptic contribution, i.e., there is no $\gamma \in \Gamma_1(N)$ with $|\tr(\gamma)|< 2$. 
\end{absatz}

\section{Green's function on cusps}
In this section we recall the definition of the canonical Green's function and derive  a formula for its evaluation on cusps essentially in terms of the function defined in \eqref{align:comparisonhypcan}. 
Our formula follows from  previous work of  A. Abbes and E. Ullmo in \cite{abbesullmo} and J. Jorgenson and J. Kramer in \cite{jorgensonkramer2}.
In the sequel we assume $g_N \geq 1$, i.e., that $N= 11$ or $N \geq 13$.
\begin{absatz}{\it{Canonical volume}.}
Let $S_2\big(\Gamma_1(N)\big)$ be the space of holomorphic cusp forms  of weight $2$  with respect to $\Gamma_1(N)$ equipped with the Petersson  inner product 
\begin{align*}
\langle f, g \rangle_{\rm Pet,2} :=  \int_{X(\Gamma_1(N))} f(z)  \overline{g}(z) \im(z)^2 \hyp(z) \hspace{1cm} \big(f, g \in S_2\big(\Gamma_1(N)\big) \big).
\end{align*}
Choosing an orthonormal basis  $\{f_1, \ldots, f_{g_{N}} \}$   of $S_2\big(\Gamma_1(N)\big)$,
 the canonical volume form on $X\big(\Gamma_1(N) \big)$ is given by  
 \begin{align*}
 \can(z) :=  \frac{i}{2 g_N} \sum_{j=1}^{g_N} | f_j(z) |^2 \dif z \wedge \dif \overline{z} = F(z) \hyp(z),
 \end{align*}
where $F(z)$ is defined by \eqref{align:comparisonhypcan}.  
Note that  using the ismomophism  $ S\big( \Gamma_1(N) \big) \cong H^0\big(X\big(\Gamma_1(N)\big), \Omega^1_{X(\Gamma_1(N))} \big)$ given by
$f(z) \mapsto f(z) \dif z$ this volume form is the one considered in \cite{arakelov}.
\end{absatz}
\begin{absatz} {\it{Canonical Green's function.}} \label{absatz:Green}
 The canonical Green's function  $g_{\mathrm{can}}$  is the unique smooth function on $X\big(\Gamma_1(N)\big) \times X\big(\Gamma_1(N)\big) \setminus \Delta_{X( \Gamma_1(N))}$, denoting
 by $\Delta_{X( \Gamma_1(N))}$ the diagonal,
 which satisfies:
 \begin{enumerate}
 \item[\rm{(i)}] $\frac{1}{\pi i}\frac{\partial^2}{\partial z \partial \overline{z}} g_{\mathrm{can}}(z,w) + \delta_w(z) = \can(z)$,
 \item[\rm{(ii)}] $\displaystyle \int_{X(\Gamma_1(N))} g_{\mathrm{can}}(z,w)\can(z) = 0$,
 \end{enumerate}
 where $ \delta_w(z)$ is the Dirac delta distribution.
 \end{absatz}
 
 \begin{absatz}
To state the formula for the canonical Green's function on cusps from \cite{abbesullmo},
 we recall that, for $\mathfrak a, \mathfrak b  \in P_{\Gamma_1(N)}$  cusps of $\Gamma_1(N)$,  the Eisenstein series $E_{\mathfrak a,0}(z,s)$ 
of weight zero with respect to $\Gamma_1(N)$ admits at the cusp $\mathfrak b$ the Fourier expansion $ E_{\mathfrak a,0}(\sigma_{\mathfrak b} z,s) = \sum_{n \in \ZZZ} a_n(y,s;\mathfrak a \mathfrak b) \exp(2 \pi i nx)$ with
$a_0(y,s; \mathfrak a \mathfrak b) = \delta_{\mathfrak a \mathfrak b} y^{s} + \varphi_{\mathfrak a \mathfrak b}(s) y^{1-s},$ 
where 
\begin{align} \label{align:fourcoef}
\varphi_{\mathfrak a \mathfrak b}(s)= \sqrt{\pi} \frac{\Gamma \left(s- \tfrac{1}{2} \right)}{\Gamma(s)}  \frac{1}{N^s} \sum_{c \in \NNN_{>0}} \frac{b_{\mathfrak a \mathfrak b }(c)}{c^{2s}}
\end{align}
with
\begin{align} \label{align:fourcoef1}
b_{\mathfrak a \mathfrak b}(c) = \# \left\{  \begin{pmatrix} \star & \star \\ c & \star \end{pmatrix}  \in g_{\mathfrak a}^{-1} \Gamma_1(N)_{\mathfrak a} \,  g_{\mathfrak a} \backslash 
g_{\mathfrak a}^{-1} \Gamma_1(N) g_{\mathfrak b} / g_{\mathfrak b}^{-1} \Gamma_1(N)_{\mathfrak b} \, g_{\mathfrak b} \right\};
\end{align}
here $g_{\mathfrak a}$ and $g_{\mathfrak b}$ denote elements of $\mathrm{SL}_2(\ZZZ)$ mapping the standard cusp $\infty$ of $\mathrm{SL}_2(\ZZZ)$ with represantative $(1:0) \in \PPP^1_{\QQQ}$
 to the cusps $\mathfrak a$ and $\mathfrak b$, respectively. Note that
$\varphi_{\mathfrak a \mathfrak b}(s)$ is a meromorphic function with a simple pole at $s=1$ and residue
$v_N^{-1}$ (see \cite{iwaniec}).
\end{absatz}
\begin{prop}  \label{prop:abul}
Let $N$ satisfy $N=11$ or $N\geq 13$. Let $X\big(\Gamma_1(N)\big)$ be the modular curve associated to the congruence subgroup $\Gamma_1(N)$. Then, we have for two different cusps $\mathfrak a, \mathfrak b \in P_{\Gamma_1(N)}$
\begin{align*}
&g_{\rm{can}}(\mathfrak a,\mathfrak b)=\\
& -  2 \pi \lim_{s \rightarrow 1} \left( \varphi_{\mathfrak a \mathfrak b}(s) - \frac{1}{v_N}  \frac{1}{ s-1} \right) +
   2 \pi \lim_{s \rightarrow 1} \bigg(  \int_{\Gamma_1(N) \backslash \HHH } F(z) E_{\mathfrak a,0}(z,s) \hyp(z) + \\ 
   &\int_{\Gamma_1(N)\backslash \HHH }  F(w) E_{\mathfrak b ,0} (w,s)  \hyp(w)  -  \frac{2}{ v_N} \frac{1}{s-1}\bigg) + O\left( \frac{1}{g_N} \right),
\end{align*}
where the error term is independent of the cusps $\mathfrak a$ and $\mathfrak b$. 
\begin{proof}
This follows from proposition E in  \cite{abbesullmo} in combination with the bound on the hyperbolic Green's function in  \cite{jorgensonkramer2}, lemma 3.7 and proposition 4.7 with the  universal constants for
$\Gamma_1(N)$ given in lemma 5.3 (c) and lemma 5.9.
\end{proof}
\end{prop}

 \begin{lem} \label{lem:red1}
Let $N$ satisfy $N=11$ or $N\geq 13$. Let $X\big(\Gamma_1(N)\big)$ be the modular curve associated to the congruence subgroup $\Gamma_1(N)$ and $0, \infty_d$, $d \in (\ZZZ/N \ZZZ)^{\times}$,
 the cusps having representatives $(0:1)$ and $(d:0)$ in $\PPP^1_{\QQQ}$, respectively. With the notation $\infty= \infty_1$, we have 
\begin{align} \label{align:redtoinfty}
    \int_{\Gamma_1(N) \backslash \HHH}  F(z) E_{0,0}(z,s) \hyp(z) = \int_{\Gamma_1(N) \backslash \HHH}  F(z) E_{\infty,0}(z,s)  \hyp(z)
\end{align}
and 
\begin{align} \label{align:redtoinfty1}
    \int_{\Gamma_1(N) \backslash \HHH}  F(z) E_{\infty_d,0}(z,s) \hyp(z) = \int_{\Gamma_1(N) \backslash \HHH}  F(z) E_{\infty,0}(z,s)  \hyp(z).
\end{align}
\begin{proof}
We choose $\sigma_0^{-1} = \frac{1}{\sqrt{N}}  W_N \in {\mathrm{SL}}_2(\RRR)$ with $W_N = \left( \begin{smallmatrix} 0 & 1 \\  -N  & 0  \end{smallmatrix} \right)$ the Atkin-Lehner involution and 
$ \sigma_d^{-1} = \left(\begin{smallmatrix} a & b \\ c & d \end{smallmatrix} \right) \in {\mathrm{SL}}_2(\ZZZ)$ with $c \equiv 0 \mod N$.
Then, $\sigma_0$ and $\sigma_d$ are scaling matrices of the cusps $0$ and $\infty_d$, respectively, i.e., we have $\sigma_0 \infty =0$ as well as
 $\sigma_{0}^{-1} \Gamma_1(N) \, \sigma_{0}  =  \Gamma_1(N)$ and
  $\sigma_{0}^{-1} \Gamma_1(N)_{0}\,  \sigma_{0}  = \Gamma_1(N)_{\infty}$ and the same for $\sigma_d$.
Hence it follows  from the definitions that $E_{0,0}(z,s) = E_{\infty,0}(\sigma_{0}^{-1}z, s)$ and 
 $E_{\infty_d}(z,s) = E_{\infty}(\sigma_d^{-1} z ,s)$.
Therefore, it suffices for the proof of equations \eqref{align:redtoinfty} and \eqref{align:redtoinfty1} to show  $F(z) =  F(\sigma_{0}^{-1} z)$ and $F(z)=F(\sigma_d^{-1} z) $, which we do now starting with the first equality.
  \newline
The Atkin-Lehner involution $W_N$ acts on the space   $S_2(\Gamma_1(N))$ of holomorphic cusp forms of weight $2$ with respect to $\Gamma_1(N)$ by
 \begin{align*}
  f|_{W_N}( z) := \det(W_N) (-Nz)^{-2} f(W_N z)= N (Nz)^{-2} f( W_Nz) \hspace{0,7cm} \big( f \in S_2\big(\Gamma_1(N)\big) \big),
  \end{align*}  
  and we have $|f  \mid_{W_N} \mid_{W_N}| = | f |$ for $ f \in S_2\big(\Gamma_1(N)\big)$
    (see \cite{atkinli}, proposition 1.1),
 from which we can deduce that $\{ f_j \vert_{W_N} \}_{j=1}^{g_N}$ remains an orthonormal basis, and so
  \begin{align} \label{align:involequal}
  \sum_{j=1}^{g_N} |f_j(z)|^2 =  \sum_{j=1}^{g_N} | f_j \vert_{W_N}(z)|^2.
  \end{align}
  Using  equation \eqref{align:involequal} we calculate 
  \begin{align*}
F(\sigma_0^{-1}z)& = F(W_N z)  =\frac{1}{g_N}  \sum_{j=1}^{g_N} {\im}(W_N z)^2 |f_j ( W_N z)|^2   = \frac{1}{g_N}   \frac{N^2}{\vert  N z  \vert^4} \sum_{j=1}^{g_N} y^2 |f_j ( W_N z) |^2 \\
& =  \frac{1}{g_N}    \sum_{j=1}^{g_N} y^2 |f_j \vert_{W_N}(z) |^2  =  \frac{1}{g_N}  \sum_{j=1}^{g_N}  y^2 | f_j (z)|^2 = F(z).
  \end{align*}
  This proves the first equality.
   For the second equality, note that 
  the space $S_2\big(\Gamma_1(N) \big)$  of cusp forms of weight $2$ decomposes in
\begin{align}
S_2\big(\Gamma_1(N) \big) = \bigoplus_{\epsilon} S_2\big(\Gamma_1(N), \epsilon \big),
\end{align}
where $\epsilon$ runs through all Dirichlet characters mod $N$ (see \cite{rohrlich}, p. 77). Thereby, we define 
$$
S_2\big(\Gamma_1(N), \epsilon \big):= \big\{ f  \in S_2\big(\Gamma_1(N)\big) \, \big| \, f|_2 \langle d \rangle = \epsilon(d) f \text{ for all } d \in  (\mathbb{Z}/N \mathbb{Z})^{\times} \big\}
$$ 
denoting by $\langle d \rangle$ the diamond operator given by
$$
f|_2 \langle d \rangle (z) := \frac{1}{(cz +d)^2} f(\gamma z) \hspace{1cm} \text{for some } \gamma = \begin{pmatrix} a & b \\ c & d' \end{pmatrix} \in \Gamma_0(N),
$$
where $d' \equiv d \mod N$. But since $ f_j |_2 \langle d \rangle = \epsilon(d) f_j $ for some Dirichlet character $\epsilon$, 
the set  $\{f_j|_2 \langle d \rangle \}_{j=1}^{g_N}$ remains an orthonormal basis of $S_2\big(\Gamma_1(N)\big)$. Further, we have
$ \im(\sigma_d^{-1} z)^2 |f_j(\sigma_d^{-1} z)|^2 =  \im (z)^2 |f_j( z)|^2$
showing  $F(\sigma_d^{-1} z) = F(z)$. This completes the proof of the lemma.
\end{proof}
\end{lem}
\begin{lem} \label{lem:sameinfty}
Let $N$ satisfy $N=11$ or $N\geq 13$. Let $X\big(\Gamma_1(N)\big)$ be the modular curve associated to the congruence subgroup $\Gamma_1(N)$ and   $0, \infty_{d}$, $d \in (\ZZZ/N \ZZZ)^{\times}$
  the cusps having representatives $(0:1)$ and $(d:0)$ in $\PPP^1_{\QQQ}$, respectively. With the notation $\infty= \infty_1$, we then have that
$\varphi_{0 \infty_{d}}(s) = \varphi_{0 \infty}(s)$ holds for all   $d \in (\ZZZ/N \ZZZ)^{\times}$.
\begin{proof}
Noting that
$g_{0}^{-1} \Gamma_1(N)_{\infty} \, g_{0} = \left\langle \left(\begin{smallmatrix} 1 & N \\ 0 & 1 \end{smallmatrix} \right) \right\rangle$  and 
$g_{d}^{-1} \Gamma_1(N)_{\infty} \, g_{d} = \left\langle \left( \begin{smallmatrix} 1 & 1 \\ 0 & 1 \end{smallmatrix} \right) \right\rangle$
with $g_0 = \left( \begin{smallmatrix} 0 & -1 \\ 1& 0 \end{smallmatrix} \right)$ and $g_{d} = \left( \begin{smallmatrix} a & b \\ N& d \end{smallmatrix} \right)  \in {\mathrm{SL}}_2(\ZZZ)$,
we have to show by formula \eqref{align:fourcoef} that the number of elements of the sets
\begin{align*}
S_{d}(c) := \left\{  \begin{pmatrix} \alpha & \beta \\ c & \delta \end{pmatrix}   \in \left\langle \begin{pmatrix} 1 & N \\ 0 & 1 \end{pmatrix} \right\rangle \Big\backslash
 g_0^{-1} \Gamma_1(N) g_{d} \Big/  \left\langle \begin{pmatrix} 1 & 1 \\ 0 & 1 \end{pmatrix} \right\rangle \right\}
\end{align*}
is independent of $d$. 
To this end, we  consider the map
\begin{align*}
\psi : S_{d}(c) \longrightarrow (\ZZZ/c \ZZZ)^{\times}
\end{align*}
induced by mapping $\left( \begin{smallmatrix} \alpha & \beta \\ c & \delta \end{smallmatrix}\right) \mapsto \delta \mod c$.
Since $(c,\delta)=1$ and the right action by $\langle \left( \begin{smallmatrix} 1 & 1 \\ 0 & 1 \end{smallmatrix} \right) \rangle$  changes $\delta$ only by $\hspace{-2mm}\mod c$, the map $\psi$ is well-defined.
We now show  that $\psi$ is bijective starting from showing that $\psi$ is injective. This will prove  the lemma. 
\newline
Let $\left( \begin{smallmatrix} \alpha_1 & \beta_1 \\ c & \delta_1 \end{smallmatrix} \right)$ and  $\left( \begin{smallmatrix} \alpha_2 & \beta_2 \\ c & \delta_2 \end{smallmatrix} \right)$ be two representatives 
of elements of $S_{d}(c)$ such that $\delta_1 \equiv \delta_2 \mod c$, i.e., $\delta_2 = \delta_1 + nc$ for some $n \in \ZZZ$.
By the right action  of $ \left( \begin{smallmatrix} 1 & n \\ 0 & 1 \end{smallmatrix} \right)$ on $\left( \begin{smallmatrix} \alpha_1 & \beta_1 \\ c & \delta_1 \end{smallmatrix} \right)$, we obtain
$\left( \begin{smallmatrix} \alpha_1 & \beta_1 \\ c & \delta_1 \end{smallmatrix} \right)  \left( \begin{smallmatrix} 1 & n \\ 0 & 1 \end{smallmatrix} \right) = 
 \left( \begin{smallmatrix} \alpha_1 & \star \\ c & \delta_2 \end{smallmatrix} \right)$.
From this we conclude $\alpha_1 \delta_1 \equiv \alpha_2 \delta_2 \mod c$, i.e., $c | (\alpha_2 - \alpha_1)$. Furthermore, being elements of $S_d(c)$, we have $N | \alpha_1$ and $N| \alpha_2$; since $(c,N)=1$, we find that
$cN | (\alpha_2 - \alpha_1)$, which shows $\alpha_2 = \alpha_1 +mcN$ for some $m \in \ZZZ$. By the left action of $\left( \begin{smallmatrix} 1 & mN \\ 0 & 1 \end{smallmatrix} \right)$ on
$ \left( \begin{smallmatrix} \alpha_1 & \star \\ c & \delta_2 \end{smallmatrix} \right)$, we obtain
$\left( \begin{smallmatrix} 1 & mN \\ 0 & 1 \end{smallmatrix} \right) \left( \begin{smallmatrix} \alpha_1 & \star \\ c & \delta_2 \end{smallmatrix} \right)=
 \left( \begin{smallmatrix} \alpha_2 & \star \\ c & \delta_2 \end{smallmatrix} \right)$,
which proves the injectivity of $\psi$.
\newline
We now show the surjectivity of $\psi$. Let $\delta \mod c$ be given in $(\ZZZ/c \ZZZ)^{\times}$; we let $\delta \in \ZZZ$ be a representative satisfying $(c,\delta)=1$. We have to find $\alpha, \beta \in \ZZZ$ such that
the representative $\left( \begin{smallmatrix} \alpha & \beta \\ c & \delta \end{smallmatrix} \right)$ 
satisfies $\alpha \equiv 0 \mod N$, $\beta \equiv   d \mod N$, and $\alpha \delta - \beta c =1$,
which, in fact, implies
$g_0 \left( \begin{smallmatrix} \alpha & \beta \\ c & \delta \end{smallmatrix} \right) g_{d}^{-1} \in \Gamma_1(N)$,
as desired. The first and second condition forces us to choose $\alpha$ and $\beta$ of the form
$\alpha = x N$ and $\beta =  d  + y N$ 
with $x,y \in \ZZZ$; we have to verify that there are $x,y \in \ZZZ$ such that also the third condition is satisfied. The three conditions imply $- d c \equiv 1 \mod N$; hence, we find
$v \in \ZZZ$ with $vN - d c =1$. Now, since $(c, \delta)=1$, there are $x, y \in \ZZZ$ such that $xd - yc = v$.
With this choice for $x,y \in \ZZZ$, we find $\alpha \delta - \beta c =1$. This completes the proof of the lemma.
\end{proof}

\end{lem}
  \begin{absatz}{{\it{Rankin-Selberg transform}}.} \label{absatz:rankin}
Let $f$ be a   $\Gamma_1(N)$-invariant function of rapid decay at the cusp $\mathfrak a$, i.e.,   the $0$-th Fourier coefficient $a_0(y; \frak a)$ of the Fourier expansion
$f(\sigma_{\frak a} z) = \sum_{n \in \ZZZ} a_n(y; \frak a) \exp(2 \pi i n x)$
of $f$ at the cusp $\mathfrak a$  satisfies $a_0(y; \frak a) = O(y^{-M})$ for all $M >0$ as $y \rightarrow \infty$.
 Then, the Rankin-Selberg transform $R_{f,\frak a}(s)$  of $f$  at the cusp $\mathfrak a$ is defined by
$$
R_{f,\frak a}(s):=\int_{\Gamma_1(N) \backslash \HHH} f(\sigma_{\frak a} z) E_{\frak a,0}(z,s) \hyp(z) = \int_0^{+ \infty} a_0(y; \frak a) y^{s-2} \dif y \hspace{1cm}(\re(s) > 1)
$$
The Rankin-Selberg transform $R_{f,\frak a}(s)$ of $f$ at the
cusp $\mathfrak a$
 can be continued meromorphically to the whole $s$-plane and has  simple poles at $s=0,1$ with residue at $s=1$ given by
\begin{align}
\res_{s=1} \big(R_{f,\frak a}(s) \big)= \frac{1}{v_N} \int_{\Gamma_1(N) \backslash \HHH} f(z) \hyp(z).
\end{align}
(see, e.g.,  \cite{gupta}, p. 9). Applying the Rankin-Selberg transform to the function $F(z)$ defined in \eqref{align:comparisonhypcan},
 which is of rapid decay at all cusps,  we have  the Laurent expansion at $s=1$, writing   $R_F(s):=R_{F, \infty}(s)$, 
\begin{align} \label{align:cfcf}
R_F(s) = \frac{1}{v_N} \frac{1}{s-1} + C_F + O(s-1),
\end{align}
denoting by $C_F$ the constant term of this expansion. 
 \end{absatz}

\begin{prop}  \label{prop:gcanwithrs}
Let $N$ satisfy $N=11$ or $N\geq 13$. Let $X\big(\Gamma_1(N)\big)$ be the modular curve associated to the congruence subgroup $\Gamma_1(N)$ and $0, \infty_d$, $d \in (\ZZZ/N \ZZZ)^{\times}$, 
the cusps having representatives $(0:1)$ and $(d:0)$, respectively.
Then, we have 
 \begin{align*}
g_{{\rm{can}}}(0,  \infty_d)  =  4 \pi  C_F  - 2 \pi \lim_{s \rightarrow 1} \left( \varphi_{0 \infty}(s)- \frac{1}{v_N} \frac{1}{ s-1} \right)  + O  \left(    \frac{1}{g_N} \right),
\end{align*}
where the error term is independent of $d$.
\begin{proof}
This follows from proposition \ref{prop:abul}, lemma \ref{lem:red1}, and lemma \ref{lem:sameinfty}, using the notation of \ref{absatz:rankin}.
\end{proof}
\end{prop}
\begin{rem} \label{rem:strategy}
In the next section we follow a strategy to determine the constant $C_F$  which
 was carried out by A. Abbes and E. Ullmo for the modular curve $X_0(N)$, $N$ 
squarefree and not divisible by $2$ and $3$ in \cite{abbesullmo} based on ideas of  D. Zagier's proof of
the Selberg trace formula in \cite{zagier2}, which goes as follows:  we
will compute the Rankin-Selberg transforms, denoted by 
\begin{align}
 R_H(t,s), \, \, \, R_P(t,s), \,\,\, R_C(t,s), \,\,\, R_D(t,s)  
 \end{align}
 of all terms displayed in  \eqref{align:wedors}  and determine their constant terms in their Laurent expansions at $s=1$.
 Letting $t$ tend
to infinity, the contribution of the discrete part $R_D(t,s)$ will vanish, and so we obtain the constant $C_F$ by formula \eqref{align:basicformula}.
\newline
It might be interesting to look at the problem to determine $C_F$ from  the adelic point of view starting with \cite{jacquetzagier}. 
\end{rem}

\section{Contribution of Rankin-Selberg  of hyperbolic part}

In this section we calculate the contribution of Rankin-Selberg transform $R_H(t,s)$ in terms of the Selberg zeta function.
\subsection{ Rankin-Selberg of hyperbolic part}
We begin with calculating  the Rankin-Selberg transforms of the  $\Gamma_1(N)$-invariant functions (see \eqref{align:hyperpart}  at the end of section 3.5)
\begin{align*}
H_{k,l}(t,z) = \sum_{\substack{\gamma \in \Gamma_1(N) \\ \tr(\gamma)=l}} K_{ k, \gamma}(t,z) \hspace{1cm}(t>0; k=0,2)
\end{align*}
 for $l \in \ZZZ$ with $|l| > 2$.
\newline
  We first note that the Rankin-Selberg transforms  of these functions exist for $s \in \CCC$ with $1<\re(s)< 1+A$ and $A$ as in \ref{absatz:growthcond}. This can be shown mutatis mutandis
 as in \cite{abbesullmo}, proposition 3.2.1: one can reduce the question to show
 \begin{align*}
\sum_{|l| > 2} \int_{2}^{+ \infty} \int_{-1/2}^{1/2}  \sum_{\substack{\gamma \in \Gamma_1(N) \\ \tr(\gamma)=l}} |K_{ 0, \gamma}(t,z)| y^{\re(s) -2}  \dif x \dif y < \infty
\end{align*} 
 for $1< \re(s)<1+ A$, and the claim follows then from \cite{abbesullmo}, lemma 3.2.1. 

\begin{absatz4}
Elements of $\Gamma_1(N)$ give rise to quadratic forms. Let us briefly discuss this link. Therefore, 
we  use  the convention to write $\lbrack a,b,c \rbrack$ for an  (integral binary) quadratic form $q(X,Y)=a X^2+ b XY + c Y^2 \in \ZZZ\lbrack X,Y \rbrack$.
 Let $Q_l$ be the set of quadratic forms with discriminant $ \discr{q} = l^2-4$.
The  modular group $\mathrm{SL}_2(\ZZZ)= \Gamma_1(1)$ acts on the set of  quadratic forms $Q_l$ via
\begin{align} \label{align:action}
\mathrm{SL}_2(\ZZZ)  \times Q_l & \longrightarrow Q_l \\ \notag
(\delta, q)& \mapsto q \circ \delta,
\end{align}
 where $(q \circ \delta)(X,Y):= q\left( (X,Y) \delta^t \right)$
 with $ \delta^t$  the transpose of $ \delta$. For $q=\lbrack a,b,c \rbrack \in Q_l$ and $\delta= \left( \begin{smallmatrix} x & y \\ z & t \end{smallmatrix} \right) \in \mathrm{SL}_2(\ZZZ)$ 
 the quadratic form $q \circ \delta$ is explicitly given by
\begin{align} \label{align:linquad}
q \circ \delta = \lbrack q(x,z), b(xt+yz)+ 2(axy+czt), q(y,t)  \rbrack.
\end{align}
 
  \end{absatz4}

 \begin{defn4} For a positive integer $N$ and $l \in \ZZZ$ with $|l|>2$, we define
\begin{align*}
 Q_l(N):= \big\{ q=\lbrack aN,bN,c \rbrack \, |  \, a,b,c \in \ZZZ ; \discr(q) = l^2 -4 \big\} \subseteq Q_l
\end{align*}
and
 \begin{align*}
 \Gamma_{1,l}(N) = \big\{ \gamma \in \Gamma_1(N) | \tr (\gamma)=l \big\} \subseteq \Gamma_1(N).
 \end{align*}
\end{defn4}
\begin{absatz4}
Let $N$ be a positive integer
and  $l \in \ZZZ$  with $|l|>2$ and $l \equiv 2  \mod N$.
We have a map
$$
\psi : \Gamma_{1,l}(N) \longrightarrow Q_l(N)
$$
defined by
\begin{align} \label{align:corres}
\gamma = \left( \begin{matrix} 1+aN & b        \\
					       cN      & 1+dN
			\end{matrix} \right)
			\mapsto q_{\gamma}:= \lbrack cN, (d-a)N, - b \rbrack.
\end{align}	
Supposed that $N$ is odd,  the map $\psi$
 defines a bijection between the sets $\Gamma_{1,l}(N)$ and $Q_l(N)$ as one verifies 
that  the map  $\psi': Q_l(N) \longrightarrow \Gamma_{1,l}(N)$ given by 
\begin{align*}
q=\lbrack aN, bN, c \rbrack \mapsto \gamma_q:= \left( \begin{matrix} \frac{l -bN}{2} & -c \\ aN & \frac{l+bN}{2} \end{matrix} \right)
\end{align*}
is well-defined. Once this is shown one easily verifies that  the  two maps $\psi$ and $\psi'$ are inverse to each other,  which establishes 
the claimed bijection.
\end{absatz4}

\begin{absatz4} \label{absatz4:matquad}
The congruence subgroup $\Gamma_1(N)$
  acts on the sets  $Q_l(N)$ and $\Gamma_{1,l}(N)$ as follows: The action of $\Gamma_1(N)$ on $Q_l(N)$ is given by
  \begin{align} \label{align:actionN}
\Gamma_1(N) \times Q_l(N) & \longrightarrow Q_l(N) \\ \notag
(\delta, q)& \mapsto q \circ \delta,
\end{align}
where equation \eqref{align:linquad} shows that this action is well-defined;
the action of $\Gamma_1(N)$  on $\Gamma_{1,l}(N)$ is given by  conjugation
\begin{align} \label{align:action1}
\Gamma_1(N) \times \Gamma_{1,l}(N) &\longrightarrow \Gamma_{1,l}(N) \\ \notag
(\delta, \gamma) &\mapsto \delta \cdot \gamma := \delta^{-1} \gamma \delta.
\end{align}
Let $N$ be an odd positive integer
and  $l \in \ZZZ$  with $|l|>2$ and $l \equiv 2  \mod N$. Then  the  diagram
$$
\begin{xy}
  \xymatrix{
       \Gamma_1(N) \times \Gamma_{1,l}(N) \ar[r] \ar[d]_{\id \times \psi }    &   \Gamma_{1,l}(N) \ar[d]^{\psi}  \\ \notag
    \Gamma_1(N) \times Q_l(N)   \ar[r]             &  Q_l(N)   
  } 
\end{xy}
$$
commutes, where the horizontal maps are the group actions \eqref{align:actionN}, \eqref{align:action1}, and the vertical map $\psi$ is given by \eqref{align:corres}. In particular, we have a bijection 
$$
\Gamma_{1,l}(N) / \Gamma_1(N) \cong Q_l(N) / \Gamma_1(N).
$$
\end{absatz4}

\begin{absatz4}
Let  $N$ be an odd positive integer.
For $0 < u < N$ with $(u,N)=1$, we set
$$
M_u(N):= \big\{(m,n) \in \ZZZ^2 \, | \, (m,n) \equiv (0,u) \! \! \! \mod \! N \big\}.
$$ 
Note that the congruence subgroup $\Gamma_1(N)$ acts on $M_u(N)$ by 
\begin{align} \label{align:actionpair}
\Gamma_1(N)  \times M_u(N) & \longrightarrow M_u(N) \\ \notag
(\delta, (m,n))& \mapsto (m,n) \delta.
\end{align}
Furthermore, for  $l \in \ZZZ$  with $|l|>2$ and 
$l \equiv 2  \mod N$, we set $\Delta^u_l(N):=  \Gamma_{1,l}(N) \times M_u(N)$ and  define
\begin{align*}
\Delta^{u\pm}_l(N) :=& \left\{\big(\gamma, (m, n) \big) \in\Delta^u_l(N)  \,  | \,  q_{\gamma}(n,-m) \gtrless 0 \right\}.
\end{align*}
 Note that $$\Delta^u_l(N)= \Delta^{u+}_l(N) \mathbin{\dot{\cup}}  \Delta^{u-}_l(N),$$ and that there is an  action of the congruence subgroup $\Gamma_1(N)$  on $\Delta^{u}_l(N)$ defined by
 \begin{align*} 
 \Gamma_1(N) \times \Delta^{u}_l(N) &\longrightarrow \Delta^{u}_l(N) \\
\big(\delta, \big(\gamma, (m, n)\big) \big)  &\mapsto \delta \cdot  \big(\gamma, (m, n)\big) := \big(\delta^{-1} \gamma \delta, (m, n) \delta \big),
\end{align*}    
which preserves the subsets $\Delta^{u+}_l(N),\Delta^{u-}_l(N) \subseteq \Delta^{u}_l(N)$.
 Observing that
\begin{align*} 
\left( q \circ \delta \right)  \cdot ((m,n) \delta) = q \cdot (m,n) 
\end{align*}
using the notation $q \cdot (m,n):= q(n,-m)$ as in \cite{zagier}, we make the following
\end{absatz4}
\begin{defn4}
 Let $N$ be an odd positive integer
and  $l \in \ZZZ$  with $|l|>2$ and 
$l \equiv 2  \mod N$.
For $0< u < N$ with $(u,N)=1$, we define for  $s \in \CCC$, $\re(s)>1$, the zeta function
\begin{align} \label{align:observ}
\zeta_{u,N}(s,l)=   \sum_{(\gamma, (m,n)) \in \Delta^{u+}_l(N)  / \Gamma_1(N)}  \frac{1}{q_{\gamma}(n,-m)^{s}}.
\end{align}
\end{defn4}
We show in the appendix that this zeta function is well-defined for $\re(s)>1$.
 Finally, we denote for any $l \in \ZZZ$ by $\gamma_l$ the matrix  $\gamma_l:= \left( \begin{smallmatrix}  \frac{l}{2} &  \frac{l^2}{4}-1 \\ 1 & \frac{l}{2}  \end{smallmatrix} \right) \in {\rm{SL}}_2(\RRR)$.
  \begin{prop4} \label{prop4:arithanaly}
 Let $N$ be an odd positive integer and $l \in \ZZZ$ with  $|l| >2$ and $l \equiv 2 \mod N$. For $s \in \CCC$ satisfying $1 < {\rm{Re}}(s) < 1+ A$, where $A$ is as in  \ref{absatz:growthcond}, and $t>0$, we have
 \begin{align*}
  \int_{\Gamma_1(N) \backslash \HHH}  H_{k,l}(t,z)  E_{\infty,0}(z,s)\hyp(z) =  I_{k,l}(t,s) \sum_{\substack{0 < u < N \\ (u,N)=1}}  c_u(s) \zeta_{u,N}(s,l)  \hspace{1cm} (k =0,2),
  \end{align*}
  where we set
 $$
I_{k,l}(t,s):= \int_{\HHH} K_{k,\gamma_l}(t,z) y^s \hyp(z) + \int_{\HHH} K_{k,\gamma_{-l}}(t,z) y^s \hyp(z) \hspace{1cm} (k=0,2)
$$
and
 $$
 c_u(s):=\sum_{\substack{d>0 \\ du \equiv 1\! \! \! \mod N}} \frac{\mu(d)}{d^{2s}}
$$
denoting by $\mu(d)$  the Moebius function.
\begin{proof}
First note that
\begin{align} \label{align:congeisen}
 E_{\infty,0}(z,s)  = &  \sum_{\Gamma_1(N)_{\infty} \backslash \Gamma_1(N)} \! \! \!\!\! \im(\gamma z)^s  =  \sum_{\left(\begin{smallmatrix} * & * \\ m & n \end{smallmatrix} \right) \in \Gamma_1(N) } \frac{y^s}{|mz + n |^{2s}} \notag \\
 = & \sum_{\substack{(m,n) \equiv (0,1)   \!\! \! \mod N  \\ (m,n)=1}} \frac{y^s}{| mz + n |^{2s}}. 
\end{align}
In order to write the Eisenstein series without the coprimality condition, we define
$$
\zeta_N(s) := \sum_{\substack{d>0 \\ (d,N)=1}}\frac{1}{d^{s}} \hspace{1cm} (\re(s) > 1),
$$
and observe that
\begin{align*}
&\zeta_N(2s) \sum_{\substack{(m,n) \equiv (0,1)  \! \!\! \! \mod N  \\ (m,n)=1}} \frac{y^s}{| mz + n |^{2s}} =  \sum_{\substack{d>0 \\ (d,N)=1}}  \sum_{\substack{(m,n) \equiv (0,1)  \! \!\! \! \mod N  \\ (m,n)=1}} \frac{y^s}{| dmz + dn |^{2s}} = \\
&  \sum_{\substack{d>0 \\ (d,N)=1}}  \sum_{\substack{(m',n') \equiv (0,d)  \! \!\! \! \mod N  \\ (m',n')=d}} \frac{y^s}{| m'z + n' |^{2s}} =  
\sum_{\substack{0<u<N \\ (u,N)=1}}  \sum_{\substack{(m',n') \in \ZZZ^2 \\ (m',n') \equiv (0,u)  \! \!\! \! \mod N}} \frac{y^s}{| m'z + n' |^{2s}}. 
\end{align*}
Multiplying by $\displaystyle \zeta_N(2s)^{-1} = \sum_{\substack{d>0 \\ (d,N)=1}}\frac{\mu(d)}{d^{2s}}$ and writing $(m,n)$ instead of $(m', n')$ yields
\begin{align*}
 E_{\infty,0}(z,s)  =
 \sum_{\substack{0 < u < N \\ (u,N)=1}}  c_u(s) \sum_{\substack{(m,n) \in \ZZZ^2  \\ (m,n) \equiv (0,u) \! \!\! \! \mod N}} \frac{y^s}{| mz + n |^{2s}}.
\end{align*}
Now, since $\Gamma_1(N)$ acts freely on $\Delta^u_l(N)$ and  $\Delta^u_l(N)= \Delta^{u+}_l(N) \mathbin{\dot{\cup}} \Delta^{u-}_l(N)$, we obtain 
\begin{align} \label{align:eh}
& \int_{\Gamma_1(N) \backslash \HHH}H_{k,l}(t,z)  E_{\infty,0}(z,s) \hyp(z) = \notag \\
 &\sum_{\substack{0 < u <  N \\ (u,N)=1}}  c_u(s) \sum_{(\gamma,(m,n)) \in \Delta^u_l(N) / \Gamma_1(N)} \int_{\HHH} K_{k, \gamma}(t,z) \frac{y^s}{| mz + n |^{2s}} \hyp(z)= \notag \\
&\sum_{\substack{0 < u < N \\ (u,N)=1}} c_u(s) \Bigg(\sum_{(\gamma,(m,n)) \in \Delta^{u+}_{l}(N) / \Gamma_1(N)}  \int_{\HHH} K_{k, \gamma}(t,z) \frac{y^s}{| mz + n |^{2s}} \hyp(z) +\notag \\  
&\hspace{2,5cm} \sum_{(\gamma,(m,n)) \in \Delta^{u-}_{l}(N) / \Gamma_1(N)}  \int_{\HHH} \! K_{k, \gamma}(t,z) \frac{y^s}{| mz + n |^{2s}} \hyp(z) \Bigg).
\end{align}
For $(\gamma,(m,n)) \in \Delta^{u+}_l(N)$ with $\gamma =  \left( \begin{smallmatrix}1+aN & b \\ cN & 1+dN \end {smallmatrix} \right)$,  we define the matrix
 $$M:= \frac{1}{q_{\gamma}(n, -m)^{\frac{1}{2}}} \left( \begin{matrix} n & -(d-a)N\frac{n}{2} -bm \\ -m & cNn-(d-a)N \frac{m}{2} \end {matrix} \right) \in {\mathrm{SL}}_2(\RRR) .$$
 One computes
 $\frac{\im(Mz)^s}{|mMz + n  |^{2s}} =  \frac{y^s}{q_{\gamma}(n, -m)^s}
 $
 (see \cite{zagier}, p. 127) 
 and an elementary calculation proves the equality $M^{-1} \gamma M = \gamma_l$. 
Recalling that
$K_{k,\gamma}(t,\delta z)= K_{k,\delta^{-1} \gamma \delta}(t,z)$ $(k=0,2)$ 
 for any  $\delta \in  {\mathrm{SL}}_2(\RRR)$, we obtain by the change of variable $z \mapsto Mz$
\begin{align*} 
\int_{\HHH} K_{k, \gamma}(t,z) \frac{y^s}{| mz + n |^{2s}} \hyp(z) = \frac{1}{q_{\gamma}(n,-m)^s}  \int_{\HHH} K_{k, \gamma_{l}}(t,z) y^s \hyp(z),
\end{align*}
and hence
\begin{align} \label{align:trans+}
\sum_{ \Delta^{u+}_{l}(N) / \Gamma_1(N)}  \int_{\HHH} K_{k, \gamma}(t,z) \frac{y^s}{| mz + n |^{2s}} \hyp(z)=
 \zeta_{u,N}(s,l)  \int_{\HHH} K_{k,\gamma_{l}}(t,z) y^s \hyp(z).
\end{align}
For $(\gamma,(m,n)) \in \Delta^{u-}_l(N)$ we define the matrix
 $$M':= \frac{1}{q_{-\gamma}(n, -m)^{\frac{1}{2}}} \left( \begin{matrix} n & (d-a)N\frac{n}{2} +bm \\ -m & -cNn+(d-a)N \frac{m}{2} \end {matrix} \right) \in {\mathrm{SL}}_2(\RRR) $$
and  one verifies  again
 $\frac{\im(M'z)^s}{|mM'z + n  |^{2s}} = \frac{y^s}{q_{-\gamma}(n, -m)^s}$
 as well as  $M'^{-1} (-\gamma) M' = \gamma_{-l}$. 
 Since we have  $K_{k, \gamma}(t,z) = K_{k, - \gamma}(t,z)$,   the change of variable $z \mapsto M'z$ implies
\begin{align*}
\int_{\HHH} K_{k, \gamma}(t,z) \frac{y^s}{| mz + n |^{2s}} \hyp(z) & = \int_{\HHH} K_{k, -\gamma}(t,z) \frac{y^s}{| mz + n |^{2s}} \hyp(z) \\ & = \frac{1}{q_{-\gamma}(n,-m)^s}  \int_{\HHH} K_{k, \gamma_{-l}}(t,z) y^s \hyp(z).
\end{align*}
Observing that $q_{-\gamma} = q_{\gamma^{-1}}$ and $q_{\gamma^{-1}}= - q_{\gamma}$, we find
\begin{align} \label{align:trans-}
&\sum_{(\gamma,(m,n)) \in \Delta^{u-}_{l}(N) / \Gamma_1(N)}  \int_{\HHH} K_{k, \gamma}(t,z) \frac{y^s}{| mz + n |^{2s}} \hyp(z)= \\ \notag
& \sum_{(\gamma, (m,n)) \in \Delta^{u-}_l(N)  / \Gamma_1(N)} \frac{1}{q_{\gamma^{-1}}(n,-m)^s} \int_{\HHH} K_{k, \gamma_{-l}}(t,z) y^s \hyp(z)= \\ \notag 
& \sum_{(\gamma, (m,n)) \in \Delta^{u+}_l(N)  / \Gamma_1(N)} \frac{1}{q_{\gamma}(n,-m)^s} \int_{\HHH} K_{k, \gamma_{-l}}(t,z) y^s \hyp(z)= \\ \notag 
& \zeta_{u,N}(s,l)  \int_{\HHH} K_{k,\gamma_{-l}}(t,z) y^s \hyp(z).
\end{align}
Equation \eqref{align:eh} together with equations \eqref{align:trans+} and \eqref{align:trans-} proves  the proposition.
\end{proof}
\end{prop4}

\subsection{Hyperbolic contribution and Selberg zeta function}
In this section  we  compute the constant term in the Laurent expansion at $s=1$ of the hyperbolic contribution $R_H(t,s)$, i.e., by proposition \ref{prop4:arithanaly}, of
\begin{align*} 
R_H(t,s)= \sum_{\substack{| l | > 2 \\ l \equiv 2 \!\! \mod N}}  \big(I_{2,l}(t,s) - I_{0,l}(t,s)\big) \sum_{\substack{0 < u < N \\ (u,N)=1}}  c_u(s) \zeta_{u,N}(s,l)  \hspace{1cm}(t>0).
\end{align*}

\begin{absatz4} \label{absatz4:pell}
Let $N$ be an odd positive integer
and  $l \in \ZZZ$  with $|l|>2$ and $l \equiv 2  \mod N$.
 We denote by $h_l(N)$ the cardinality of set $Q_l(N)/ \Gamma_1(N)$. The finiteness of $h_l(N)$
follows from the finiteness of the class number $h_l(1)$ of properly equivalent quadratic forms with discriminant $l^2-4$.
\newline
Note that for a quadratic form  $q=\lbrack aN,bN,c \rbrack \in Q_l(N)$ the  stabilizer  $\Gamma_1(N)_{q}= \mathrm{SL}_2(\ZZZ)_q \cap \Gamma_1(N) $ of  $q \in Q_l(N)$ is an infinite cyclic group generated by
$\alpha_q:= \alpha_0^k$, where $\alpha_0$ is the generator of $ \mathrm{SL}_2(\ZZZ)_q$ and $k$ is the least positive integer such that $ \alpha_0^k \in \Gamma_1(N)$.
We denote by $\varepsilon_{q}$ the eigenvalue of $\alpha_q$ with $\varepsilon_{q}^2 >1$ which is, in fact, the $k$-th power of 
the fundamental unit $\varepsilon_0$ in the real quadratic field  $\QQQ(\sqrt{l^2-4})$. In particular, $\varepsilon_q$ is independent of the choice $q$ in $ Q_l(N)$.
Let us also mention that observation \ref{absatz4:matquad} implies $\Gamma_1(N)_{q}= Z(\gamma_q)$, where  $Z(\gamma_q)$ is the centralizer  of $\gamma_q$ in $\Gamma_1(N)$.
\end{absatz4}

\begin{prop4} \label{prop4:continuation}
Let $N$ be  an odd positive integer,  $ l \in \ZZZ$ with $|l|>2$ and $l \equiv 2 \mod N$, and $0< u < N$ with $(u,N)=1$.
Then, the zeta function $\zeta_{u,N}(s,l)$ defines for $\re(s)>1$ a holomorphic function and  has a meromorphic continuation to the whole complex plane with a simple pole at $s=1$ with residue 
\begin{align*}
\res_{s=1}\zeta_{u,N}(s,l) & = \sum_{q \in Q_l(N) / \Gamma_1(N)}  \frac{2}{N^2 \sqrt{l^2 -4}}   \log(\varepsilon_{q}) 
\end{align*}

\end{prop4}

We prove this proposition in the appendix which is a slight variant of a result of  E. Landau in \cite{landau}.

\begin{cor4} \label{cor4:residue}
Let  $N$  be an odd positive integer and $l \in \ZZZ$ with $|l|>2$ and $l \equiv 2 \mod N$.
Then the series  
\begin{align*}
\sum_{\substack{0 < u < N \\ (u,N)=1}} c_u(s) \zeta_{u,N}(s,l) 
\end{align*}
 has a meromorphic continuation to the half plane $\re(s)>1/2$ with a simple pole at $s=1$ and residue 
\begin{align} \label{align:residue}
\res_{s=1}  \sum_{\substack{0 < u < N \\ (u,N)=1}}c_u(s)  \zeta_{u,N}(s,l)
 =  \frac{2}{\pi v_N} 
 \sum_{q \in Q_l(N) / \Gamma_1(N)} \! \! \frac{\log(\varepsilon_q) }{\sqrt{l^2 -4}} .
\end{align}
\begin{proof}
Note that $c_u(s)$ with $0<u < N$ and $(u,N)=1$ is holomorphic for  $\re(s)> 1/2$. Therefore, by  proposition  \ref{prop4:continuation}, we have 
a  meromorphic continuation of $$ \sum_{\substack{0 < u < N \\ (u,N)=1}}c_u(s) \zeta_{u,N}(s,l)$$ to the half plane $\re(s)> 1/2$ and
\begin{align} \label{align:coreq1}
\res_{s=1}   \sum_{\substack{0 < u < N \\ (u,N)=1}} c_u(s) \zeta_{u,N}(s,l) = 
 \sum_{\substack{0 < u < N \\ (u,N)=1}}  c_u(1) \cdot \frac{1}{N^2}\sum_{q \in Q_l(N) / \Gamma_1(N)} \frac{2}{\sqrt{l^2 -4}}   \log(\varepsilon_{q}).
\end{align}
We have
\begin{align} \label{align:coreq7}
\sum_{\substack{0 < u < N \\ (u,N)=1}}  c_u(1)    = \sum_{\substack{d>0 \\  (d, N)=1 }} \frac{\mu(d)}{d^{2}} = \frac{1}{\zeta_N(2)}
   = \frac{1}{ \zeta(2)} \prod_{ p | N} \frac{1}{1- \frac{1}{p^{2}}} = \frac{6}{\pi^2}  \prod_{ p | N} \frac{1}{1- \frac{1}{p^2}}.
\end{align}
Plugging equation \eqref{align:coreq7} in \eqref{align:coreq1} and taking the formula \eqref{align:volume}  into account, the claim of the corollary follows immediately.
\end{proof}
\end{cor4}

Recall that we defined in section 3.3 the test function $h(t,r)$ to be the function
$$h(t,r)= \exp \left( -t \left( \tfrac{1}{4} + r^2 \right) \right) \hspace{1cm} (t>0, r \in \RRR)$$
 with Fourier transform  given by
 \begin{align} \label{align:fouriertransform0}
g(t,w) = \frac{1}{ \sqrt{4 \pi t}} \exp \left( -\tfrac{t}{4} - \tfrac{w^2}{4t} \right) \hspace{1cm} (t>0, w \in \RRR).
\end{align}
\begin{lem4} \label{lem4:I2-I0}
For  any $t>0$, we have the Laurent expansion
\begin{align} \label{align:laurentana}
I_{2,l}(t,s) - I_{0,l}(t,s) = (s-1) A_l(t) + O_t \left( (s-1)^2 \right) \hspace{1cm} (s \rightarrow 1),
\end{align}
where
\begin{align*}
A_l(t) := - \frac{ \pi}{2n_l} + \frac{1}{4} \int_{- \infty}^{+ \infty} \frac{h(t,r)}{\tfrac{1}{4} + r^2} \exp \big( - 2ir \log(n_l)\big) \dif r
\end{align*}
with $n_l:= \frac{l + \sqrt{l^2 -4}}{2}$ and $n_l^{-1}$  the eigenvalues of $\gamma_{l}$.
Furthermore, we have 
\begin{align}  \label{align:laurentana1}
|A_l(t)|  \leq  \frac{\pi}{2 (\log(n_l))^2} n_l^{-  \log (n_l)/t}.
\end{align}
\begin{proof}
For the first statement see \cite{abbesullmo}, proposition 3.2.2. For the second statement, we first observe that
\begin{align} \label{align:thetarel}
- \frac{1}{\pi} A_l(t) = \frac{1}{2n_l} - \frac{1}{4 \pi} \int_{- \infty}^{+ \infty} \frac{h(t,r)}{\tfrac{1}{4} + r^2} \exp \! \big( - 2ir  \log(n_l) \big) \dif r=
\frac{1}{2} \int_{0}^{t} g( \xi, 2 \log (n_l)) \dif \xi.
\end{align}
By definition, we have
\begin{align*}
\frac{1}{2} \int_{0}^{t} g( \xi, 2 \log (n_l)) \dif \xi = \frac{1}{2} \int_{0}^{t} \frac{1}{\sqrt{4 \pi \xi}} \exp\left( - \tfrac{\xi}{4} - \tfrac{(\log(n_l))^2}{\xi} \right)\dif \xi \hspace{0,5cm} \in \RRR_{>0},
\end{align*}
and the change of variable $x:=\frac{1}{\xi}$ yields then
\begin{align*}
- \frac{1}{\pi} A_l(t) =  \frac{1}{2} \int_{1/t}^{+ \infty} \frac{\sqrt{x}}{\sqrt{4 \pi}x^2} \exp\left( - \tfrac{1}{4x} - (\log(n_l))^2x \right) \dif x.
\end{align*}
Since $ \tfrac{\sqrt{x}}{\sqrt{4 \pi}x^2} \exp\left( - \tfrac{1}{4x} \right) < 1$ for $x>0$,
we find 
\begin{align} \label{align:anaconstbound}
- \frac{1}{\pi} A_l(t) \leq \frac{1}{2} \int_{1/t}^{+ \infty}  \exp\left( - (\log(n_l))^2x \right) \dif x  =
\frac{1}{2 (\log(n_l))^2} n_l^{- \log (n_l)/t}. 
\end{align}
from which  the claimed bound follows.
\end{proof}
\end{lem4}
\begin{prop4} \label{prop4:value}
For any $t>0$, the function $R_H(t,s)$ is holomorphic at $s=1$. Furthermore,  we have with  $n_l= \frac{l+ \sqrt{l^2-4}}{2}$
\begin{align*}
R_H(t,1)= & \frac{1}{ \pi v_N}   \sum_{\substack{| l | > 2 \\ l \equiv 2 \!\! \mod N}}
 \sum_{q \in Q_l(N) / \Gamma_1(N)} \frac{2}{\sqrt{l^2 -4}}   \log(\varepsilon_q) \times \\ \notag 
&\left(-\frac{\pi}{2n_l}  + \frac{1}{4} \int_{- \infty}^{+ \infty} \frac{h(t,r)}{\tfrac{1}{4} + r^2} \exp\big(- 2 ir \log(n_l)\big) \dif r \right).
\end{align*}
\begin{proof}
It suffices to prove that $R_H(t,s)$ is bounded at $s=1$.
 From corollary \ref{cor4:residue} and  lemma \ref{lem4:I2-I0}  we find
the inequality 
 \begin{align} \label{align:second}
 \bigg|\big(I_{2,l}(t,1) - I_{0,l}(t,1)\big)  \sum_{\substack{0 < u < N, \\ (u,N)=1}} c_u(1) \zeta_{N,u}(1,l)  \bigg| \leq 
   \frac{n_l^{- \log (n_l)/t}}{2 (\log(n_l))^2} \frac{2 h_l(N)\log(\varepsilon_q)}{ v_N \sqrt{l^2-4}}.
\end{align}
We first note that $h_l(N)$ can be bounded by the classical class number $h_{l}(1)$  by $h_l(N) \leq h_{l}(1) \, \lbrack {\rm SL}_2(\ZZZ) : \Gamma_1(N) \rbrack$
and  that $\varepsilon_{0}^k= \varepsilon_q$ holds for some $0 \leq k \leq \varphi(N)$.
 Since by  Siegel's theorem (see \cite{zagier}, p. 85)
\begin{align} \label{align:siegel}
\sum_{0<D< l^2-4}   \frac{h_D \log(\varepsilon_0)}{\sqrt{D}} = O\big(l^2\big) \hspace{1cm} (l \rightarrow \infty), 
\end{align}
 we obtain
\begin{align} \label{align:supersiegel}
\bigg| \sum_{\substack{|l|>2 \\ l \equiv 2 \! \! \mod N}} \big(I_{2,l}(t,1) - I_{0,l}(t,1)\big) \sum_{\substack{0 < u < N \\ (u,N)=1}} c_u(1)  \zeta_{N,u}(1,l)   \bigg| \leq C_N \sum_{|l|>2} \frac{ l^2 n_l^{-  \log (n_l)/t}}{2 (\log(n_l))^2} 
\end{align}
with $C_N$ some constant depending solely on $N$. As for $l > \exp(3t)$ we have 
$$
 \frac{l^2 n_l^{-  \log (n_l)/t}}{2 (\log(n_l))^2} < l^{-1-\varepsilon}
$$
 with $\varepsilon >0$ small enough,
the series on the right hand side of inequality \eqref{align:supersiegel} converges, which
  proves the holomorphicity of $R_H(t,s)$ at $s=1$.
\newline
The claimed value of $R_H(t,s)$ at $s=1$  follows now  from corollary \ref{cor4:residue} and lemma \ref{lem4:I2-I0}.
\end{proof}
\end{prop4}
%
\begin{defn4}
Let $\gamma \in \Gamma_1(N)$ be a hyperbolic element, i.e., $ | \tr(\gamma)| >2$. The {\it norm} $N(\gamma)$ {\it of} $\gamma$ is defined by 
 $\ N(\gamma) := v^2$, where $v$ is the eigenvalue of $\gamma$ with $v^2 >1$.  We denote by $\gamma_0$ the generator of the centralizer $Z(\gamma)$ of $\gamma$, and
 call $\gamma$ primitive if $\gamma= \gamma_0$ holds.
 \end{defn4}

  \begin{absatz4}
The   Selberg zeta function $Z_{\Gamma_1(N)}(s)$ associated to $\Gamma_1(N)$ is defined via the Euler product expansion
$$
Z_{\Gamma_1(N)}(s): = \prod_{\lbrack \gamma \rbrack \in H(\Gamma_1(N))} Z_{\gamma}(s) \hspace{1cm} ({\re}(s) > 1),
$$
where $ H \big(\Gamma_1(N) \big)$ denotes the set of  primitive conjugacy classes of hyperbolic elements in $\Gamma_1(N)$  and the local factors $ Z_{\gamma}(s)$ are given by
$$
 Z_{\gamma}(s):= \prod_{n=0}^{\infty} \left(1- N(\gamma)^{-(s+n)} \right).
 $$  
The Selberg zeta function $Z_{\Gamma_1(N)}(s)$ is known to have a meromorphic continuation 
to all of $\CCC$ and satisfies a functional equation (cf. \cite{iwaniec}, section 10.8).
\end{absatz4}

\begin{absatz4}
The hyperbolic contribution $\Theta_{\Gamma_1(N)}(t)$ in the Selberg trace formula (cf. \cite{iwaniec}, theorem 10.2) is given by
 \begin{align*}
\Theta_{\Gamma_1(N)}(t) = \sum_{\substack{ \lbrack \gamma \rbrack \\ \mathrm{hyperbolic}}}  \frac{\log\big( N( \gamma_0)\big)}{N(\gamma)^{1/2}- N(\gamma)^{-1/2}} g \big(t, \log \big(N( \gamma) \big) \big).
\end{align*}
For $\gamma \in \Gamma_{1,l}(N)$, one easily verifies, using the fact $\tr(\gamma) = N(\gamma)^{1/2} + N(\gamma)^{-1/2}$, the formulas
\begin{align*}
N(\gamma)^{1/2} - N(\gamma)^{-1/2} = \sqrt{l^2-4} \hspace{0,5cm}  \text{and} \hspace{0,5cm} N(\gamma) = \left( \frac{l + \sqrt{l^2-4}}{2} \right)^2 = n_l^2.
\end{align*}
Further, note that the primitive element $\gamma_0$ associated to the hyperbolic conjugacy class $\lbrack \gamma \rbrack$ generates the centralizer $Z(\gamma)=Z(\gamma_0)$ (see \cite{iwaniec}, p. 137), and hence
equals the generator $\alpha_{q_{\gamma}}$  of the stabilizer $\Gamma_1(N)_{q_{\gamma}}$ of $q_{\gamma}$.
 Therefore, we have 
\begin{align} \label{align:thetaneu}
\Theta_{\Gamma_1(N)}(t) = \sum_{\substack{ | l |>2 \\ l \equiv 2 \! \! \mod  N}} \sum_{q \in Q_l(N) / \Gamma_1(N)} \frac{2 \log\big(  \varepsilon_q)}{\sqrt{l^2 -4}} g \big(t,2 \log( n_l) \big).
\end{align}
For any $t>0$, we have  by formula \eqref{align:thetarel} 
\begin{align} \label{align:Theta}
R_H(t,1)= - \frac{1}{2 v_N}  \int_0^{t} \Theta_{\Gamma_1(N)}(\xi) \dif \xi.
\end{align}
\end{absatz4}
\begin{prop4} \label{prop4:hyper2}
The integral $ \int_0^{t} \Theta_{\Gamma_1(N)}(\xi) \dif \xi$ is asymptotically equivalent to $t$ for $t \rightarrow \infty$, i.e.,
\begin{align*}
 \int_0^{t} \Theta_{\Gamma_1(N)}(\xi) \dif \xi \sim t
\end{align*}
holds. Furthermore, we have 
 \begin{align} \label{align:hypzeta}
\int_0^{+ \infty} \left( \Theta_{\Gamma_1(N)}(t) - 1 \right) \dif t = \lim_{s \rightarrow 1} \left( \frac{Z_{\Gamma_1(N)}^{'}(s)}{Z_{\Gamma_1(N)}(s)} - \frac{1}{s-1} \right)  - 1.
 \end{align}
\begin{proof}
Having McKean's formula (see \cite{mckean}, p. 239)
$$
\frac{1}{2s-1}  \frac{Z_{\Gamma_1(N)}^{'}(s)}{Z_{\Gamma_1(N)}(s)} = \int_0^{+\infty} \exp \! \big(\! -s(s-1)t \big) \Theta_{\Gamma_1(N)}(t) \dif t,
$$
the proof is exactly the same as the proof in \cite{abbesullmo}, proposition 3.3.3, for $\Gamma_0(N)$. 
\end{proof}
\end{prop4}
\begin{cor4} \label{cor4:hyper}
The hyperbolic contribution $R_H(t,1)$ for $t \rightarrow \infty$ is given by
\begin{align*}
R_H(t,1) = -  \frac{t}{2 v_N}  -   \frac{1}{2 v_N } \lim_{s \rightarrow 1}
 \left(  \frac{Z_{\Gamma_1(N)}^{'}(s)}{Z_{\Gamma_1(N)}(s)} -   \frac{1}{s-1} \right) + \frac{1}{2 v_N} + o(1).
\end{align*}
\begin{proof}
This follows now immediately  from observation \eqref{align:Theta} and proposition \ref{prop4:hyper2}.
\end{proof}
\end{cor4}


\section{Contribution of Rankin-Selberg of  spectral and parabolic part}
In this section we determine the contribution of the Rankin-Selberg transforms of the remaining $\Gamma_1(N)$-invariant functions  $P_k(t,s)$,  $C_k(t,s)$, and $D(t,s)$  $(k=0,2)$  of formula \eqref{align:basicformula}.
The contribution of $P_k(t,s)$ can be taken from \cite{abbesullmo} which we cite at the end of the second section of this chapter.
\newline
First observe that we have (for suitable $s$)
\begin{align*}
\int_{\Gamma_1(N) \backslash \HHH} P_{k}(t,z) E_{\infty,0}(z,s) \hyp(z) = \int_0^{+ \infty} p_{k}(t,y) y^{s-2} \dif y \hspace{1cm} (k =0,2),
\end{align*}
and 
\begin{align*}
\int_{\Gamma_1(N) \backslash \HHH} C_{k}(t,z) E_{\infty,0}(z,s) \hyp(z) = \int_0^{+ \infty} c_{k}(t,y) y^{s-2} \dif y \hspace{1cm} (k =0,2),
\end{align*}
where $p_{k}(t,y)$ and $c_{k}(t,y)$
are  the $0$-th Fourier coefficients  of $P_{k}(t,z)$ and $C_k(t,z)$ with respect to the cusp $\infty$, respectively.
Proceeding as in \cite{abbesullmo} we provide  a further decomposition of  the sums $p_{k}(t,y) + c_k(t,y)$, $k =0,2$. Recall that the $0$-th Fourier coefficient 
of the Eisenstein series $E_{\mathfrak a,0}(z,s)$, $\mathfrak a \in P_{\Gamma_1(N)}$, is given by
\begin{align*} 
a_0(y,s; \frak a \infty, 0) = \delta_{\mathfrak a \infty} y^s + \varphi_{\mathfrak a \infty}(s) y^{1-s},
\end{align*} 
where $ \varphi_{\mathfrak a \infty}(s)$ is defined as in  \eqref{align:fourcoef}. Then, 
for $s=\tfrac{1}{2} + ir \in \CCC$, we have
\begin{align} \label{align:zeroweightsplit}
\sum_{\frak a \in P_{\Gamma_1(N)}} \big| a_0\left(y,\tfrac{1}{2} + ir; \frak a \infty, 0 \right) \big|^2 = 2 y + \varphi_{\infty \infty} \left( \tfrac{1}{2} - ir \right)  y^{1+2ir} +  \varphi_{\infty \infty} \left( \tfrac{1}{2} + ir \right) y^{1-2ir}
\end{align} 
 since the scattering matrix $\Phi(s)= \big(\varphi_{\mathfrak a \mathfrak b}(s) \big)_{\mathfrak a, \mathfrak b \in P_{\Gamma_1(N)}}$  is unitary for $\re(s) =\tfrac{1}{2}$
 (see \cite{iwaniec}, theorem 6.6). 
 By observation \eqref{align:spectral}, we find
\begin{align*} 
&\left(\tfrac{1}{2} + ir \right) a_0 \left(y,\tfrac{1}{2} + ir; \frak a \infty, 2\right)=\Lambda_0 \big(a_0 \left(y,\tfrac{1}{2} + ir; \frak a \infty, 0 \right) \big) = \\ \notag
& \left( \tfrac{1}{2} +ir \right) \delta_{\frak a \infty} y^{\frac{1}{2} +ir} + \varphi_{\frak a \infty}\left(\tfrac{1}{2} +ir \right)\left( \tfrac{1}{2} -ir \right)  y^{\frac{1}{2} -ir},
\end{align*}
which implies
\begin{align} \label{align:twoweightsplit}
\sum_{\frak a  \in P_{\Gamma_1(N)}} \big| a_0\left(y,\tfrac{1}{2} + ir; \frak a \infty, 2\right) \big|^2 = 
2 y + & \varphi_{\infty \infty} \left( \tfrac{1}{2} - ir \right) \frac{\tfrac{1}{2} +ir}{\tfrac{1}{2} - ir}  y^{1+2ir}+ \\ &\varphi_{\infty \infty} \left( \tfrac{1}{2} + ir \right) \frac{\tfrac{1}{2} -ir}{\tfrac{1}{2} + ir} y^{1-2ir}. \notag
\end{align}
 Let us define for $k=0,2$ the  functions
  \begin{align} \label{align:p1}
 p_{k,1}(t,y):=  \int_{-\frac{1}{2}}^{+\frac{1}{2}} \sum_{\substack{\gamma \in \Gamma_1(N) \\  |\tr (\gamma) |  =2, \gamma \notin \Gamma_1(N)_{\infty}}} K_{k, \gamma} (t,x + i y) \dif x,
  \end{align}
  \begin{align} \label{align:p2}
 p_{k,2}(t,y):=  \int_{-\frac{1}{2}}^{+\frac{1}{2}} \sum_{ \gamma \in \Gamma_1(N)_{\infty}} K_{k, \gamma} (t,x + i y) \dif x - \frac{y}{2 \pi} \int_{- \infty}^{+\infty} h(t,r) \dif r
  \end{align}
  and
  \begin{align}\label{align:p3}
 c_{k,1}(t,y) :=   - \frac{y}{2 \pi} \int_{- \infty}^{+\infty}   h(t,r) \varphi_{\infty \infty} \! \left(\tfrac{1}{2} - ir \right)   \left( \frac{\tfrac{1}{2} +ir}{\tfrac{1}{2} - ir} \right)^{ \frac{k}{2}}\!\! y^{2ir}  \dif r -
 \frac{2- k}{2}  \frac{1}{v_N},
  \end{align}
  \begin{align} \label{align:p4}
 c_{k,2}(t,y):=  - \frac{1}{4 \pi}  \sum_{\frak a \in P_{\Gamma_1(N)}} \int_{-\frac{1}{2}}^{+\frac{1}{2}}  \int_{- \infty}^{+\infty} h(t,r)  \left| \widetilde{E}_{\mathfrak a, k}\left(x +iy, \tfrac{1}{2} + ir \right) \right|^2 \dif r \, \dif x,
  \end{align} 
with $\widetilde{E}_{\mathfrak a, k}\left(x +iy,  \tfrac{1}{2} + ir \right) := E_{\mathfrak a, k}\left(x +iy,  \tfrac{1}{2} + ir \right) - a_0\left(y, \tfrac{1}{2} +ir; \frak a \infty, k \right)$.
With the above notation,  observations \eqref{align:zeroweightsplit} and \eqref{align:twoweightsplit}  allow us to write
 \begin{align*}
 p_{k}(t,y) + c_{k}(t,y) = p_{k,1}(t,y)+ p_{k, 2}(t,y) + c_{k,1}(t,y)+ c_{k,2}(t,y) \hspace{1cm} (k=0,2).
 \end{align*}
Defining the integrals
 \begin{align*} 
R_{P_{k,j}}(t,s):= \int_0^{+ \infty} p_{k, j}(t,y)  y^{s-2} \dif y \hspace{1cm} (k=0,2;  j=1,2).
\end{align*}
 \begin{align*} 
R_{C_{k,j}} (t,s) := \int_0^{+ \infty} c_{k, j}(t,y)  y^{s-2} \dif y \hspace{1cm} (k=0,2; j=1,2)
 \end{align*}
 we find
\begin{align} \label{align:rankselpar}
 R_{P}(t,s)=  
 \left( R_{P_{2,1}}(t,s) - R_{P_{0,1}}(t,s) \right) + \left( R_{P_{2,2}}(t,s) - R_{P_{0,2}}(t,s) \right) 
\end{align}
and
\begin{align} \label{align:rankselspe}
 R_{C}(t,s) =
\left( R_{C_{2,1}}(t,s)- R_{C_{0,1}}(t,s) \right) +\left( R_{C_{2,2}}(t,s)- R_{C_{0,2}}(t,s) \right). 
\end{align}


 \subsection{Rankin-Selberg of spectral contribution}
 
%
 \begin{lem4}  \label{lem4:RSp_3}
 For $s \in \CCC$ with $1< {\rm{Re}}(s) < A$, where $A$ is as in \ref{absatz:growthcond}, and $t>0$, we have
 \begin{align*}
R_{C_{0,1}} (t,s)=
  -\frac{1}{2} h\left (t, \tfrac{is}{2} \right) \varphi_{\infty \infty} \left( \tfrac{1+s}{2} \right)  
 \end{align*}
 and
 \begin{align*}
R_{C_{2,1}} (t,s)=
  -\frac{1}{2} h\left (t,  \tfrac{is}{2} \right) \varphi_{\infty \infty} \left( \tfrac{1+s}{2} \right)   \frac{1-s}{1+s}. 
 \end{align*}
 \begin{proof}
 The proof is similar to the proof in  \cite{abbesullmo}, lemma 3.2.17 for the congruence subgroup $\Gamma_0(N)$. Since loc. cit.  there is no  proof in the case $k=0$, we give for the convenience of the reader
 a proof here, and refer to \cite{abbesullmo} for the case $k =2$.
 \newline 
The idea is to apply  Mellin's inversion theorem.
  To this end, we consider, for a fixed $t>0$,  the function $h(t,r) \varphi_{\infty \infty} \left(\tfrac{1}{2} - ir \right) y^{2ir}$
 as a function in the complex variable $r$. This function is holomorphic in the strip $ 0 < {\im}(r) < \frac{A}{2}$ except for a simple pole at $r= \frac{i}{2}$, as
$ \varphi_{\infty \infty}  \left(\tfrac{1}{2} - ir \right)$ has a simple pole there with residue $\frac{i}{v_N}$.
 Furthermore, this function
 tends to $0$ for ${\re}(r) \rightarrow \infty$ by  the fact that $\varphi_{\infty \infty}  \left(\tfrac{1}{2} - ir \right)$ is uniformly
 bounded by \cite{hejhal}, theorem 12.9. Now, for  $c \in \RRR$ such that $ \frac{1}{2} < c < \frac{A}{2}$, the residue theorem implies
 \begin{align*}
 & - \frac{y}{2 \pi}  \int_{- \infty}^{ + \infty}  h(t,r) \varphi_{\infty \infty} \left( \tfrac{1}{2} - ir \right) y^{2ir} \dif r  = \\ \notag
 &- \frac{y}{2 \pi} \int_{- \infty +ic }^{+  \infty + ic } \!\!\!  h(t,r) \varphi_{\infty \infty} \left( \tfrac{1}{2} - ir \right) \! y^{2ir}  \dif r 
 \! - \! \frac{y}{2 \pi} \! \left(\frac{ 2\pi  i}{ y} h\left(t, \tfrac{i}{2} \right) {\rm res}_{r=\frac{i}{2}}\left(\varphi_{\infty \infty} \left( \tfrac{1}{2} - ir \right)\right) \! \! \right) \! \! = \\ \notag
 &-\frac{y}{2 \pi}  \int_{ - \infty + ic}^{ +  \infty + ic}\! \!\!   h(t,r) \varphi_{\infty \infty} \left( \tfrac{1}{2} - ir \right) y^{2ir} \dif r 
 + \frac{1 }{v_N}.
\end{align*}
Hence we have
 \begin{align*}
  c_{0,1}(t,y) = -\frac{y}{2 \pi}  \int_{- \infty +ic }^{ +  \infty +ic} h(t,r) \varphi_{\infty \infty} \left( \tfrac{1}{2} - ir \right)  y^{2ir}  \dif r.
 \end{align*}
 If we change the variable by $s = -2ir$, we obtain
 \begin{align*}
  c_{0,1}(t,y) = - \frac{y}{4 \pi i}  \int_{2c - i\infty}^{2c + i \infty} h\left(t, \tfrac{is}{2} \right) \varphi_{\infty \infty} \left( \tfrac{1+s}{2}  \right)  y^{-s}  \dif s,
  \end{align*}
  and from the inverse Mellin transform we deduce the claimed equation
  \begin{align*}
  &R_{C_{0,1}} (t,s)=\int_{0}^{+ \infty} c_{0,1}(t,y) y^{s-2} \dif y =    -\frac{1}{2} h\left (t,\tfrac{i s}{2} \right) \varphi_{\infty \infty} \left( \tfrac{1+s}{2} \right). 
  \end{align*}
 \end{proof}
 \end{lem4}
\begin{absatz4} \label{absatz4:RSEisenstein}
 Let  $f(z)$ be an automorphic function of weight $0$ with respect to $\Gamma_1(N)$ 
with eigenvalue $\lambda$ and of rapid decay at the cusp $\infty$. 
Then, we have
\begin{align*}
\int_{\Gamma_1(N) \backslash \HHH} \!\! \! \! \! \! \! | \Lambda_0(f(z)) |^2 E_{\infty,0}(z,s) \hyp(z) \! = \! \left(\! \lambda + \frac{s(s-1)}{2} \right) \int_{\Gamma_1(N) \backslash \HHH} \!\! \!\! \! \! \!\! |f(z)|^2 E_{\infty,0}(z,s) \hyp(z).
\end{align*}
 The proof is analogous to the proof  of  \cite{abbesullmo},  lemma 3.2.18, where this statement is formulated for  the congruence subgroup $\Gamma_0(N)$.
Namely, the Rankin-Selberg method implies as in \cite{abbesullmo} 
\begin{align} \label{align:shift}
&\int_{\Gamma_1(N) \backslash \HHH}  | \Lambda_0f(z) |^2 E_{\infty,0}(z,s) \hyp(z) = \notag \\ 
& - \frac{1}{2} \int_{\Gamma_1(N) \backslash \HHH} \Delta_0 \left(|f(z)|^2 \right) E_{\infty,0}(z,s) \hyp(z) + \lambda  \int_{\Gamma_1(N) \backslash \HHH} |f(z)|^2 E_{\infty,0}(z,s) \hyp(z). 
\end{align}
 Since  $s(1-s)$ is the eigenvalue of $ E_{\infty,0}(z,s)$, Green's  second identity  implies
\begin{align} \label{align:local}
- \frac{1}{2}\int_{\Gamma_1(N) \backslash \HHH} \!\! \!\!\!\! \Delta_0 \left(|f(z)|^2 \right) E_{\infty,0}(z,s) \hyp(z) = \frac{s(s-1)}{2} \int_{\Gamma_1(N) \backslash \HHH} \!\!\!\!\!\! \!  |f(z)|^2  E_{\infty,0}(z,s) \hyp(z).
\end{align}
Plugging equation \eqref{align:local} in equation \eqref{align:shift} the claim follows.
\end{absatz4}
\begin{lem4} \label{lem4:RSp_4}
 For $s \in \CCC$ with $ {\rm{Re}}(s) > 1$ and $t>0$ the integrals
 \begin{align*}
R_{C_{k,2}}(t,s)= \int_{0}^{+ \infty} c_{k,2}(t,y) y^{s-2} \dif y \hspace{1cm} (k=0,2) 
 \end{align*}
exist and have a meromorphic continuation to the whole $s$-plane with a simple pole at $s=1$. Furthermore, we have
 \begin{align*}
 R_{C_{2,2}}(t,s) =  R_{C_{0,2}}(t,s) + \frac{s(s-1)}{2}  \int_{0}^{+ \infty} c^{*}_{0,2}(t,y) y^{s-2} \dif y,
 \end{align*}
 where $c^{*}_{0,2}(t,y)$ is defined as $c_{0,2}(t,y)$, but with the function $h(t,r)$ replaced by   $\displaystyle h^*(t,r) := \frac{h(t,r)}{\tfrac{1}{4} +r^2}$. 
 \begin{proof}
Recall that we have by definition
  \begin{align*}
 & \int_{0}^{+ \infty} c_{k,2}(t,y) y^{s-2} \dif y =  \\
  & -  \frac{1}{4 \pi} \int_{0}^{+ \infty} \left(   \sum_{\frak a \in P_{\Gamma_1(N)}} \int_{-1/2}^{+1/2}  \int_{- \infty}^{+\infty}  h(t,r)
  \left| \widetilde{E}_{\mathfrak a, k}\left(x +iy, \tfrac{1}{2} + ir \right) \right|^2 \dif r \,  \dif x \right) y^{s-2} \dif y.
  \end{align*}
Since $\widetilde{E}_{\mathfrak a, 0}(x +iy, \tfrac{1}{2} + ir)$  and also $\widetilde{E}_{\mathfrak a, 2}(x +iy, \tfrac{1}{2} + ir)$ are of rapid decay as 
$y \rightarrow \infty$ by \cite{roelcke2}, lemma 10.2, we find that the expressions
\begin{align*}
 \sum_{\frak a \in P_{\Gamma_1(N)}} \int_{-1/2}^{1/2}  \int_{- \infty}^{+\infty}  h(t,r)  \left| \widetilde{E}_{\mathfrak a, k}\left(x +iy, \tfrac{1}{2} + ir \right) \right|^2 \dif r \,  \dif x \hspace{1cm} (k=0,2)
 \end{align*}
satisfy the growth condition of  \ref{absatz:rankin}, and  the first statement of the lemma follows.
\newline
For the second statement, note that we are allowed to interchange the integrals for $\re(s)>1$. Hence we have for $k=0,2$
 \begin{align} \label{align:4to0}
 & \int_{0}^{+ \infty} c_{k,2}(t,y) y^{s-2} \dif y =  \\ \notag
  & -  \frac{1}{4 \pi} \int_{-\infty}^{+\infty}  h(t,r) \left(   \sum_{\frak a \in P_{\Gamma_1(N)}} \int_{\Gamma_1(N) \backslash \HHH}   \left| \widetilde{E}_{\frak a, k}\left(z, \tfrac{1}{2} + ir \right) \right|^2  E_{\infty,0}(z,s) \hyp(z) \right) \dif r.   
  \end{align}
  Since we have
  \begin{align} \label{align:eigenvalueshift}
  \Lambda_0 \left( \widetilde{E}_{\frak a, 0}\left(z, \tfrac{1}{2} + ir \right)\right) = \left(\tfrac{1}{2} + ir\right) \widetilde{E}_{\frak a, 2}\left(z, \tfrac{1}{2} + ir \right),
  \end{align}
  and as $\widetilde{E}_{\frak a, 0}\left(z, \tfrac{1}{2} + ir \right)$ is an eigenfunction of the hyperbolic Laplacian $\Delta_0$
 with eigenvalue $\tfrac{1}{4} + r^2$, which is of rapid decay at the cusp $\infty$, equation \eqref{align:eigenvalueshift} and observation
   \ref{absatz4:RSEisenstein} imply
\begin{align*}
 & \sum_{\frak a \in P_{\Gamma_1(N)}} \int_{\Gamma_1(N) \backslash \HHH}   \left| \widetilde{E}_{\frak a, 2}\left(z, \tfrac{1}{2} + ir \right) \right|^2  E_{\infty,0}(z,s) \hyp(z) = \\
 &\left( 1 + \frac{s(s-1)}{2(\tfrac{1}{4} + r^2)} \right) \sum_{\frak a \in P_{\Gamma_1(N)}} \int_{\Gamma_1(N) \backslash \HHH}   \left| \widetilde{E}_{\frak a, 0}\left(z, \tfrac{1}{2} + ir \right) \right|^2  E_{\infty,0}(z,s) \hyp(z). 
\end{align*}
Therefore, we have by equations \eqref{align:4to0}
\begin{align*}
& \int_{0}^{+ \infty} c_{2,2}(t,y) y^{s-2} \dif y =   \int_{0}^{+ \infty} c_{0,2}(t,y) y^{s-2} \dif y + \frac{s(s-1)}{2} \int_{0}^{+ \infty} c^{*}_{0,2}(t,y) y^{s-2} \dif y, 
\end{align*}
where
\begin{align} \label{align:testneu}
c^{*}_{0,2}(t,y) :=  - \frac{1}{4 \pi}  \sum_{\frak a \in P_{\Gamma_1(N)}} \int_{-\frac{1}{2}}^{+\frac{1}{2}}  \int_{- \infty}^{+\infty} h^{*}(t,r)  \left| \widetilde{E}_{\mathfrak a, 0}\left(x +iy, \tfrac{1}{2} + ir \right) \right|^2 \dif r \, \dif x 
\end{align}
with $ \displaystyle h^{*}(t,r) = \frac{h(t,r)}{\tfrac{1}{4} + r^2}$. This completes the proof of the lemma.
 \end{proof}
 \end{lem4}
 \begin{lem4} \label{lem4:cuspspect}
For $s \in \CCC$ with $\re(s)>1$ and $t>0$, we have 
\begin{align*}
R_D(t,s) = \frac{s(s-1)}{2} \sum_{j =1}^{\infty} \frac{h(t,r_j)}{\lambda_j} R_{| u_j|^2}(s).
\end{align*}
\begin{proof}
Since the $u_j$'s are cusp forms they are of rapid decay at all cusps and the Rankin-Selberg transforms exist for $\re(s)>1$.
By  observation \ref{absatz4:RSEisenstein} we find
\begin{align*}
R_{|\Lambda_0(u_j)|^2}(t,s) = \left( \lambda_j + \frac{s(s-1)}{2} \right) R_{|u_j|^2} (t,s) 
\end{align*}
which implies the claim of the lemma.
\end{proof}
\end{lem4}
\subsection{Spectral and parabolic contribution}

\begin{absatz4}
Note that for an odd and squarefree positive integer $N$, we have  by \cite{hejhal}, p.~566,  proposition 6.3
\begin{align*}
\varphi_{\infty  \infty}(s) =&  2 \sqrt{\pi} \frac{\Gamma \left(s- \tfrac{1}{2} \right) \zeta(2s-1)}{\Gamma(s) \zeta(2s)} N^{-2s} \prod_{p \mid N} \frac{1}{1 - p^{-2s}} 
\end{align*}
where $\varphi_{\infty \infty}(s)$ is the function in the constant term of the $0$-th Fourier expansion of the Eisenstein series $E_{\infty, 0}(z,s)$.
(In the notation loc. cit. we have to choose $x_1=x_2=1$ and $A_1=A_2=1$, which corresponds to the cusp $\infty$). One easily computes its Laurent expansion to be
\begin{align} \label{align:expansion}
\varphi_{\infty  \infty}(s) = \frac{1}{v_N} \frac{1}{s-1} +
 \frac{1}{v_N} \left( 2 \gamma +  \frac{a\pi}{6} - 2 \sum_{p \mid N} \frac{ p^2 \log (p)}{p^2-1}   \right) + O(s-1), 
\end{align}
where $\gamma$ is the Euler constant and $a$ is the derivative of $\sqrt{\pi}  \frac{\Gamma \left(s-\tfrac{1}{2} \right)}{\Gamma(s) \zeta(2s)}$ at $s=1$.
\end{absatz4}
\begin{prop4} \label{prop4:cRSc_3}
Let $N$ be an  odd and squarefree positive integer. For any $t>0$, we have the following Laurent expansion  in a neighbourhood of $s=1$
\begin{align*} 
& R_{C_{2,1}}(t,s)- R_{C_{0,1}}(t,s) = \\
& \frac{1}{v_N} \frac{1}{s-1} + \frac{t}{ 2 v_N} +  
\frac{1}{2 v_N} \left(1+  2 \gamma +  \frac{a\pi}{6} - 2 \sum_{p \mid N} \frac{p^2 \log(p)}{p^2-1}  \right) + O(s-1). 
\end{align*} 
\begin{proof}
By lemma \ref{lem4:RSp_3}, we have
\begin{align*}
R_{C_{2,1}}(t,s)= -\frac{1}{2} h\left( t, \tfrac{i s}{2} \right) \varphi_{\infty \infty} \big( \tfrac{1+s}{2} \big)  \frac{1-s}{1+s}, 
\end{align*}
which is holomorphic at $s=1$. Then, since $h(t,\tfrac{i}{2}) = 1$, we obtain by observation \eqref{align:expansion} the Laurent expansion
\begin{align} \label{align:L23}
R_{C_{2,1}}(t,s)=  \frac{1}{2 v_N} + O(s-1).
\end{align} 
Again by lemma \ref{lem4:RSp_3}, we have
 \begin{align*}
 - R_{C_{0,1}}(t,s)= \frac{1}{2}  h \left(t, \tfrac{i s}{2} \right) \varphi_{\infty \infty} \big( \tfrac{1+s}{2} \big),
 \end{align*}
which  has a pole of order one at $s=1$. Then, since 
 \begin{align*}
\frac{1}{2} h \left( t, \tfrac{is}{2} \right) = \frac{1}{2} \exp\! \big( \!-t \left(\tfrac{1}{4} + (\tfrac{is}{2})^2\right)\big) = \frac{1}{2} + \frac{1}{4}t (s-1) + O\left( (s-1)^2 \right), 
 \end{align*} 
we obtain again  by  \eqref{align:expansion} the Laurent expansion
 \begin{align*}
  - R_{C_{0,1}}(t,s)= \frac{1}{v_N} \frac{1}{s-1} +  \frac{t}{2 v_N} + 
  \frac{1}{2 v_N} \left( 2\gamma +  \frac{a\pi}{6} - 2\sum_{p \mid N} \frac{p^2}{p^2-1} \log p \right) + O(s-1)
  \end{align*}
  which implies with
  the Laurent expansion \eqref{align:L23} the claim of the proposition.
\end{proof}
\end{prop4}
\begin{prop4} \label{prop4:cRSc_4}
For any $t>0$, we have  in a neighbourhood  of  $s=1$  the following Laurent expansion
\begin{align*}
R_{C_{2,2}}(t,s) -R_{C_{0,2}}(t,s)=  C_4(t) + O(s-1),
\end{align*}
where the  constant $C_4(t)$ tends to zero as $t \rightarrow \infty$.  
\begin{proof}
 By lemma \ref{lem4:RSp_4} we have
\begin{align*}
R_{C_{2,2}}(t,s) -R_{C_{0,2}}(t,s) = \frac{s(s-1)}{2}   \int_{0}^{+ \infty} c^{*}_{0,2}(t,y) y^{s-2} \dif y,
\end{align*}
where we set
\begin{align*}
c^{*}_{0,2}(t,y)= -  \frac{1}{4 \pi} \sum_{\frak a \in P_{\Gamma_1(N)}} \int_{-1/2}^{1/2}  \int_{- \infty}^{+\infty}  \frac{h(t,r)}{\tfrac{1}{4} +r^2}  \left| \widetilde{E}_{\frak a, 0}(x +iy, \tfrac{1}{2} + ir) \right|^2 \dif r \,  \dif x. 
\end{align*} 
By means of  \ref{absatz:rankin}  the integral
\begin{align*}
 \int_{0}^{+ \infty} c^{*}_{0,2}(t,y) y^{s-2} \dif y
\end{align*}
has a meromorphic continuation to the whole $s$-plane with a simple pole at $s=1$. Therefore 
$R_{C_{2,2}}(t,s) -R_{C_{0,2}}(t,s)$
  is holomorphic  at $s=1$ for any $t>0$.  
Since we have chosen $h(t,r) = \exp\! \big( \!-t(\tfrac{1}{4} +r^2) \big)$ the constant $C_4(t)$ in the Laurent expansion tends to $0$ as $t \rightarrow \infty$, as we claimed.
\end{proof}
\end{prop4}
\begin{prop4} \label{prop4:cdiscrete}
For any $t>0$, we have    in a neighbourhood  of  $s=1$  the following Laurent expansion
\begin{align} \label{align:constantdiskret}
R_D(t,s)= \frac{1}{2 v_N} \sum_{j =1}^{\infty} \frac{h(t,r_j)}{\lambda_j} + O(s-1).
\end{align}
\begin{proof}
We know from  \ref{absatz:rankin} that all $R_{| u_j|^2}(s)$ have a 
meromorphic continuation to the whole $s$-plane with a simple pole  at $s=1$ with residue  $v_N^{-1}$.
Hence the summands $$ \frac{s(s-1)}{2}  \frac{h(t,r_j)}{\lambda_j} R_{| u_j|^2}(s)$$ are holomorphic  at $s=1$, and since the series
\begin{align*}
R_D(t,s)= \frac{s(s-1)}{2} \sum_{j =1}^{\infty} \frac{h(t,r_j)}{\lambda_j} R_{| u_j|^2}(s)
\end{align*}
converges uniformly in a neighbourhood of $s=1$, the series  is holomorphic at $s=1$.
From this the Laurent expansion  \eqref{align:constantdiskret} follows immediately.
\end{proof}
\end{prop4}
We conclude this section by recalling some known results from \cite{abbesullmo}.
Let $N$ be a positive integer.  We denote by 
$d(N)$  the number of positive divisors of $N$ and by $\displaystyle \sigma_s(N):= \sum_{d |N} d^s$, $s \in \CCC$, the divisor sum of $N$.
\begin{prop4} \label{prop4:cRSp_1}
 For any $t>0$, we have  in a neighborhood of $s=1$  the following Laurent expansion
\begin{align*}
& R_{P_{2,1}}(t,s) - R_{P_{0,1}}(t,s) = \notag \\
 &  \Bigg(- \frac{1}{2} + \frac{1}{4 \pi} \int_{-\infty}^{+ \infty} \frac{h(t,r)}{\tfrac{1}{4}+r^2} \dif r \Bigg) \frac{\varphi(N) d(N)}{ v_N} \frac{1}{s-1}  + \\ 
 &\frac{\varphi(N) d(N)C_1(t)}{ v_N}  - 
  \frac{\varphi(N)}{v_N} \Bigg( \frac{1}{2} - \frac{1}{4 \pi} \int_{-\infty}^{+ \infty} \frac{h(t,r)}{\tfrac{1}{4}+r^2} \dif r \Bigg) \times \\
  & d(N) \Bigg(3 \gamma + \frac{a \pi}{6}  - \sum_{p | N} \frac{2p +1}{p+1} \log (p)   + \left(1 - \frac{1}{d(N)} \right) \sigma_{-1}(N) \Bigg) + O(s-1),
  \end{align*}
 where   $a$ is the derivative  of $\sqrt{\pi} \frac{\Gamma(s - \frac{1}{2})}{\Gamma(s) \zeta(2s)}$ at $s=1$ and 
  $C_1(t)$ is a function, which converges for $t \rightarrow \infty$. 
  \begin{proof}
Every element $\gamma \in \Gamma_0(N)$ with  $|\tr(\gamma) |= 2$ lies already in $\Gamma_1(N)$ as we have 
$\gamma = \left( \begin{smallmatrix} 1+a & b \\ c & 1-a \end{smallmatrix} \right)$
with $a^2 = -bc \equiv 0  \mod N$.
Hence we can just apply 
  \cite{abbesullmo}, p. 57, equation $(28)$, and  p. 60, lemma 3.3.10.
  \end{proof}
\end{prop4}
%
%
%
\begin{prop4}  \label{prop4:cRSp_2}
 For any $t>0$, we have  in a neighborhood of $s=1$  the following Laurent expansion
\begin{align*}
& R_{P_{2,2}}(t,s) -R_{P_{0,2}}(t,s) =  \\
& \left(\frac{1}{4 \pi} + C_2(t) \right) \frac{1}{s-1} + \left(\frac{\Gamma'(2) + \gamma -\log(4 \pi)}{4 \pi } +\gamma C_2(t) + C_3(t)  \right) + O(s-1),
\end{align*}
where $C_2(t)$ and $C_3(t)$ are functions in $t$ tending to $0$ as $t\rightarrow \infty$ and $\gamma$ denotes the Euler constant.
\begin{proof}
Also here, as every element $\gamma \in \Gamma_0(N)$ with  $|\tr(\gamma) |= 2$ lies already in $\Gamma_1(N)$, we can apply  \cite{abbesullmo}, p. 57, equation $(28)$  and p. 60 on the top.
\end{proof}
\end{prop4}
\section{Analytic part of $\overline{\omega}_N^2$}
In this section we determine the analytic part of the stable arithmetic self-intersection number of the relative dualizing sheaf, i.e., 
we obtain an asymptotic formula in $N$ for the Greens function
evaluated at the cusps $0$ and $\infty_d$, $d \in (\ZZZ / N \ZZZ)^{\times}$, of $\Gamma_1(N)$.
\begin{lem} \label{lem:eisensteinconst}
For an odd and  squarefree positive integer $N$, we have 
\begin{align*}
\lim_{s \rightarrow 1} \left( \varphi_{0\infty}(s) - \frac{1}{v_N} \frac{1}{s-1}\right) = 
  \frac{1}{v_N}   \left( 2 \gamma + \frac{a \pi}{6} +  \sum_{p|N} \frac{-p^2 +2p +1 }{p^2-1} \log(p)  \right),
\end{align*}
where $\gamma$ is the Euler constant and $a$ is the derivative of $\sqrt{\pi}  \frac{\Gamma \left(s-\tfrac{1}{2} \right)}{\Gamma(s) \zeta(2s)}$ at $s=1$.
\begin{proof}
From
\cite{keil}, Satz 1, p. iv, we have for any $j$
\begin{align*}
v_{0}^{j} b^{\Gamma_0(N)}_{0 \infty }(c) = \sum_{\ell=1}^{r_{\infty}} b_{0_j \infty_{\ell}}(c); 
\end{align*}
here $v_{0}^{j}$ denotes the ramification index of the cusp $0_j$ of $\Gamma_1(N)$ 
lying over the cusp $0$ of the congruence subgroup $\Gamma_0(N)$ contained in $\Gamma_1(N)$, $r_{\infty}$ refers to the number of cusps $\infty_{\ell}$ lying over the cusp
$\infty$ of $\Gamma_0(N)$, and, for $c \in \NNN$,
\begin{align*}
b^{\Gamma_0(N)}_{0 \infty}(c) = \# \left\{  \begin{pmatrix} \star & \star \\ c & \star \end{pmatrix}  \in g_{0}^{-1} \Gamma_0(N)_{0} \,  g_{0} \backslash 
g_{0}^{-1} \Gamma_0(N) g_{\infty} / g_{\infty}^{-1} \Gamma_0(N)_{\infty} \, g_{\infty} \right\}
\end{align*}
with $g_{0}= \left( \begin{smallmatrix} 0 & -1 \\  1 & 0 \end{smallmatrix} \right)$ and
 $g_{\infty}= \left( \begin{smallmatrix} 1 & 0 \\  0 & 1 \end{smallmatrix} \right)$ of $\mathrm{SL}_2(\ZZZ)$ mapping the standard cusp $\infty$ of $\mathrm{SL}_2(\ZZZ)$ to the cusps $0$ and $\infty$, respectively.
The quantity $b_{0_j \infty_{\ell} }(c)$ is defined by \eqref{align:fourcoef1}. 
\newline
In the sequel we will choose for $0_j$ the cusp $0$ of $\Gamma_1(N)$ as the cusp lying over the cusp $0$ of $\Gamma_0(N)$; 
 since we have $\Gamma_0(N)_{0} = \Gamma_1(N)_{0}$,
the cusp $0$ of $\Gamma_1(N)$ is unramified over the cusp $0$  of $\Gamma_0(N)$, which shows $v_{0}^j=1$.
Furthermore, as the group $\Gamma_1(N)$ is normal in $\Gamma_0(N)$ 
 and since again $\Gamma_0(N)_{\infty} = \Gamma_1(N)_{\infty}$,
all the cusps $\infty_{\ell}$ lying over the cusp $\infty$ of $\Gamma_0(N)$ are unramified, whence
\begin{align} \label{align:kramer}
r_{\infty} =  \lbrack \Gamma_0(N) : \Gamma_1(N) \rbrack = \varphi(N).
\end{align}
Note that we have (in the obvious notation) 
\begin{align*}
\vol_{\rm hyp}\! \left(\Gamma_1(N)  \backslash \HHH \right) = \lbrack \Gamma_0(N) : \Gamma_1(N) \rbrack  \vol_{\rm hyp}\! \left(\Gamma_0(N) \backslash \HHH \right).
\end{align*}
Then, the lemma follows immediately from lemma \ref{lem:sameinfty} and the Laurent expansion for odd and squarefree $N$ (see \cite{abbesullmo}, p. 67)
\begin{align*}
&\varphi^{\Gamma_0(N)}_{0\infty}(s) =  \frac{1}{\vol_{\rm hyp}\! \left(\Gamma_0(N) \backslash \HHH \right)} \frac{1}{s-1}  + \\
 & \frac{1}{\vol_{\rm hyp}\! \left(\Gamma_0(N) \backslash \HHH \right)} \!   \left( 2 \gamma + \frac{a \pi}{6} + \sum_{p|N} \frac{-p^2 +2p +1 }{p^2-1} \log(p)  \right)  +O(s-1).
\end{align*} 
\end{proof}
\end{lem}

\begin{thm} \label{thm:C_F} Let $N$ be an  odd  and squarefree positive integer.
Then the constant term $C_F$ in the Laurent expansion of the Rankin-Selberg transform $R_F(s)$ of $F$ at $s=1$ is given by 
\begin{align*}
C_F =&  - \frac{1}{2 g_N  v_N} \lim_{s \rightarrow 1}  \Bigg(  \frac{Z_{\Gamma_1(N)}^{\prime}(s)}{Z_{\Gamma_1(N)}(s)} - \frac{1}{s-1} \Bigg) + \\ 
&  \frac{1}{2 g_{N }  v_N} \! \Bigg( \! 2 + 2 \gamma +  \frac{a\pi}{6} - 2 \sum_{p \mid N} \frac{p^2 \log( p)}{p^2-1} \! \Bigg) +    \frac{\Gamma^{'}(2) + \gamma -\log(4 \pi)}{ 4 \pi g_N} +\\ 
&\frac{\varphi(N) d(N)}{2 g_N  v_N} \! \Bigg(\!  2 C_1 \!  - \!  3 \gamma \! - \! \frac{a \pi}{6}\! + \! \sum_{p \mid N} \frac{2p+1}{p+1} \log(p)\!  -  \! \left(\! 1 \! - \! \frac{1}{d(N)} \! \right)  \! \sigma_{-1}(N) \! \! \Bigg), 
\end{align*}
where $\gamma$ is the Euler constant, $a$ is the derivative of $\sqrt{\pi}  \frac{\Gamma \left(s-\tfrac{1}{2} \right)}{\Gamma(s) \zeta(2s)}$ at $s=1$, and 
$C_1$  is the limit $C_1:= \lim_{t \rightarrow \infty} C_1(t)$ (coming from proposition \ref{prop4:cRSp_1}).
\begin{proof}
The idea of the proof is to use equation \eqref{align:basicformula} and to determine the constant terms in the Laurent expansions at $s=1$ of the Rankin-Selberg transforms
 on the right-hand side as $t \rightarrow \infty$. This gives an expression of $C_F$, since 
 by proposition \ref{prop4:cdiscrete}, we have
 \begin{align*}
 R_D(t,s) =\frac{1}{2 v_N} \sum_{j =1}^{\infty} \frac{h(t,r_j)}{\lambda_j} + O(s-1),
 \end{align*} 
 and, with  $h(t,r)= \exp \left( -t( \tfrac{1}{4} + r^2) \right)$, we find $\lim_{t \rightarrow \infty} R_D(t,1)=0$ showing that on the left-hand side there is no contribution  from the discrete part adding to $C_F$.
\newline
Now we determine the constants in the Laurent expansions of the Rankin-Selberg transforms on the right-hand side of equation \eqref{align:basicformula} as $t \rightarrow \infty $. From 
corollary \ref{cor4:hyper} and proposition \ref{prop4:cRSc_3} we find
\begin{align} \label{align:emain21}  
& \lim_{t \rightarrow \infty}  \lim_{s \rightarrow 1} \left( R_H(t,s) + \big(R_{C_{2,1}}(t,s)- R_{C_{0,1}}(t,s) \big) -  \frac{1}{v_N} \frac{1}{s-1}  \right)= \notag \\
 &-\frac{1}{2  v_N} \lim_{s \rightarrow 1}  \left(  \frac{Z_{\Gamma_1(N)}^{\prime}(s)}{Z_{\Gamma_1(N)}(s)} - \frac{1}{s-1} \right)+ 
  \frac{1}{2 v_N} \Bigg(2+ 2 \gamma +  \frac{a\pi}{6} - 2 \sum_{p \mid N} \frac{p^2 \log (p)}{p^2-1} \Bigg).
\end{align}
Denoting by $R_1$
the residue of $R_{P_{2,1}}(t,s)- R_{P_{0,1}}(t,s)$ at $s=1$, we find 
from proposition \ref{prop4:cRSp_1} 
\begin{align} \label{align:emain22}
& \lim_{t \rightarrow \infty}  \lim_{s \rightarrow 1} \left(R_{P_{2,1}}(t,s)- R_{P_{0,1}}(t,s)   -  \frac{R_1}{s-1}  \right)=  \\ \notag
&\frac{\varphi(N) d(N)}{2 v_N} 
 \! \Bigg(\!  2 C_1 \!  - \!  3 \gamma \! - \! \frac{a \pi}{6}\! + \! \sum_{p \mid N} \frac{2p+1}{p+1} \log(p)\!  -  \! \left(\! 1 \! - \! \frac{1}{d(N)} \! \right)  \! \sigma_{-1}(N) \! \! \Bigg), 
\end{align} 
 where $C_1$ is the limit $C_1:= \lim_{t \rightarrow \infty} C_1(t)$.
From  proposition \ref{prop4:cRSp_2} we find
\begin{align} \label{align:emain23}
& \lim_{t \rightarrow \infty}  \lim_{s \rightarrow 1} \left( \! R_{P_{2,2}}(t,s) \! -R_{P_{0,2}}(t,s) -  \left(\frac{1}{4 \pi}  +  C_2(t) \right)\frac{1}{s-1}  \right)\!=  \! \frac{\Gamma'(2) + \gamma -\log(4 \pi)}{ 4 \pi},
\end{align}
and proposition \ref{prop4:cRSc_4} implies
\begin{align} \label{align:emain24}
\lim_{t \rightarrow \infty}  \lim_{s \rightarrow 1} \left(  R_{C_{2,2}}(t,s) -R_{C_{0,2}}(t,s)  \right)= 0.
\end{align}
Collecting  all constants  \eqref{align:emain21}, \eqref{align:emain22}, \eqref{align:emain23}, and \eqref{align:emain24} and 
dividing them by $g_N$,  proves the claim of the theorem.
\end{proof}
\end{thm}
\begin{cor}  \label{cor:greenbound}
Let $N$ be  an odd and squarefree positive integer  satisfying $N=11$ or $N\geq 13$. Let $X(\Gamma_1(N)) \cong X_1(N)(\CCC)$ be the modular curve with genus
$g_N > 0$. Then for the  cusps $0, \infty_d \in X_1(N)(\CCC)$, $d \in (\ZZZ / N \ZZZ)^{\times}$, we have
\begin{align*}
& g_{\rm can}(0, \infty_d) =   - \frac{2 \pi}{ g_N  v_N} \lim_{s \rightarrow 1}  \left(  \frac{Z_{\Gamma_1(N)}^{\prime}(s)}{Z_{\Gamma_1(N)}(s)} - \frac{1}{s-1} \right) +    \\
&   \frac{2 \pi}{ g_N  v_N} \Bigg( \! 2+ 2 \gamma +  \frac{a\pi}{6} - 2 \sum_{p \mid N} \frac{p^2 \log( p)}{p^2-1} \! \Bigg) +  \frac{\Gamma'(2) + \gamma -\log(4 \pi)}{ g_N }+ \\ \notag
&\frac{2 \pi \varphi(N) d(N)}{ g_N \! v_N}
 \! \Bigg(\!  2 C_1 \!  - \!  3 \gamma \! - \! \frac{a \pi}{6}\! + \! \sum_{p \mid N} \frac{2p+1}{p+1} \log(p)\!  -  \! \left(\! 1 \! - \! \frac{1}{d(N)} \! \right)  \! \sigma_{-1}(N) \! \! \Bigg) - \\
& \frac{2 \pi }{v_N}    \Bigg( 2 \gamma + \frac{a \pi}{6} + \sum_{p|N} \frac{-p^2 + 2p +1}{p^2-1}  \log(p) \Bigg) + 
O\! \left( \frac{1}{g_N} \right) \! ,
\end{align*}
where $\gamma$ is the Euler constant, $a$ is the derivative of $\sqrt{\pi}  \frac{\Gamma \left(s-\tfrac{1}{2} \right)}{\Gamma(s) \zeta(2s)}$ at $s=1$ and 
$C_1$ is the constant from theorem \ref{thm:C_F}.
\begin{proof}
By proposition \ref{prop:gcanwithrs}, we have
\begin{align*}
g_{{\rm{can}}}(0,  \infty_d)  =  4 \pi  C_F  - 2 \pi \lim_{s \rightarrow 1} \left( \varphi_{0 \infty}(s)- \frac{1}{v_N} \frac{1}{ s-1} \right) + O  \left(\frac{1}{g_N}  \right) .
\end{align*}
Hence  lemma \ref{lem:eisensteinconst} and theorem \ref{thm:C_F} imply the statement of the corollary.
\end{proof}
\end{cor}

\begin{absatz}
A bound for the constant term in the Laurent expansion of the logarithmic derivative of the Selberg zeta function at $s=1$ is provided by a result of J.~Jorgenson and J. Kramer in \cite{jorgensonkramer}, p. 29, namely
\begin{align} \label{align:selbergzeta}
\lim_{s \rightarrow 1}  \left(  \frac{Z_{\Gamma}^{'}(s)}{Z_{\Gamma}(s)} - \frac{1}{s-1} \right) = O_{\varepsilon}(N^{\varepsilon}),
\end{align}
where the implied constant depends only on $\varepsilon$.
There is the  weaker bound $O_{\varepsilon}\big(N^{\tfrac{7}{8} +\varepsilon} \big)$ due to P.~Michel and E.~Ullmo in \cite{michelullmo}, corollary 1.4.
\end{absatz}

\begin{thm} \label{thm:analyticpart}
Let $N$ be  an odd and squarefree positive integer  satisfying $N=11$ or $N\geq 13$. Let $X(\Gamma_1(N)) \cong X_1(N)(\CCC)$ be the modular curve with genus
$g_N > 0$. Then for the  cusps $0, \infty_d \in X_1(N)(\CCC)$, $d \in (\ZZZ / N \ZZZ)^{\times}$, we have
\begin{align*}
4g_N(g_N -1) g_{\mathrm{can}}(0,\infty_d) = 2 g_N + o \big(g_N \log (N) \big). 
\end{align*}
\begin{proof}
From corollary \ref{cor:greenbound}, we deduce
\begin{align} \label{align:asympgreenbound}
&4\, g_N (g_N-1)  g_{\rm can}(0, \infty_d) =   - \frac{8 \pi (g_N-1) }{v_N} \lim_{s \rightarrow 1}  \left(  \frac{Z_{\Gamma_1(N)}^{\prime}(s)}{Z_{\Gamma_1(N)}(s)} - \frac{1}{s-1} \right) + \notag \\ \notag
&   \frac{8 \pi (g_N-1)}{ v_N} \Bigg( \! 2+ 2 \gamma +  \frac{a\pi}{6} - 2 \sum_{p \mid N} \frac{p^2 \log( p)}{p^2-1} \! \Bigg) + 
 4 \left( g_N - 1 \right)  \big( \Gamma'(2)  +  \gamma  -  \log(4 \pi) \big) + \notag \\ \notag
&\frac{8 \pi \varphi(N) d(N) (g_N \! - \! 1) }{v_N}
 \! \Bigg(\!  2 C_1 \!  - \!  3 \gamma \! - \! \frac{a \pi}{6}\! + \! \sum_{p \mid N} \frac{2p+1}{p+1} \log(p)\!  -  \!\! \left(\! 1 \! - \! \frac{1}{d(N)} \! \right)  \! \sigma_{-1}(N) \! \! \Bigg) - \notag \\ 
& \frac{8 \pi g_N (g_N-1) }{v_N}    \Bigg( \! 2 \gamma + \frac{a \pi}{6} + \sum_{p|N} \frac{-p^2 +2p+1}{p^2-1}  \log(p) \Bigg) + 
O \! \left( g_N \right) \! .
\end{align}
The asymptotic 
\begin{align*} 
 \frac{24 (g_N-1)}{\prod_{p|N} (p^2 - 1)} = 1 + o(1)
\end{align*}
 together with $v_N = \frac{\pi}{6} \prod_{p|N} (p^2 - 1)$ imply
\begin{align} \label{align:hilfsasymp}
\frac{8 \pi \left( g_N -1 \right)}{ v_N} = \frac{48 \left( g_N -1 \right)}{\prod_{p|N} (p^2 - 1)} = 2 + o(1).
\end{align}
Hence, since we have
\begin{align*} 
 2 \gamma + \frac{a \pi}{6} + \sum_{p|N} \frac{-p^2 +2p+1}{p^2-1}  \log(p) =  \log(N) + O\big(\log \log(N)\big),
\end{align*}
we find with observation \eqref{align:hilfsasymp}
\begin{align*}  
\frac{8 \pi g_N (g_N-1) }{v_N}    \Bigg( \! 2 \gamma + \frac{a \pi}{6} + \sum_{p|N} \frac{-p^2 +2p+1}{p^2-1}  \log(p) \Bigg) =   2g_N  \log(N)  + o \big(g_N \log(N) \big).
\end{align*}
Noting that all other summands on the right-hand side of \eqref{align:asympgreenbound} vanish in the little $o$-term  (for the first summand we use the bound \eqref{align:selbergzeta})  we obtain
\begin{align*} 
4\, g_N (g_N-1)  g_{\rm can}(0, \infty_d) =    2g_N  \log(N)  + o \big(g_N \log(N) \big).
\end{align*}
\end{proof}
\end{thm}

\section{Geometric part of $\overline{\omega}_N^2$ and main theorem} \label{section:geometry}

\begin{absatz}
Let $X_1(N)/ \mathbb{Q}$ be the smooth projective algebraic curve over $\mathbb{Q}$ that classifies elliptic curves equipped with a point of exact order $N$. 
There exists a canonical analytic isomorphism
\begin{align*}
j:  X\big( \Gamma_1(N) \big) = \Gamma_1(N) \backslash \big(\mathbb{H} \cup \mathbb{P}^1_{\mathbb{Q}} \big) \longrightarrow X_1(N)(\mathbb{C}),
\end{align*} 
and we say that a  $K$-rational point $x \in X_1(N)(K)$ is a cusp if $x \in j( \Gamma_1(N) \backslash \mathbb{P}^1_{\mathbb{Q}})$. Note that the cusp $0$ of  $X_1(N)/ \mathbb{Q}$  is $\QQQ$-rational and
the cusp $\infty$ of  $X_1(N)/ \mathbb{Q}$ is $\QQQ(\zeta_N)$-rational (see \cite{ogg}, proposition 1).   
\end{absatz}

\begin{absatz}
Let $X_1(N)/\mathbb{Q}(\zeta_N) = X_1(N) \times_{\mathbb{Q}} \mathbb{Q}(\zeta_N)$ be the modular curve over the cyclotomic field $ \mathbb{Q}(\zeta_N)$.
 Let   $\mathcal X_1(N)/\mathbb{Z} \lbrack \zeta_N \rbrack$ be the minimal  regular model of $X_1(N)/\mathbb{Q}(\zeta_N)$, i.e.,
   a regular, projective, and flat $\mathbb{Z} \lbrack \zeta_N \rbrack$-scheme with
generic fiber isomorphic to $X_1(N) /  \mathbb{Q}(\zeta_N)$. Under the assumption that $g_N \geq1$,  minimality means by Castelnuovo's criterion that the canonical divisor 
 $K_N$ corresponding to relative dualizing sheaf 
 is numerically
effective, i.e.,  $K_N \cdot V  \geq 0$ holds for every vertical prime divisor $V$ of $\mathcal X_1(N) / \mathbb{Z} \lbrack \zeta_N \rbrack$.
\newline
Integral models of  modular curves were intensively studied by many people. We collect in the following proposition some
 facts from  \cite{katzmazur}  we need in the sequel.
\end{absatz}

\begin{prop} \label{prop:katzmazur}
Let $N$ be a squarefree positive integer of the form $N= N^{'}qr$  with $q$ and $r$ two relative prime integers  satisfying $q,r \geq 4$. 
The minimal regular model $\mathcal X_1(N) / \mathbb{Z}\lbrack \zeta_N \rbrack$ has smooth fibers over prime ideals $\mathfrak p$ of  $ \mathbb{Z}\lbrack \zeta_N \rbrack$ with $\mathfrak p \nmid N$. 
For  $\mathfrak p | N$ the fiber of $ \mathcal X_1(N) / \mathbb{Z}\lbrack \zeta_N \rbrack$ over $\mathfrak p$ is the union
of two irreducible,  smooth, and proper $k(\mathfrak p)$-curves $C_{1,\mathfrak p}$ and $C_{2,\mathfrak p}$, $k(\mathfrak p)$ the residue field at $\mathfrak p$, intersecting transversally in
\begin{align} \label{align:ssp}
s_{\mathfrak p}:= \frac{p-1}{24} \cdot  \frac{\varphi(N/p) N}{p} \prod_{q \mid \frac{N}{p}} \left(1  + \frac{1}{q} \right)
\end{align} 
$k(\mathfrak p)$-rational points. Moreover, the curves $C_{1, \mathfrak p}$ and $C_{2, \mathfrak p}$ are isomorphic.
\begin{proof}
First we suppose that there is a regular model of $X_1(N)/\mathbb{Q}(\zeta_N)$ with fibers as described in the proposition. Then, the adjunction formula for arithmetic surfaces immediately implies 
$K_N \cdot V  \geq 0$  for every vertical prime divisor $V$, hence the regular model is minimal. As a minimal regular model is unique up to isomorphism, it suffices to find such a regular model. 
\newline
Let $\overline{\mathfrak M}\big(\Gamma_1(N) \big)/ \mathbb{Z}\lbrack \zeta_N \rbrack$ be the  compactified coarse moduli scheme of canonical balanced $\Gamma_1(N)$-structures as described in \cite{katzmazur}, chap. 9,
which exists only if $N$ is a  positive integer of the form $N= N^{'}qr$ with $q$ and $r$ two relative  prime integers satisfying $q,r \geq 4$. 
It follows from the modular interpretation that $\mathfrak M(\Gamma_1(N))/ \mathbb{Z}\lbrack \zeta_N \rbrack$ is a model
for $X_1(N)/ \mathbb{Q}(\zeta_N)$. 
That $\overline{\mathfrak M}\big(\Gamma_1(N) \big)/ \mathbb{Z}\lbrack \zeta_N \rbrack$ is a regular model having smooth fibers over $\mathfrak p \nmid N$ 
follows from \cite{katzmazur}, theorem 5.5.1, theorem 10.9.1, and the summarizing table on p. 305. For  $\mathfrak p | N$ the fiber
 of $ \mathcal X_1(N) / \mathbb{Z}\lbrack \zeta_N \rbrack$ over $\mathfrak p$ is the union
of two irreducible,  smooth, and proper $k(\mathfrak p)$-curves which are isomorphic and intersect transversally, follows from \cite{katzmazur}, theorem 13.11.4.
 The formula for the  intersection number $s_{\mathfrak p}$ follows
from \cite{katzmazur}, corollaries 5.5.3 and 12.9.4.
\end{proof}
\end{prop}
\begin{absatz} \label{absatz:horizontaldiv}
Let $0, \infty \in X_1(N)\big(\mathbb{Q}(\zeta_N)\big)$ be the cusps with representatives $(0:1)$, $(1:0)$ in $\mathbb{P}^1_{\mathbb{Q}}$, respectively. We let $H_0, H_{\infty}$ be the horizontal divisors obtained 
by taking the Zariski closure of $0, \infty$ in $\mathcal X_1(N) / \mathbb{Z}\lbrack \zeta_N \rbrack$, respectively. Note that for $m=0, \infty$ there exists an open subschemes containing $H_m$, which is smooth
over $\mathbb{Z}\lbrack \zeta_N \rbrack$ (see \cite{katzmazur}, theorem 10.9.1).
\end{absatz}

\begin{prop} \label{prop:falhirassumption}
Let $N$ be a squarefree positive integer of the form $N= N^{'}qr$ with $q$ and $r$ two relative prime integers satisfying $q,r \geq 4$. 
Then, there exist vertical divisors $V_{m} \in \mathrm{Div}(\mathcal{X}_1(N))_{\mathbb{Q}}$ $(m=0,\infty)$ satisfying  
\begin{align} \label{align:perpend}
\left( \overline{\omega}_{\mathcal X_1(N)/\mathbb{Z}\lbrack \zeta_N \rbrack} \otimes \overline{\mathcal O} _{\mathcal{X}_1(N)}(H_m)^{\otimes -(2g_N-2)} 
\otimes \overline{\mathcal O}_{\mathcal{X}_1(N)}(V_{m}) \right) \cdot
 \overline{\mathcal O}_{\mathcal X_1(N)}(V) =0 
\end{align}
 for all vertical divisors $V$ of $\mathcal{X}_1(N) / \mathbb{Z}\lbrack \zeta_N \rbrack$.
  Furthermore, we have the following intersection
 numbers
\begin{align} \displaystyle
 (V_{0},V_{0})_{\rm{fin}}  =  (V_{\infty},V_{\infty})_{\rm{fin}} =
-  (V_{0},V_{\infty})_{\rm{fin}}   = - \varphi(N) \frac{24(g_N-1)^2}{\prod_{p| N} (p^2 -1)}  \sum_{p|N} \frac{p+1}{p-1} \log(p).
\end{align}
\begin{proof}
We start by considering a fiber over a closed point $\mathfrak p \in \Spec(\ZZZ\lbrack \zeta_N\rbrack)$ with $\mathfrak p | N$.
 The fiber consists by proposition \ref{prop:katzmazur}  of the two irreducible components $C_{1,\mathfrak p}$ and $C_{2,\mathfrak p}$.
 The horizontal divisors $H_0$ and $H_{\infty}$ intersect the fiber in a smooth $k(\mathfrak p)$-rational point, and,
by the cusp and component labeling in \cite{katzmazur}, p. 296, they do not intersect the same component. We denote
by $C_{0,\mathfrak p}$ and $C_{\infty,\mathfrak p}$ the component intersected by $H_0$ and $H_{\infty}$, respectively. Let us define the $\QQQ$-divisors 
\begin{align*}
V_{0, \mathfrak p} := - \frac{g_N-1}{ s_{\mathfrak p}} C_{0,\mathfrak p}
\hspace{1cm} \text{and} \hspace{1cm}
V_{\infty, \mathfrak p} := - \frac{g_N-1}{ s_{\mathfrak p}} C_{\infty,\mathfrak p}.
\end{align*}
 Then, we claim  that
\begin{align*}
V_{0} := \sum_{\substack{ \mathfrak p \in \Spec(\ZZZ \lbrack \zeta_N \rbrack) \\ \mathfrak p | N}}
 V_{0, \mathfrak p} \hspace{0,5cm} \text{ and } \hspace{0,5cm} V_{\infty}:= \sum_{\substack{ \mathfrak p \in \Spec(\ZZZ \lbrack \zeta_N \rbrack) \\ \mathfrak p | N}} V_{\infty, \mathfrak p}
\end{align*}
fulfill the conditions stated in \eqref{align:perpend}. Noting that
\begin{align*}
K_N  \cdot C_{0,\mathfrak p} = \deg \omega_{C_{0,\mathfrak p}/ k(\mathfrak p)} = \deg \omega_{C_{\infty,\mathfrak p}/ k(\mathfrak p)} 
= K_N \cdot C_{\infty,\mathfrak p}
\end{align*}
and $K_N \cdot (C_{0,\mathfrak p} + C_{\infty,\mathfrak p}) = 2g_N -2$, we find $K_N  \cdot C_{0,\mathfrak p}= K_N \cdot C_{\infty,\mathfrak p} = g_N-1$. We have to consider the following three cases:
\begin{enumerate}
\item[$(i)$] $V = V_{\mathfrak p}$  with $\mathfrak p \nmid N$. In this case we calculate, using the adjunction formula
\begin{align*}
&\left( \overline{\omega}_{\mathcal X_1(N)/\mathbb{Z}\lbrack \zeta_N \rbrack} \otimes \overline{\mathcal O} _{\mathcal{X}_1(N)}(H_m)^{\otimes -(2g_N-2)} 
\otimes \overline{\mathcal O}_{\mathcal{X}_1(N)}(V_{m}) \right) \cdot
 \overline{\mathcal O}_{\mathcal X_1(N)}(V) = \\
 &\big( 2g_N-2 - (2g_N-2)  \big) \log\big(\sharp k(\mathfrak p) \big) =0.
\end{align*}
\item[$(ii)$] $V =  C_{m,\mathfrak p}$ with $\mathfrak p | N$. In this case we calculate
\begin{align*}
&\left( \overline{\omega}_{\mathcal X_1(N)/\mathbb{Z}\lbrack \zeta_N \rbrack} \otimes \overline{\mathcal O} _{\mathcal{X}_1(N)}(H_m)^{\otimes -(2g_N-2)} 
\otimes \overline{\mathcal O}_{\mathcal{X}_1(N)}(V_{m}) \right) \cdot
 \overline{\mathcal O}_{\mathcal X_1(N)}(V) = \\
 &\big( g_N-1 - (2g_N-2) + g_N-1 \big) \log\big(\sharp k(\mathfrak p)\big) =0.
\end{align*}
\item[$(iii)$] $V \not= C_{m,\mathfrak p}$ with $\mathfrak p | N$. In this case we calculate
\begin{align*}
&\left( \overline{\omega}_{\mathcal X_1(N)/\mathbb{Z}\lbrack \zeta_N \rbrack} \otimes \overline{\mathcal O} _{\mathcal{X}_1(N)}(H_m)^{\otimes -(2g_N-2)} 
\otimes \overline{\mathcal O}_{\mathcal{X}_1(N)}(V_{m}) \right) \cdot
 \overline{\mathcal O}_{\mathcal X_1(N)}(V) = \\
 &\big( g_N-1 -(g_N-1) \big) \log \big( \sharp k(\mathfrak p)\big) =0.
\end{align*}
\end{enumerate}
This proves the first part of the proposition. Now,  proposition \ref{prop:katzmazur} implies
\begin{align*} 
&(V_{0},V_{0})_{\rm{fin}} = (V_{\infty},V_{\infty})_{\rm{fin}} =
- (V_{0},V_{\infty})_{\rm{fin}} =-\sum_{\mathfrak p | N} \frac{(g_N-1)^2}{ s^2_{\mathfrak p}} C_{0,\mathfrak p} \cdot C_{\infty,\mathfrak p}  = \\ 
 & -\sum_{\mathfrak p | N}
 \frac{(g_N-1)^2 \log\big(\sharp k(\mathfrak p) \big)}{ s_{\mathfrak p} } 
  = - 24  (g_N-1)^2\sum_{p|N} \frac{1}{p-1} \sum_{\mathfrak p| p} \frac{\log\big(\sharp k(\mathfrak p)\big)}{\prod_{q| N/p} \left(q^2-1 \right)}.
\end{align*}
Noting that
$ \displaystyle \sum_{\mathfrak p| p} \log\big(\sharp k(\mathfrak p)\big) = \varphi(N/p) \log (p)$,
we obtain 
\begin{align*}
(V_{0},V_{0})_{\rm{fin}} = (V_{\infty},V_{\infty})_{\rm{fin}} =
- (V_{0},V_{\infty})_{\rm{fin}} =
 - \varphi(N) \frac{24(g_N-1)^2}{\prod_{p| N} (p^2 -1)}  \sum_{p|N} \frac{p+1}{p-1} \log(p),
\end{align*}
which proves the proposition.
\end{proof}
\end{prop}

\begin{prop} \label{prop:geometricpart}
Let $N$ be a squarefree positive integer of the form $N= N^{'}qr$ with $q$ and $r$ two relative prime integers satisfying $q,r \geq 4$.
Let $0, \infty$ be the cusps of $X\big( \Gamma_1(N) \big)$ with representatives $(0:1)$, $(1:0)$ in $\mathbb{P}^1_{\mathbb{Q}}$, respectively, and let $V_0, V_{\infty}$ be the vertical
divisors of proposition \ref{prop:falhirassumption}. Then, we have
\begin{align} \label{align:geombound}
\overline{\omega}^{2}_N =   4 g_{N} (g_N -1)  g_{{\rm{can}}}(0, \infty) + 
\frac{1}{\varphi(N)} \frac{g_N+1}{g_N - 1}  \left( V_{0}, V_{\infty} \right)_{\rm{fin}}.
\end{align}
\begin{proof}
Formula \eqref{align:geombound} is analogous  to the corresponding formula in proposition D of \cite{abbesullmo} and, in the smooth case, the formula
 given in \cite{szpiro}, p. 241. 
 \newline
We start by noting that a multiple of the line bundle
 \begin{align*}
 \Big( \omega_{\mathcal{X}_1(N) / \ZZZ\lbrack \zeta_N \rbrack} \otimes
  \mathcal O_{\mathcal{X}_1(N)}(H_0)^{\otimes -(2g_N-2)} \otimes \mathcal O_{\mathcal{X}_1(N)}(V_{0})\Big) \big|_{X_1(N)} \in J_1(N)\big(\QQQ(\zeta_N)\big)
 \end{align*}
 has support in the cusps.
  Hence a well-known theorem of Manin and Drinfeld (see \cite{drinfeld}) says, that it is a torsion element in the Jacobian $J_1(N)/ \QQQ\big(\zeta_N\big)$.
  Now, as condition \eqref{align:perpend} is satisfied, a theorem of Faltings and Hriljac  (see  \cite{faltings},  theorem 4), and the fact that the 
    N${\rm \acute{e}}$ron-Tate height vanishes at torsion points, imply as in \cite{abbesullmo}
\begin{align} \label{align:horiz}
 \overline{\omega}^2_{\mathcal X_1(N)/\ZZZ\lbrack \zeta_N \rbrack}  = &- 2g_N(g_N-1) 
 \left( \overline{\mathcal O}_{\mathcal{X}_1(N)}(H_0)^2+\overline{\mathcal O}_{\mathcal{X}_1(N)}(H_{\infty})^2 \right) + \notag \\
&\frac{1}{2} \left(\overline{\mathcal O}_{\mathcal{X}_1(N)}(V_{0})^2 + \overline{\mathcal O}_{\mathcal{X}_1(N)}(V_{\infty})^2\right).
\end{align}
Now we consider the admissible metrized line bundle 
\begin{align*}
 \overline{\mathcal O}_{\mathcal{X}_1(N)}(H_{\infty}) \otimes \overline{\mathcal O}_{\mathcal{X}_1(N)}(H_0)^{\otimes-1} \otimes
\left( \overline{\mathcal O}_{\mathcal{X}_1(N)}(H_0) \otimes  \overline{\mathcal O}_{\mathcal{X}_1(N)}(H_0)^{\otimes -1} \right)^{\otimes1/(2g_N-2)},
\end{align*}
which is orthogonal to all vertical divisors $V$ of $\mathcal X_1(N)$ because of the conditions \eqref{align:perpend} and 
has the generic fiber with support in the cusps. Then, a similar argument as above shows
 \begin{align} \label{align:horiz1}
\overline{\mathcal O}_{\mathcal{X}_1(N)}(H_0)^2+\overline{\mathcal O}_{\mathcal{X}_1(N)}(H_{\infty})^2  =& 
 2 \overline{\mathcal O}_{\mathcal{X}_1(N)}(H_0)  \cdot  \overline{\mathcal O}_{\mathcal{X}_1(N)}(H_{\infty})+ \notag \\ 
& \frac{(V_{0},V_{\infty})_{\rm{fin}}  -2\left(V_{0}, V_{\infty} \right)_{\rm{fin}} \! + (V_{\infty},V_{\infty})_{\rm{fin}} }{(2g_N-2)^2} . 
\end{align}
Since by \ref{absatz:horizontaldiv} the horizontal divisors $H_0$ and $H_{\infty}$ do not intersect, lemma \ref{prop:falhirassumption}, 
substituting equation \eqref{align:horiz1} into \eqref{align:horiz}, implies
\begin{align} \label{align:fundamentalb}
 \overline{\omega}_{\mathcal X_1(N)/\ZZZ\lbrack \zeta_N \rbrack}^2 =
   4  g_N (g_N -1) \sum_{\sigma : \QQQ(\zeta_N) \rightarrow \CCC } g^{\sigma}_{{\rm{can}}}\big(0, \infty_{\sigma} \big) + \frac{g_N+1}{g_N - 1}  \left( V_{0}, V_{\infty} \right)_{\rm{fin}} 
\end{align}
Note that   the modular curve $X_1(N)/\QQQ$ is defined over $\QQQ$, and hence 
 the Riemann surfaces $X_1(N)_{\sigma}$ are in fact all equal to $X(\Gamma_1(N))$. Therefore,  dividing both sides of equality \eqref{align:fundamentalb}  by $\varphi(N)$, it follows from proposition \ref{prop:gcanwithrs}
noting that $\infty_{\sigma}= \infty_d$ for some $d \in \left( \ZZZ / N \ZZZ \right)^{\times}$  that
 \begin{align*} 
\overline{\omega}_N^2 = 4 g_N (g_N-1)  g_{{\rm{can}}}(0, \infty) + 
\frac{1}{\varphi(N)} \frac{g_N+1}{g_N - 1}  \left( V_{0}, V_{\infty} \right)_{\rm{fin}} .
\end{align*}
\end{proof}
\end{prop}
\begin{thm} \label{thm:main}
Let $N$ be an odd and squarefree positive integer of the form $N= N^{'}qr$ with $q$ and $r$ two relative prime integers satisfying $q,r \geq 4$. Then, we have
\begin{align*}
\overline{\omega}_N^2 = 3 g_N \log(N) + o\big(g_N \log(N) \big).
\end{align*}
\begin{proof}
By proposition \ref{prop:geometricpart}, using  theorem \ref{thm:analyticpart} and proposition \ref{prop:falhirassumption}, we have
\begin{align*} 
\overline{\omega}_N^2 = 2 g_N  + 
 \frac{24(g_N+1)(g_N-1)}{\prod_{p|N} (p^2 - 1)} \sum_{p|N} \frac{p+1}{p-1} \log(p) + o\big(g_N \log(N) \big).
\end{align*}
Noting that
\begin{align*} 
 \frac{24 (g_N-1)}{\prod_{p|N} (p^2 - 1)} = 1 + o(1),
\end{align*}
 the claimed asymptotic follows.
\end{proof}
\end{thm}


\section{Arithmetic applications}

\begin{absatz} {\it{Stable Faltings height}}. 
Let $J_1(N)/ \QQQ$ be the Jacobian variety of the modular curve $X_1(N)/ \QQQ$ and let $ h_{\rm{Fal}} \big(J_1(N) \big)$ be the stable Faltings height of $J_1(N)/ \QQQ$.
The arithmetic Noether formula (see  \cite{moret-bailly}, theorem 2.5)
implies
\begin{align} \label{align:noether}
 12  h_{\rm{Fal}} \big(J_1(N) \big) =  \overline{\omega}^2_{N}   + \sum_{p|N} \frac{s_p}{p-1} \log (p)  + \delta_{\rm{Fal}} \big(X\big(\Gamma_1(N)\big) \big) - 4g_{N} \log (2 \pi),
\end{align}
where $s_p$ is given by the formula \eqref{align:ssp} and $\delta_{\rm{Fal}}\big(X\big(\Gamma_1(N)\big)\big)$
 denotes the Faltings's delta invariant of $ X\big(\Gamma_1(N)\big)$
 (for the definition see  \cite{faltings}, theorem 1,  or, for another approach due to J.-B. Bost, see \cite{soule}).  In
 \cite{jorgensonkramer1} it is proved (see loc. cit.  theorem 5.3 and remark 5.8) that
\begin{align} \label{align:faltingsdelta}
\delta_{\rm{Fal}}\big(X\big(\Gamma_1(N)\big)\big) = O(g_N).
\end{align}

\begin{thm} \label{thm:falheight}
Let $N$ be an odd and squarefree positive integer of the form $N= N^{'}qr$ with $q$ and $r$ two relative prime integers satisfying $q,r \geq 4$.
Then, we have
\begin{align*}
h_{\rm{Fal}}\big(J_1(N) \big) = \frac{g_{N}}{4} \log (N) + o\big(g_{N} \log (N) \big).
\end{align*}
\begin{proof}
Noting that
\begin{align*}
 \frac{1}{3} \sum_{p | N} \frac{s_p}{ p-1} \log( p) = \frac{g_N}{3} \sum_{p|N} \frac{\log(p)}{p^2-1} + \frac{\varphi(N) d(N)}{12} \sum_{p|N} \frac{\log(p)}{p^2-1} + O\big(\log(N)\big)
\end{align*}
and 
\begin{align*}
 \frac{g_N}{3} \sum_{p|N} \frac{\log(p)}{p^2-1} + \frac{\varphi(N) d(N)}{12} \sum_{p|N} \frac{\log(p)}{p^2-1} = o\big(g_N \log(N)\big),
\end{align*}
we find 
\begin{align} \label{align:asympmichel}
 \frac{1}{3} \sum_{p | N} \frac{s_p}{ p-1} \log( p) = o\big(g_N \log(N)\big).
\end{align}
Hence observations \eqref{align:noether} and \eqref{align:faltingsdelta} together with theorem \ref{thm:main} imply statement of the theorem.
\end{proof}
\end{thm}
\begin{rem} \label{rem:edixhoven}
In \cite{edixhoven}, p. 83, theorem 16.7,  the authors found the bound $$h_{\rm Fal}\big(J_1(pl) \big) = O\big((pl)^2 \log(pl)\big),$$ 
where $p$ and $l$ are  distinct primes with $l>5$. In theorem \ref{thm:falheight} the leading term is explicit.
\newline
We note that theorem \ref{thm:falheight} might lead to a
similar bound on the minimal number of congruences of modular forms  with respect to $\Gamma_1(N)$  as in \cite{jorgensonkramer},  remark 6.6, p. 37, in the case of  $\Gamma_0(N)$.
\end{rem}
\end{absatz}

\begin{absatz}{\it{Admissible self-intersection number}}.
We assume that the reader is familiar with the theory of the admissible pairing in \cite{zhang}.

 Let $N$ be  an odd and squarefree integer  of the form $N=N^{'} qr$ with $q$ and $r$ relative prime integeres satisfying  $ q,r \geq 4$. 
The dual reduction graph
  $G_{\mathfrak p } $ of the fiber of $\mathcal X_1(N)/ \ZZZ \lbrack \zeta_N \rbrack$ over the prime ideal $\mathfrak p $ of $\ZZZ \lbrack \zeta_N \rbrack$ with $\mathfrak p \nmid N$
 consists by proposition \ref{prop:katzmazur} of two vertices which are connected by
 $s_{\mathfrak p}=\frac{p-1}{24} \cdot  \frac{\varphi(N/p) N}{p} \prod_{q \mid \frac{N}{p}} \left(1  + \frac{1}{q} \right)$ edges  of length $1$.
Hereby the two vertices correpond to the irreducible components $C_{0, \mathfrak p}$ and $C_{\infty, \mathfrak p}$ over $\mathfrak p$.
Note that the genera of the two components $C_{0, \mathfrak p}$ and $C_{\infty, \mathfrak p}$ are the same by proposition \ref{prop:katzmazur}, which we denote by $g_{\mathfrak p}$. Then, we have
$g_N = 2 g_{\mathfrak p} + s_{\mathfrak p} -1$ and
the    canonical divisor 
$  K_{G_{\mathfrak p}}$ on $G_{\mathfrak p}$ (for the definition see \cite{zhang}, p. 175) is in this case given by
\begin{align*}
 K_{G_{\mathfrak p}}= (g_N -1) \lbrack 0 + \infty \rbrack. 
\end{align*}
 We set $a_{\mathfrak p} :=\frac{s_{\mathfrak p}}{ s_{\mathfrak p} -1} g_{\mathfrak p}$ and $l_{\mathfrak p} :=2a_{\mathfrak p} + s_{\mathfrak p}$, and the admissible  measure
with respect to $  K_{G_{\mathfrak p}}$ on  $G_{\mathfrak p}$ (for the definition see \cite{zhang}, theorem 3.2) is given by
\begin{align} \label{align:admmeasure}
\mu_{\mathfrak p} := \frac{a_{\mathfrak p}}{l_{\mathfrak p}} \big( \delta_0 + \delta_{\infty} \big) + \frac{1}{l_{\mathfrak p}} \dif x.
\end{align}
 Recalling the notation from   \cite{abbesullmo}, we have for
   $n = s_{\mathfrak p}$ that the admissible Green's function $g^{G_{\mathfrak p}}_{\mu_{\mathfrak p}}$ with respect to $\mu_{\mathfrak p}$ of the graph  $G_{\mathfrak p}$ is given by
 (cf. \cite{abbesullmo}, p. 65)
\begin{align} \label{align:admgreen}
g^{G_{\mathfrak p}}_{\mu_{\mathfrak p} }(x_j,0)&=\frac{1}{2l_{\mathfrak p}} x_j^2 - \frac{a_{\mathfrak p}+s_{\mathfrak p}}{s_{\mathfrak p}l_{\mathfrak p}}x_j - \frac{3a_{\mathfrak p}+2s_{\mathfrak p}}{6s_{\mathfrak p} l_{\mathfrak p}} ,\notag \\
g^{G_{\mathfrak p}}_{\mu_{\mathfrak p} }(x_j, \infty)&= \frac{1}{2 l_{\mathfrak p}} (1-x_j)^2 - \frac{a_{\mathfrak p}+s_{\mathfrak p}}{s_{\mathfrak p} l_{\mathfrak p}}(1- x_j) - \frac{3a_{\mathfrak p}+2s_{\mathfrak p}}{6s_{\mathfrak p} l_{\mathfrak p}}, \notag \\
g^{G_{\mathfrak p}}_{\mu_{\mathfrak p} }(x_j, x_j)&= \left(1- \frac{1}{s_{\mathfrak p}} -\frac{1}{l_{\mathfrak p}} \right) x_j(1-x_j) - \frac{3a_{\mathfrak p}+2s_{\mathfrak p}}{6s_{\mathfrak p} l_{\mathfrak p}},
\end{align}
where $0 \leq x_j \leq 1$ ($j= 1 , \ldots, s_{\mathfrak p}$) denotes the coordinate along an edge.
Note that for prime ideals $\mathfrak p , \mathfrak p'$ of $\ZZZ\lbrack \zeta_N \rbrack$ with $\mathfrak p | p $ and $\mathfrak p' | p$, we have 
$G_{\mathfrak p} = G_{\mathfrak p'}$,
 and we simply write  $G_p$ for this graph;  further, we write  in this case also simply $ s_p,  a_p, l_p, K_{G_p},  \mu_p$ for  
 $ s_{\mathfrak p}, a_{\mathfrak p}, l_{\mathfrak p}, K_{G_{\mathfrak p}}, \mu_{\mathfrak p} $, respectively. Let $\overline{\omega}_{a,N}$ be the admissible metrized relative dualizing sheaf of the curve
$X_1(N)/ \QQQ$ as defined in \cite{zhang}, p. 181 and p. 188.
\begin{thm} \label{thm:lower}
 Let $N$ be  an odd and squarefree integer  of the form $N=N^{'} qr$ with $q$ and $r$ relative prime integeres satisfying  $ q,r \geq 4$. 
Then, we have 
\begin{align*}
 \overline{\omega}^{2}_{a,N}= 3 g_N  \log(N) + o \big(g_N  \log(N) \big).
\end{align*}
\begin{proof}
By \cite{zhang}, theorem 5.5, we have
\begin{align*}
 \overline{\omega}^{2}_{N}  -  \overline{\omega}^{2}_{a,N}  
 = \frac{1}{\varphi(N)} \sum_{\mathfrak p | N} r_{\mathfrak p} \log\big( \sharp k(\mathfrak p) \big) =   \sum_{p | N} \frac{r_{p}}{p-1} \log( p),
\end{align*}
where
\begin{align*}
r_{\mathfrak p}:= r_p:= \int_{G_p} g^{G_p}_{\mu_p}(x,x) \big( (2g_N -2)\mu_p(x) + \delta_{K_{G_p}} \big).
\end{align*}
 Now we calculate (as in \cite{abbesullmo}, p. 65) the value $r_p$ for the graph $G_p$ with $p | N$. The definition of $r_p$ and  the formula for the admissible measure $\mu_p$ of \eqref{align:admmeasure} yield
\begin{align*}
r_p = - \frac{(g_N-1)(3g_N+ s_p -1)}{3s_pg_N} +  \frac{(g_N-1)^2}{3 g_N^2} s_p - \frac{(2s_p-1) (g_N-1)^2}{s_pg_N^2}.
\end{align*}
Noting that $s_p <g_N $ and 
$\frac{g_N}{s_p} = p+1 + O(1)
$
holds for $p | N$, we obtain
\begin{align*}
r_p = - (p+1) +  \frac{s_p}{3} + O(1).
\end{align*}
Therefore, we have
\begin{align*}
  \sum_{p | N} \frac{r_{p}}{p-1} \log( p) 
=   \sum_{p | N} \frac{s_p}{ 3(p-1)} \log( p)  + O\big( \log(N)\big). 
\end{align*}
Then, theorem \ref{thm:main} and the asymptotic \eqref{align:asympmichel} imply
\begin{align*}
\overline{\omega}_{a,N}^2 = 3 g_N \log(N) + o\big(g_N \log(N)\big).
\end{align*}
This completes the proof of the theorem.
\end{proof}
\end{thm}
\end{absatz}

\begin{absatz} {\it{Effective Bogomolov}}.
Let $X/K$ be a smooth projective and geometrically connected curve over a number field $K$ of genus $g_X>1$.  
For a divisor $D \in {\mathrm{Div}}(X)$ of degree $1$, let 
$$
\varphi_D : X \longrightarrow \jac(X)
$$
be the embedding of the curve $X/K$ into its Jacobian $\jac(X)/K$ defined by the mapping $ x \mapsto \lbrack \mathcal O_X( x - D) \rbrack$.
Then
there exists an $\varepsilon >0$ such that the set of algebraic points
$$
\left \{ x \in X(\overline{K})\, \big| \, h_{\rm{NT}} \big( \varphi_{D}(x) \big) \leq \varepsilon \right \}
$$ 
is finite, where $h_{\rm{NT}}$ denotes the N$\mathrm{\acute{e}}$ron-Tate height on  $\jac(X)/K$ (Bogomolov's conjecture).
The conjecture was proved by E. Ullmo   in the case of curves and by S.-W. Zhang, more general, for any non-torsion subvariety $X$ of an abelian variety $A/K$.

Due to L. Szpiro and S.-W. Zhang it is known  (see \cite{zhang}, theorem 5.6) that the  set of algebraic points
\begin{align} \label{align:bogomolov}
\left \{ x \in X(\overline{K}) \, \big| \, h_{\rm{NT}} \big( \varphi_{D}(x) \big) < \frac{\overline{\omega}_{a}^2}{4(g_X-1)} \right \}
\end{align}
is finite, and that the set of algebraic points
\begin{align} \label{align:bogomolov1}
\left \{ x \in X(\overline{K}) \, \big| \, h_{\rm{NT}} \big( \varphi_{D}(x) \big) \leq \frac{\overline{\omega}^2_a }{2(g_X-1)} \right \}
\end{align}
is infinite if  $\lbrack \mathcal O_X \big( K_{X} - (2g_{X} -2)D \big) \rbrack$ is a torsion element in $\jac(X)/K$.

\begin{thm} \label{thm:bogomolov}
Let $N$ be  an odd and squarefree positive integer  of the form $N=N^{'} qr$ with $q$ and $r$ relative prime integers satisfying $ q,r \geq 4$. 
Then, 
for any $\varepsilon > 0$, there is   a  sufficiently large $N$ such that the set of algebraic points
\begin{align*}
\left \{ x \in X_1(N)(\overline{\QQQ})\, \big| \, h_{\rm{NT}}\big( \varphi_D(x) \big) < \left(\frac{3}{4}- \varepsilon \right) \log( N) \right \}
\end{align*}
is finite, and  
the set of algebraic points
\begin{align*}
\left \{ x \in X_1(N)(\overline{\QQQ}) \, \big| \, h_{\rm{NT}} \big( \varphi_D(x) \big) \leq \left( \frac{3}{2} + \varepsilon \right) \log (N) \right \} 
\end{align*}
is infinite, if $\lbrack \mathcal O_{X_1(N)} \big(K_{X_1(N)} - (2g_N -2)D \big) \rbrack$ is a torsion element in $J_1(N)/\QQQ$.
\begin{proof}
The first and second statement of the theorem follow  immediately from the height bounds    \eqref{align:bogomolov} and 
\eqref{align:bogomolov1} in conjunction with theorem \ref{thm:lower}.
\end{proof}
\end{thm}
\begin{rem}
If $D=\lbrack 0 \rbrack$ then $\lbrack \mathcal O_{X_1(N)} \big(K_{X_1(N)} - (2g_N -2)D \big) \rbrack$ 
is a torsion element in $J_1(N)/ \QQQ$, which
gives an example for the second statement in theorem \ref{thm:bogomolov}.
\end{rem}
\end{absatz}

\section{Appendix: Meromorphic Continuation}

In this appendix we give a proof of proposition \ref{prop4:continuation}. Recall that for
$N$   an odd positive integer
and  $l \in \ZZZ$  with $|l|>2$ and 
$l \equiv 2  \mod N$, as well as
 $0< u < N$ with $(u,N)=1$,  we define the zeta function
\begin{align*} 
\zeta_{u,N}(s,l)=   \sum_{(\gamma, (m,n)) \in \Delta^{u+}_l(N)  / \Gamma_1(N)}  \frac{1}{q_{\gamma}(n,-m)^{s}} \hspace{1cm} (s \in \CCC, \, \, \re(s)>1).
\end{align*}
We will show that $\zeta_{u,N}(s,l)$ is well-defined and has a meromorphic continuation to the whole $s$-plane. Furthermore, we will determine its residue at $s=1$.
\newline
Throughout this section we keep  the above assumptions on $N$, $l$, and $u$.
\begin{absatz}
 Let $q=\lbrack aN, bN, c \rbrack \in Q_l(N)$ with $r:= {\mathrm{gcd}}(aN, bN, c)$ the greatest common divisor of the coefficients of $q$. We define the set
$$
M^q_u(N):= \big\{(m,n) \in \ZZZ^2 \, | \, (m,n) \equiv (0,u) \!\! \! \mod \! N; \, q(n,-m) >0 \big\} \subseteq M_u(N).
$$ 
Let us set $ \theta:=\frac{bN+\sqrt{l^2-4}}{2c}$ and $\overline{\theta}:=\frac{bN-\sqrt{l^2-4}}{2c}$
 such that
 \begin{align} \label{align:quadtheta}
 q(n,-m)= c(m-\theta n) (m - \overline{\theta} n).
 \end{align}
Further, we set  $(t_q, u_q):=(t_0, \frac{u_0}{r})$ with  $(t_{0}, u_{0})$ the smallest positive solution  of Pell's equation
$X^2 - \frac{\discr(q)}{r^2} Y^2 = 4$ such that $\frac{t_{q} - bN u_{q}}{2} \equiv 1 \mod N$.
Then, the generator $\alpha_q$ of $\Gamma_1(N)_q$ is explicitly given by (cf. \cite{zagier}, p. 63, Satz 2)
\begin{align*}
\alpha_q = 
\left( \begin{matrix}
\frac{t_{q} - bN u_{q}}{2} & - c u_{q} \\ aN u_{q} & \frac{t_{q} + bN u_{q}}{2}
\end{matrix} \right)  \in \Gamma_1(N),
\end{align*}
and the power $\varepsilon_q$ of the fundamental unit $\varepsilon_0$ of \ref{absatz4:pell}  is of the form $\varepsilon_q= \frac{t_{q} +u_{q} \sqrt{l^2-4}}{2}$.
Therefore, for $(m',n')= (m,n) \alpha_q^k  \in M_u(N)$,  $k \in \ZZZ$, we find  (cf. \cite{zagier1}, p. 70)
 \begin{align} \label{align:observe}
 m' - \theta n' = \varepsilon_{q}^k (m - \theta n) \hspace{1cm} \text{and} \hspace{1cm} m' - \overline{\theta} n' = \overline{\varepsilon}_{q}^{k} (m -  \overline{\theta} n) ,
 \end{align}
  where $\overline{\varepsilon}_{q}:=  \frac{t_{q} -u_{q} \sqrt{l^2-4}}{2}$ is the conjugate of  $ \varepsilon_{q}$ in the quadratic field $\QQQ(\sqrt{l^2-4})$.

As $\varepsilon_{q}\overline{\varepsilon}_{q}>0$, equations \eqref{align:quadtheta} and \eqref{align:observe} imply $q(n',-m')>0$; this gives a  well-defined action
\begin{align} \label{align:actionpair1}
\Gamma_1(N)_q  \times M^q_u(N) & \longrightarrow M^q_u(N) \\ \notag
(\delta, (m,n))& \mapsto (m,n) \delta.
\end{align}
With the above notation and the isomorphy $\Gamma_{1,l}(N) / \Gamma_1(N) \cong Q_l(N) / \Gamma_1(N)$ in \ref{absatz4:matquad},  we find
\begin{align*}
\zeta_{u,N}(s,l) =
 \sum_{q \in Q_l(N)  / \Gamma_1(N)} \sum_{(m,n) \in M^q_u(N) / \Gamma_1(N)_{q}} \frac{1}{q(n,-m)^{s}} . 
\end{align*}
 Finally, we set $E_{q}:= \frac{t_{q} + bN u_{q}}{2cu_{q}}$ with $(t_q, u_q)$ as above, and define the set
\begin{align*}
R^{q}_{u}(N):=\left\{(m,n) \in \ZZZ^2 \, | \, (m,n) \equiv (0,v) \text{ {\rm mod} } N; n > 0, m \geq E_{q} n \right\}.
\end{align*}
\end{absatz}
\begin{rem}
By observation \eqref{align:linquad} we may  assume that the coefficients of some representative $q=\lbrack aN,bN,c \rbrack$ of an equivalence class of $ Q_l(N) / \Gamma_1(N)$  satisfies $aN>0, bN<0,$ and $c>0$.
We will exploit this fact in the sequel.
\end{rem}
\begin{lem} \label{lem4:zeta}
 Let $N$ be an odd positive integer and  $l \in \ZZZ$  with $|l|>2$ and $l \equiv 2  \mod N$.  Further, let 
  $0< u < N$ with $(u,N)=1$ as well as $0<u' < N$ with $u' \equiv -u \mod N$. Then, we have for ${\rm{Re}}(s) > 1$
\begin{align*}
\zeta_{u,N}(s,l) = \sum_{q \in Q_l(N) / \Gamma_1(N)} \left(\sum _{ (m,n) \in R^q_u(N)} \frac{1}{q(n,-m)^s}  +  \sum _{ (m,n) \in R^q_{u'}(N)} \frac{1}{q(n,-m)^s} \right).
\end{align*}
\begin{proof}
We prove that the zeta function is well-defined in \ref{absatz:continuation}.
Let us define for $q=\lbrack aN, bN, c \rbrack \in Q_l(N)$ with $c>0$ (which we can assume) the sets
\begin{align*}
M^{q \pm}_u(N):= \left\{(m,n) \in M^q_u(N) \, | \, m- \theta n \gtrless 0 \right\}. 
\end{align*}
This allows us to write $M^q_u(N)$ as a disjoint union
$M^{q}_u(N) = M^{q+}_u(N) \mathbin{\dot{\cup}} M^{q-}_u(N)$,
which descends by observation \eqref{align:observe} to
\begin{align*}
M^{q}_u(N)/ \Gamma_1(N)_q = M^{q+}_u(N)/ \Gamma_1(N)_q \mathbin{\dot{\cup}} M^{q-}_u(N)/ \Gamma_1(N)_q.
\end{align*}
Now, we define  a map
$
\varphi^+_u: M^{q+}_u(N)/\Gamma_1(N)_q \longrightarrow R^q_u(N)
$
as well as a map
$
\varphi^-_u: M^{q-}_u(N)/\Gamma_1(N)_q \longrightarrow R^q_{u'}(N)
$
which will turn out  to be bijections. This proves then the  statement of the lemma. 
\newline
To define the map $\varphi^+_u$, let $(m,n) \in M^{q+}_u(N)$.
By  observation \eqref{align:observe} we find for $(m',n') = (m,n) \alpha_q^k$, $k \in \ZZZ$,
   the equation (cf. \cite{zagier1}, p. 70)
$$
\frac{m' - \overline{\theta} n'}{m' - \theta n'} = \left( \frac{\overline{\varepsilon}_{q}}{\varepsilon_{q}} \right)^{k} \frac{m - \overline{\theta} n}{m - \theta n}
= \left( \varepsilon_{q} \right)^{-2k} \frac{m - \overline{\theta} n}{m - \theta n}.
$$
Hence in each orbit of $M^{q+}_u(N) $ by $ \Gamma_1(N)_{q} $, there is   an element $(m',n')$ which satisfies 
$
 1<  \frac{m' - \overline{\theta} n'}{m' - \theta n'} \leq \varepsilon_{q}^2.
$
We find 
\begin{align*}
 1<  \frac{m' - \overline{\theta} n'}{m' - \theta n'} \Longleftrightarrow n' > 0 \hspace{0,5cm} \text{and} \hspace{0,5cm}
\frac{m' - \overline{\theta} n'}{m' - \theta n'} \leq \varepsilon_{q}^2 \Longleftrightarrow m' \geq \frac{\varepsilon_{q}^2\theta - \overline{\theta}}{\varepsilon_{q}^2 -1} n',
\end{align*}
 and a little computation shows $\frac{\varepsilon_{q}^2\theta - \overline{\theta}}{\varepsilon_{q}^2 -1}=E_{q}$, from which we obtain  a well-defined map 
 $
  \varphi^+_u: M^{q+}_u(N) / \Gamma_1(N)_q \longrightarrow R^q_u(N)
 $
 by mapping $(m,n) \cdot \Gamma_1(N)_q \mapsto (m',n')$. 
 Now, we define the map
$
 (\varphi^+_u)^{'} : R_u^q(N) \longrightarrow M^{q+}_u(N) / \Gamma_1(N)_q 
$
   by  $(m,n) \mapsto (m,n) \cdot \Gamma_1(N)_q$. It is straightforward that this map is well-defined.
  One easily verifies that the maps $ \varphi^+_u$ and $(\varphi^+_u)^{'}$ are inverse to each other, which proves that $\varphi^+_u$ is a bijection.
 \newline
To define the map $\varphi^-_u$, we let  
 $
 \varphi : M^{q-}_u(N) / \Gamma_1(N)_q  \xrightarrow{\,\,\, \sim \,\,\,}  M^{q+}_{u'}(N) / \Gamma_1(N)_q
 $
  be the bijection induced by mapping $(m,n) \mapsto (-m,-n)$. Defining
  the map $\varphi^-_u:= \varphi^+_{u'} \circ \varphi$ gives rise to a bijection 
 $
  \varphi^-_u: M^{q-}_u(N) / \Gamma_1(N)_q \longrightarrow R^q_{u'}(N).
 $
 Since  $\varphi^+_u$ and $ \varphi^-_u$ are bijections, the statement of the lemma follows.
\end{proof}
\end{lem}

\begin{absatz}
Let $N$ be an odd positive integer,  $l \in \ZZZ$  with $|l|>2$ and $l \equiv 2  \mod N$, $0< u < N$ with $(u,N)=1$, and $q \in Q_l(N)$. Then, we define
 the theta series $\vartheta^{q}_{u,N}(t)$   by
\begin{align} \label{align:ts}
\vartheta^{q}_{u,N}(t) := \sum _{(m,n) \in R^q_u(N)} \exp \big(-t q(n,-m) \big) \hspace{1cm}(t \in \RRR_{>0}),
\end{align}
which is a little variant of the theta series studied by E. Landau in \cite{landau}  defined as follows:
Let $E \in \RRR \setminus \{0\}$  and  $L_E$ be the lattice defined by
 \begin{align*}
 L_E:= \ZZZ \omega_1  + \ZZZ \omega_2 
 \end{align*}
 with $\omega_1=(1,0)$ and $\omega_2:=(E,1)$.
 Let  $S_{E,P}$ be the truncated and shifted lattice
defined by
\begin{align*}
S_{E,P}:= \left\{(x,y) \in \RRR^2 \, | \, (x,y) \in P + L_E; y >0, x\geq Ey \right\}.
\end{align*}
and  $P= (x_0, y_0) \in \RRR^2$  lying inside the parallelogram determined by the four points $(0,0), (1,0), (E,1)$, and $(E+1, 1)$.
 Let  $q= \lbrack a,b,c \rbrack \in Q_l(1)$ be  a quadratic form with $a>0,b>0,c>0$, and discriminant $D = l^2-4$. Then, E. Landau considers the  
 theta series 
\begin{align*} 
\vartheta^{q}_{E,P}(t) := \sum _{(x,y) \in S_{E,P}} \exp \big( \! - t \, q(x,y) \big) \hspace{1cm}(t \in \RRR_{> 0}).
\end{align*}
 If  $E$ satisfies $E> \frac{-b+ \sqrt{D}}{2a}$,
 then the theta series $\vartheta^{q}_{E,P}(t)$
 converges for $t>0$, and    we have 
\begin{align} \label{align:scal}
 \vartheta^{q}_{E,P}(t) = 
\frac{\alpha_{-2}}{t } +  \frac{\alpha_{-1}}{t^{\frac{1}{2}}}  +  \alpha_{0} +  \alpha_{1} t^{\frac{1}{2}}  + \ldots +  \alpha_{k}t^{\frac{k}{2}} + O_k \left( t^{\frac{k+1}{2}} \right) \hspace{0,7cm} (\alpha_j \in \RRR, k \ge -2)
 \end{align}
 as $t \rightarrow 0$,
where all appearing constants depend on the choice of $E,P$, and $q$. Furthermore, $\alpha_{-2}$ is given by
\begin{align} \label{align:res}
\alpha_{-2} =   \frac{1}{2 \sqrt{D}} \log{\frac{2aE+b+\sqrt{D}}{2aE+b-\sqrt{D}}}.
\end{align}
(see \cite{landau}, Hilfssatz 11; note that in \cite{landau} the quadratic forms are of the form $q(X,Y)=aX^2+2bXY+cY^2$ 
and  the discriminant is defined by $\discr(q) =b^2-ac$. 
This causes the factor of $2$ appearing in the expression of $\alpha_{-2}$). 

\end{absatz}

\begin{lem} \label{lem4:landau}
Let $N$ be an odd positive integer,  $l \in \ZZZ$  with $|l|>2$ and $l \equiv 2  \mod N$, $0< u < N$ with $(u,N)=1$,
 and $q = \lbrack aN,bN,c \rbrack \in Q_l(N)$ a  quadratic form with   $aN>0$,
$bN<0, c>0$. 
Then, the theta series
$\vartheta^{q}_{u,N}(t)$  converges for $t>0$,   and we have 
\begin{align} \label{align:thetaexp}
\vartheta^{q}_{u,N}(t)=  \frac{\beta_{-2}}{t } +  \frac{\beta_{-1}}{t^{\frac{1}{2}}}  +  \beta_{0} +  \beta_{1} t^{\frac{1}{2}} +   \ldots +  \beta_{k}t^{\frac{k}{2}} + O_k \left( t^{\frac{k+1}{2}} \right) \hspace{1cm} (\beta_j \in \RRR, k \ge -2) 
\end{align}  
as $t \rightarrow 0$, where all appearing constants depend on the choice of $u,N,$ and $q$. Furthermore, $\beta_{-2}$ is  given by 
\begin{align}
\beta_{-2}= \frac{1}{N^2 \sqrt{l^2-4}} \log(\varepsilon_{q}).
\end{align}
\begin{proof}
We  show that $\vartheta^{q}_{u,N}(t)$ is a sum of Landau's theta series $\vartheta^{q}_{E,P}(t)$.  Let be $L:= cu_{q} N$.
First note that the points $(m,n) \equiv (0, u) \mod N$  with $m \geq E_q n$ and $n>0$ do lie in the interior of the parallelograms
 of width and height  $L$ of the cone defined by $y > 0$  and $x \geq E_q y$. 
Now,  let be $\xi:= \frac{x}{L}$ and $\eta:= \frac{y}{L}$. Then, we  obtain $c^2u_{q}^2$ points $P_h=(x_h, y_h)$, $h=1, \ldots , c^2u_{q}^2$, in the parallelogram now of  width and height $1$
 with respect to the coordinates $\xi$ and $\eta$.
 We set
  $$
  A:=c L^2, \hspace{0,5cm} B: = -bNL^2, \hspace{0,5cm} C:= aN L^2
  $$
  and define $q':= \lbrack A,B,C \rbrack$ such that 
$
q(y,-x) = q'(\xi, \eta)
$
and
$
 \discr(q') = B^2 - 4AC= \discr(q)L^4.
 $
 By the definition of the theta series,  we  have 
 \begin{align} \label{align:thetas}
 \vartheta^{q}_{u,N}(t)= \sum_{h=1}^{c^2u_{q}^2} \vartheta^{q'}_{E_q, P_h}(t).
 \end{align}
  Since
 $$E_{q}= \frac{t_{q} + bNu_{q}}{2c u_{q}}> \frac{bN + \sqrt{l^2-4}}{2c}= \frac{-B + \sqrt{\discr(q')}}{2A},$$
we can apply Landau's result and obtain
  \begin{align*}
 \left | \sum _{h=1}^{c^2u_{q}^2}  \left ( \vartheta^{q'}_{E_q, P_h}(t) - \frac{\alpha_{-2,h}}{t }-  \frac{\alpha_{-1,h}}{t^{\frac{1}{2}}}  -   \alpha_{0,h}  - \alpha_{1,h} t^{\frac{1}{2}} -   \ldots
  -  \alpha_{n,h}t^{\frac{k}{2}} \! \right)  \right |   <  C t^{\frac{k+1}{2}}
\end{align*}
with $\displaystyle C := \sum_{h=1}^{c^2 u_{q}^2} C_h$, where the $C_h$ are the implied
  constants of formula \eqref{align:scal}.
This implies with $\beta_j := \displaystyle  \sum_{h=1}^{c^2 u_{q}^2} \alpha_{j,h}$ the claimed growth behaviour
 \eqref{align:thetaexp} of $\vartheta^q_{u,N}(t)$ as $t \rightarrow 0$. 
 \newline
 It remains to determine $\beta_{-2}$. By   formula \eqref{align:res} we have
\begin{align*}
\beta_{-2} =& c^2u_{q}^2 \frac{1}{2 \sqrt{\discr(q')}} \log{\frac{2AE_{q}+B+\sqrt{\discr(q')}}{2AE_q+B-\sqrt{\discr(q')}}} \\ =  & c^2u_{q}^2 \frac{1}{2 L^2\sqrt{l^2-4}} \log{\frac{2cE_q-bN+\sqrt{l^2-4}}{2cE_{q}-bN-\sqrt{l^2-4}}} 
  =\frac{1}{2 N^2 \sqrt{l^2-4}} \log{\frac{t_{q}+u_{q}\sqrt{l^2-4}}{t_{q}- u_{q}\sqrt{l^2-4}}} \\ = &\frac{1}{ N^2 \sqrt{l^2-4}} \log{\frac{t_{q}+u_{q}\sqrt{l^2-4}}{2}} = \frac{1}{ N^2 \sqrt{l^2-4}} \log(\varepsilon_{q}).
\end{align*}
 This finishes the proof of the lemma.
\end{proof}
\end{lem}
\begin{absatz}{{\it{Proof of proposition \ref{prop4:continuation}}}. \label{absatz:continuation}}
Since $q(n,-m) > 0$ holds for $(m,n) \in R^q_u(N)$,  we can use the well-known  integral representation  
\begin{equation} \label{equation:theta}
\sum _{ (m,n) \in R^q_u(N)}\frac{1}{q(n,-m)^s}= \frac{1}{\Gamma(s)} \int_0^{+ \infty} t^{s-1} \vartheta^{q}_{u,N}(t) \dif t \hspace{1cm} ({\rm{Re}}(s) > 1).
 \end{equation}
Using the asymptotics \eqref{align:thetaexp}, we find
\begin{align}  \label{align:split}
\frac{1}{\Gamma(s)} \int_0^{+ \infty} t^{s-1} \vartheta^{q}_{u,N}(t) \dif t  = &
 \frac{1}{\Gamma(s)} \int_0^{1} t^{s-1} \left( \frac{\beta_{-2}}{t } +  \frac{\beta_{-1}}{t^{\frac{1}{2}}}  +  \beta_{0} +   \ldots +  \beta_{k}t^{\frac{k}{2}} \right) \dif t + \notag \\ 
&\frac{1}{\Gamma(s)} \int_0^{1} t^{s-1} R_k(t) \dif t + \frac{1}{\Gamma(s)} \int_1^{+ \infty}  t^{s-1} \vartheta^{q}_{u,N}(t) \dif t
\end{align} 
with $R_k(t) = O_k(t^{\frac{k+1}{2}})$. 
 For $\re(s)>1$ the first term on the right hand side is well-defined and is equal to
\begin{align} \label{align:pole}
\frac{1}{\Gamma(s)} \left( \frac{\beta_{-2}}{s-1} + \frac{\beta_{-1}}{s-1/2} + \frac{\beta_0}{s} + \frac{\beta_1}{s+1/2}  + \ldots + \frac{\beta_k}{s+k/2} \right),
\end{align}
which extends to a meromorphic function on the whole complex plane.
 For  ${\rm{Re}}(s) > \frac{-k+1}{2}$ with $k \geq 2$ 
 the integral of second term on the right hand side converges absolutely and uniformly, since in this case we have
 \begin{align*}
 t^{s-1} R_k(t) = O(t) \hspace{1cm} (t  \rightarrow 0).
 \end{align*}
 The third integral converges absolutely and uniformly for any $s \in \CCC$. 
  So we can deduce from equations \eqref{equation:theta} and \eqref{align:split}, applied to $u$ and $u'$,  that 
$\zeta_{u,N}(s,l)$ is well-defined for $\re(s)>1$ and  has a meromorphic continuation to the half plane 
${\rm{Re}}(s) > \frac{-k+1}{2}$ with $k \geq 2$. This finishes the proof  of the first part of the proposition as $k$ can be arbitrarily large. 
  \newline
 It remains to calculate the residue of $\zeta_{u,N}(s,l)$ at $s=1$. 
  We deduce from  equation \eqref{align:pole} that $\zeta_{u,N}(s,l)$ 
  has a simple pole
at $s=1$ with residue  
\begin{align*}
\res_{s=1}\zeta_{u,N}(s,l) = 2 h_l(N) \beta_{-2}=\frac{2h_l(N)}{N^2 \sqrt{l^2-4}} \log(\varepsilon_{q})
\end{align*}
as $\beta_{-2}$ does not depend on $u$ and $u'$. 
 This proves the proposition \ref{prop4:continuation}.
\end{absatz}

\bibliography{ref}

\def\cprime{$'$}
\providecommand{\bysame}{\leavevmode\hbox to3em{\hrulefill}\thinspace}
\providecommand{\MR}{\relax\ifhmode\unskip\space\fi MR }
\providecommand{\MRhref}[2]{%
  \href{http://www.ams.org/mathscinet-getitem?mr=#1}{#2}
}
\providecommand{\href}[2]{#2}
\begin{thebibliography}{BMMB90}

\bibitem[AL78]{atkinli}
A.~O.~L. Atkin and W.-C.~W. Li, \emph{Twists with newforms and
  pseudo-eigenvalues of {$W$}-operators}, Invent. Math. \textbf{48} (1978),
  221--243.

\bibitem[Ara74]{arakelov}
S.~J. Arakelov, \emph{Intersection theory of divisors on an arithmetic
  surface}, Math. USSR Izvestija \textbf{8} (1974), no.~6, 1167--1180.

\bibitem[AU97]{abbesullmo}
A.~Abbes and E.~Ullmo, \emph{Auto-intersection du dualisant relatif des courbes
  modulaires {$X_0(N)$}}, J. {R}eine {A}ngew. {M}ath. \textbf{484} (1997),
  1--70.

\bibitem[BMMB90]{bostmestremoret-bailly}
J.-B. Bost, J.-F. Mestre, and L.~Moret-Bailly, \emph{Sur le calcul explicite
  des ``classes de {C}hern'' des surfaces arithm{\'e}tiques de genre {$2$}},
  Ast{\'e}risque (1990), no.~183, 69--105, S{\'e}minaire sur les Pinceaux de
  Courbes Elliptiques (Paris, 1988).

\bibitem[CK09]{curillakuehn}
C.~Curilla and U.~Kuehn, \emph{On the arithmetic self-intersection numbers of
  the dualizing sheaf for fermat curves of prime exponent}, arXiv:0906.3891v1
  (2009).

\bibitem[Dri73]{drinfeld}
V.~G. Drinfeld, \emph{Two theorems on modular curves}, Funct. Anal. Appl.
  \textbf{7} (1973), no.~3, 155--156.

\bibitem[DS05]{diamondshurman}
F.~Diamond and J.~Shurman, \emph{A first course in modular forms}, Graduate
  Texts in Mathematics, vol. 228, Springer-Verlag, New York, 2005. \MR{2112196
  (2006f:11045)}

\bibitem[EC11]{edixhoven}
B.~Edixhoven and J.-M. Couveignes (eds.), \emph{Computational aspects of
  modular forms and {G}alois representations}, Annals of Mathematics Studies,
  vol. 176, Princeton University Press, Princeton, NJ, 2011, How one can
  compute in polynomial time the value of Ramanujan's tau at a prime.

\bibitem[Fal84]{faltings}
G.~Faltings, \emph{Calculus on arithmetic surfaces}, Ann. of Math. (2)
  \textbf{119} (1984), no.~2, 387--424.

\bibitem[Gup00]{gupta}
S.~D. Gupta, \emph{The {R}ankin-{S}elberg method on congruence subgroups},
  Illinois J. Math. \textbf{44} (2000), no.~1, 95--103.

\bibitem[Hej76]{hejhal0}
D.~A. Hejhal, \emph{The {S}elberg trace formula for {${\rm PSL}(2,R)$}. {V}ol.
  {I}}, Lecture Notes in Mathematics, Vol. 548, Springer-Verlag, Berlin, 1976.

\bibitem[Hej83]{hejhal}
\bysame, \emph{The {S}elberg trace formula for {${\rm PSL}(2,\,{\bf R})$}.
  {V}ol. 2}, Lecture Notes in Mathematics, vol. 1001, Springer-Verlag, Berlin,
  1983.

\bibitem[Iwa02]{iwaniec}
H.~Iwaniec, \emph{Spectral methods of automorphic forms}, second ed., Graduate
  Studies in Mathematics, vol.~53, American Mathematical Society, Providence,
  RI, 2002.

\bibitem[JK01]{jorgensonkramer}
J.~Jorgenson and J.~Kramer, \emph{Bounds for special values of {S}elberg zeta
  functions of {R}iemann surfaces}, J. Reine Angew. Math. \textbf{541} (2001),
  1--28.

\bibitem[JK06]{jorgensonkramer2}
\bysame, \emph{Bounds on canonical {G}reen's functions}, Compos. Math.
  \textbf{142} (2006), no.~3, 679--700.

\bibitem[JK09]{jorgensonkramer1}
\bysame, \emph{Bounds on {F}altings's delta function through covers}, Ann. of
  Math. (2) \textbf{170} (2009), no.~1, 1--43.

\bibitem[JZ87]{jacquetzagier}
H.~Jacquet and D.~Zagier, \emph{Eisenstein series and the {S}elberg trace
  formula. {II}}, Trans. Amer. Math. Soc. \textbf{300} (1987), no.~1, 1--48.

\bibitem[Kei06]{keil}
C.~Keil, \emph{Die {S}treumatrix f{\"u}r {U}ntergruppen der {M}odulgruppe},
  Ph.D. thesis, Universit{\"a}t Frankfurt am Main, 2006.

\bibitem[KM85]{katzmazur}
N.~M. Katz and B.~Mazur, \emph{{A}rithmetic {M}oduli of {E}lliptic {C}urves},
  Princeton University Press, Princeton, NJ, 1985.

\bibitem[Lan03]{landau}
E.~Landau, \emph{{N}euer {B}eweis eines analytischen {S}atzes des {H}errn de la
  {V}all{\'e}e {P}oussin}, Math. Ann. \textbf{56} (1903), no.~4, 645--670.

\bibitem[May12]{mayer}
H.~Mayer, \emph{{S}elf-{I}ntersection of the {R}elative {D}ualizing {S}heaf of
  {M}odular {C}urves {$ X_1(N)$}}, Ph.D. thesis, Humboldt-Universit{\"a}t zu
  Berlin, 2012.

\bibitem[MB89]{moret-bailly}
L.~Moret-Bailly, \emph{La formule de {N}oether pour les surfaces
  arithm{\'e}tiques}, Invent. Math. \textbf{98} (1989), 491--498.

\bibitem[McK72]{mckean}
H.~P. McKean, \emph{{S}elberg's trace formula as applied to a compact {R}iemann
  surface}, Comm. Pure Appl. Math. \textbf{25} (1972), 225--246.

\bibitem[MU98]{michelullmo}
P.~Michel and E.~Ullmo, \emph{Points de petite hauteur sur les courbes
  modulaires {$X_0(N)$}}, Invent. Math. \textbf{131} (1998), 645--674.

\bibitem[Ogg73]{ogg}
A.~P. Ogg, \emph{Rational points on certain elliptic modular curves}, Analytic
  number theory (Proc. Sympos. Pure Math., Vol XXIV, St. Louis Univ., St.
  Louis, Mo., 1972), Amer. Math. Soc., Providence, RI, 1973, pp.~221--231.

\bibitem[Roe66]{roelcke}
W.~Roelcke, \emph{Das {E}igenwertproblem der automorphen {F}ormen in der
  hyperbolischen {E}bene {I}}, Math. Ann. \textbf{167} (1966), 292--337.

\bibitem[Roe67]{roelcke2}
\bysame, \emph{Das {E}igenwertproblem der automorphen {F}ormen in der
  hyperbolischen {E}bene {II}}, Math. Ann. \textbf{168} (1967), 261--324.

\bibitem[Roh97]{rohrlich}
D.~E. Rohrlich, \emph{Modular curves, {H}ecke correspondence, and
  {$L$}-functions}, Modular forms and {F}ermat's last theorem ({B}oston, {MA},
  1995), Springer, New York, 1997, pp.~41--100.

\bibitem[Sou89]{soule}
C.~Soul{\'e}, \emph{G{\'e}om{\'e}trie d'{A}rakelov des surfaces
  arithm{\'e}tiques}, Ast{\'e}risque (1989), no.~177-178, Exp.\ No.\ 713,
  327--343, S{\'e}minaire Bourbaki, Vol. 1988/89.

\bibitem[Szp90]{szpiro}
L.~Szpiro, \emph{Sur les propri{\'e}t{\'e}s num{\'e}riques du dualisant relatif
  d'une surface arithm{\'e}tique}, The {G}rothendieck {F}estschrift, {V}ol.\
  {III}, Progr. Math., vol.~88, Birkh{\"a}user Boston, Boston, MA, 1990,
  pp.~229--246.

\bibitem[Zag77]{zagier}
D.~Zagier, \emph{Modular forms whose {F}ourier coefficients involve
  zeta-functions of quadratic fields}, Modular functions of one variable, {VI}
  ({B}onn, 1976), Springer-Verlag, Berlin, 1977, pp.~105--169.

\bibitem[Zag81a]{zagier2}
\bysame, \emph{Eisenstein series and the {S}elberg trace formula. {I}},
  Automorphic forms, representation theory and arithmetic ({B}ombay, 1979),
  Tata Inst. Fund. Res. Studies in Math., vol.~10, Tata Inst. Fundamental Res.,
  Bombay, 1981, pp.~303--355.

\bibitem[Zag81b]{zagier1}
\bysame, \emph{Zetafunktionen und quadratische {Z}ahlk{\"o}rper},
  Springer-Verlag, Berlin-Heidelberg-New York, 1981.

\bibitem[Zha93]{zhang}
S.-W. Zhang, \emph{Admissible pairing on a curve}, Invent. Math. \textbf{112}
  (1993), 171--193.

\end{thebibliography}

\vspace{1cm}

\noindent
Hartwig Mayer \\
Fakult\"at f\"ur Mathematik \\
Universit\"at Regensburg \\
Universit\"atsstrasse 31 \\
93053 Regensburg \\
{\it{E-mail:}} hartwig.mayer@mathematik.uni-regensburg.de

\end{document}